\documentclass[11pt]{article}
\usepackage[T1]{fontenc}
\usepackage{lmodern}
\usepackage{amsmath,amssymb,amsthm}
\usepackage{esint}
\usepackage{graphicx}
\usepackage{float}
\usepackage[margin=1in]{geometry}
\usepackage[hidelinks]{hyperref}
\numberwithin{equation}{section}
\newtheorem{theorem}{Theorem}[section]
\newtheorem{proposition}[theorem]{Proposition}
\newtheorem{lemma}[theorem]{Lemma}
\newtheorem{corollary}[theorem]{Corollary}
\theoremstyle{definition}
\newtheorem{definition}[theorem]{Definition}
\theoremstyle{remark}

\newcommand{\R}{\mathbb R}

\newcommand{\av}[1]{\overline{#1}}

\title{Sharp condition for Superharmonicity of positive singular solutions
to the critical polyharmonic equation}
\author{Xia Huang\thanks{School of Mathematical Sciences, Key Laboratory of MEA (Ministry of Education) \& Shanghai Key Laboratory of PMMP, East China Normal University. Email: \texttt{xhuang@cpde.ecnu.edu.cn}.}
\and Yuan Li\thanks{School of Mathematics and Statistics, Yunnan University. Email: \texttt{liyuan5397@163.com}.}
\and Yixing Pu\thanks{School of Mathematics and Statistics, and Key Laboratory of Nonlinear Analysis \& Applications (Ministry of Education), Central China Normal University. Email: \texttt{yxpu@stu.ecnu.edu.cn}.}}
\date{}
\begin{document}
\maketitle
\begin{abstract}
We determine the sharp condition for local superharmonicity of positive singular solutions to
the critical polyharmonic equation \((-\Delta)^m u=u^{\frac{n+2m}{n-2m}}\), where \(m\geq2\) is an integer and
\(n>2m\). For a solution in a punctured ball with a nonremovable singularity at zero, local superharmonicity
holds if and only if \(\liminf_{x\to0}u(x)>0\). The threshold is sharp. In the biharmonic
case \(m=2\), for every \(n\geq5\), we construct smooth positive solutions of \(\Delta^2u=u^{\frac{n+4}{n-4}}\)
in a punctured ball with a nonremovable singularity at \(0\), satisfying \(\liminf_{x\to0}u(x)=0\) and
\(\limsup_{x\to0}u(x)=+\infty\), while \(-\Delta u\) converges to a finite negative limit along a fixed ray. The
construction relies on an infinite array of concentrating bubbles arranged on concentric polygons,
with a Lyapunov--Schmidt reduction controlling their interactions against a vanishing
background.
\end{abstract}



\section{Introduction and main results}\label{sec:introduction}
Isolated singularities of conformally invariant elliptic equations are governed by a close interaction
between positivity, blow-up rates, and symmetry. For the critical Yamabe equation
\[
 -\Delta u=u^{\frac{n+2}{n-2}}\quad\text{in }\mathbf B_1\setminus\{0\},\,\ n\geq3,
\]
a positive solution defines, up to normalization, a conformally flat metric of constant positive scalar
curvature. Here \(\mathbf B_r\) denotes the ball of radius \(r\) centered at zero. Caffarelli, Gidas and Spruck \cite{CGS}
proved that a positive solution with a nonremovable singularity at zero satisfies
\[
 C^{-1}|x|^{-\frac{n-2}{2}}\leq u(x)\leq C|x|^{-\frac{n-2}{2}},
 \,\ u(x)=\av u(|x|)\bigl(1+O(|x|)\bigr),
\]
where \(\av u(r)\) is the spherical average of \(u\) on \(\partial\mathbf B_r\). More precisely, the solution is asymptotic, with
relative error \(o(1)\), to a radial singular solution. These radial profiles, usually called Fowler solutions,
have a periodic, possibly constant, rescaling in the logarithmic radial variable. In this second-order
setting, superharmonicity follows directly from the equation and is built into the singularity
analysis.

For \(m\geq2\), however, positivity of \(u\) does not directly determine the signs of the lower-order
iterates \((-\Delta)^j u\), \(1\leq j<m\). The local potential representation contains a smooth polyharmonic
remainder whose lower-order iterates have no prescribed sign. We consider
\begin{equation}\label{eq:original}
 (-\Delta)^m u=u^{\frac{n+2m}{n-2m}}\quad\text{in }\mathbf B_1\setminus\{0\},
 \,\ u\in C^{2m}(\mathbf B_1\setminus\{0\}),\,\ u>0,
\end{equation}
where \(m\geq2\) is an integer and \(n>2m\). Jin and Xiong \cite{JX} extended the local singularity theory
to this equation under the additional assumptions
\[
 (-\Delta)^j u\geq0\quad\text{in }\mathbf B_1\setminus\{0\},\,\ j=1,\ldots,m-1.
\]
For a nonremovable singularity, they obtained two-sided bounds of order \(|x|^{-\frac{n-2m}{2}}\) and the
asymptotic symmetry \(u(x)=\av u(|x|)(1+O(|x|))\). The lower-order sign assumptions play an essential
role in their analysis: locally, they yield an integral representation with a positive Riesz potential
and a positive smooth remainder, which permits the use of integral-equation methods for blow-up
estimates and symmetry.

This raises a natural question: does positivity of a solution of \eqref{eq:original}, together with nonremovability
of its singularity, already imply
\begin{equation}\label{eq:target}
 (-\Delta)^j u>0\quad\text{in }\mathbf B_{r_0}\setminus\{0\},\,\ j=1,\ldots,m-1,
\end{equation}
for some \(r_0\in(0,1)\)? We call \eqref{eq:target} local super-polyharmonicity; when \(m=2\), it is local superharmonicity.
The corresponding nonnegative conditions imply strict positivity by successive
applications of the strong maximum principle. Throughout the paper, a singularity is removable if
the solution extends to a \(C^{2m}\) solution across zero. No boundary condition is imposed on \(\partial\mathbf B_1\).

Several global results support an affirmative answer, but rely on the equation being defined on an
unbounded domain. Wei and Xu \cite{WX} classified positive entire solutions of the critical polyharmonic
equation as standard bubbles. More generally, Chen and Li \cite{CL} proved that positive classical
solutions of \((-\Delta)^m u\geq u^q\) in \(\mathbb R^n\), with \(q>1\), satisfy \((-\Delta)^j u>0\) for \(j=1,\ldots,m-1\). These
signs also hold for positive solutions of the critical equation on \(\mathbb R^n\setminus\{0\}\); see Ng\^o and Ye \cite{NY}. In
the biharmonic case, Lin \cite{Lin} proved that such solutions with a nonremovable singularity at zero
are radial about zero. Guo, Huang, Wang and Wei \cite{GHWW} studied their radial profiles and constructed
Delaunay-type singular solutions, while Frank and K\"onig \cite{FK} obtained a complete classification.
Recently, Huang, Li and Zhou \cite{HLZ} established further symmetry results on \(\mathbb R^n\setminus Z\), for singular
sets contained in a hyperplane and satisfying suitable capacity or dimension conditions, without
assuming \(-\Delta u\geq0\). These results collectively suggest that superharmonicity may be a universal
feature of nonremovable singularities. However, these results concern the whole space \(\mathbb R^n\) or the
complement of a prescribed singular set, and do not settle the local problem under positivity and
nonremovability alone.

In the local theory, doubling and rescaling arguments connect Liouville theorems with singularity
and decay estimates \cite{PQS}. For the strictly subcritical equation
\[
 (-\Delta)^m u=u^p\quad\text{in }\mathbf B_1\setminus\{0\},
 \,\ \frac{n}{n-2m}<p<\frac{n+2m}{n-2m},
\]
the lower-order signs are known to follow from positivity and nonremovability alone. Yang \cite{Yang}
established the isolated-singularity classification and asymptotic profile in the biharmonic case
\(m=2\) without assuming \(-\Delta u\geq0\), and Huang, Li and Yang \cite{HLY} extended this conclusion to general
\(m\): every positive classical solution with a nonremovable singularity satisfies \eqref{eq:target} near zero.

The critical exponent \(p=\frac{n+2m}{n-2m}\) is markedly harder. At this endpoint, Jin and Xiong \cite{JX}
proved that \eqref{eq:target} follows from the two-sided growth bound
\(C^{-1}|x|^{-\frac{n-2m}{2}}\leq u(x)\leq C|x|^{-\frac{n-2m}{2}}\)
near zero for some \(C>1\)---specifically, from
the lower bound \(u\geq c|x|^{-\frac{n-2m}{2}}\), or from the upper bound when the singularity is nonremovable.
Whether the lower-order signs hold without any a priori growth assumption has remained
open. The obstruction is structural: the integral representation of \(u\) splits into a Riesz potential of
\(u^p\) and a smooth polyharmonic remainder, and the sign of the latter is not prescribed.

Our first result identifies a condition on the lower limit that recovers all the lower-order signs.
It requires neither a prescribed blow-up rate nor radial symmetry.

\begin{theorem}\label{thm:positive-lower-limit}
Let \(m\geq2\), \(n>2m\), and let \(u\) satisfy \eqref{eq:original} with a nonremovable singularity at zero.
If
\begin{equation}\label{eq:B}
 \liminf_{x\to0}u(x)>0,
\end{equation}
then there exists \(r_0\in(0,1)\) such that
\[
 (-\Delta)^j u>0\quad\text{in }\mathbf B_{r_0}\setminus\{0\},\,\ j=1,\ldots,m-1.
\]
\end{theorem}

Condition \eqref{eq:B} is also necessary: if \(-\Delta u\geq0\) in a punctured neighborhood of the origin, comparison
with harmonic functions on annuli forces \(\liminf_{x\to0}u(x)>0\). Thus local superharmonicity
is equivalent to \eqref{eq:B} in this class.

It is worth noting that every positive radial solution of \eqref{eq:original} with a nonremovable singularity at
the origin satisfies \eqref{eq:B} and is therefore superharmonic in a punctured neighborhood of the origin
(see Corollary~\ref{cor:necessary}(i)).

The lower-limit dichotomy in Proposition~\ref{thm:dichotomy} makes this condition more explicit: at a nonremovable
singularity, the lower limit is either zero or \(+\infty\). In the latter case, all lower-order iterates tend
to \(+\infty\). Our second result realizes the zero-limit alternative, with a failure of superharmonicity,
for the critical biharmonic equation
\begin{equation}\label{eq:biharmonic-problem}
 \Delta^2u=u^{\frac{n+4}{n-4}}\quad\text{in }\mathbf B_1\setminus\{0\}.
\end{equation}

\begin{theorem}\label{thm:counterexample}
For every integer \(n\geq5\), there exists a positive function \(u\in C^\infty(\mathbf B_1\setminus\{0\})\) satisfying
\eqref{eq:biharmonic-problem}, with a nonremovable singularity at zero and
\[
 \liminf_{x\to0}u(x)=0,\,\ \limsup_{x\to0}u(x)=\infty.
\]
Moreover,
\[
 \int_{\mathbf B_1}|x|^{2-n}u(x)^{\frac{n+4}{n-4}}\,dx<\infty,
 \,\ -\infty<\lim_{s\to0^+}(-\Delta u)(se_n)<0,
\]
where \(e_n=(0,\ldots,0,1)\).
\end{theorem}

Consequently, positivity and nonremovability do not imply local superharmonicity at the critical
biharmonic exponent, in contrast with the strictly subcritical results above. The examples also show
that the lower-order signs in the higher-order asymptotic-symmetry theory cannot be discarded in
general. Indeed, in the zero-limit alternative the spherical average satisfies \(\av u(r)\leq Cr^2\), whereas
our solution has unbounded spikes approaching zero---so the radial-profile asymptotics \(u(x)=
\av u(|x|)(1+o(1))\), known for superharmonic singular solutions of \eqref{eq:biharmonic-problem}, cannot hold.

The construction draws on methods for solutions with many concentrating bubbles. For higher-order
prescribed curvature problems on the sphere, Guo and Li \cite{GL} obtained nonradial positive
solutions with arbitrarily large energy under suitable assumptions on the curvature. For critical
polyharmonic equations with suitable periodic coefficients, Guo, Peng and Yan \cite{GPY} constructed
infinitely bubbling solutions and proved local uniqueness. Our bubbles instead accumulate at an
interior singular point, where the solution must remain positive while its lower limit vanishes. We
arrange their centers at the vertices of regular polygons on concentric circles shrinking to zero,
extending the finite polygonal geometry of \cite{WY}; related polygonal and lattice configurations appear
in \cite{DMPP,LWX}.

The main difficulty is to control the correction relative to a positive background that itself
vanishes at the singularity. An exterior potential and a quadratic term determine the behavior
away from the bubble centers; along a normal ray, the quadratic term produces the zero limit of \(u\)
and the negative limit of \(-\Delta u\). The reference geometry cancels the leading background value and
radial derivative at each center while preserving positivity between centers. In dimensions \(n=5,6\),
the relevant lattice asymptotics require finite-part cancellations because the limiting sums do not
converge separately. A Lyapunov--Schmidt reduction and local Pohozaev identities \cite{GLN,PWY} then lead
to the equations for the infinitely many bubble parameters. Uniform estimates and a Schauder
fixed point argument complete the construction.

Section~\ref{sec:auxiliary-structure} proves the lower-limit dichotomy, Theorem~\ref{thm:positive-lower-limit}, and its consequences.
Sections~\ref{sec:exact-counterexample}--\ref{sec:ce-reduction} develop the counterexample construction, from the ansatz and the constrained equation to the
choice of bubble parameters. Section~\ref{sec:ce-final} verifies the singular behavior and proves Theorem~\ref{thm:counterexample}.
Appendix~\ref{app:balanced-geometry} establishes the geometric balance and positivity estimates, including the finite-part
analysis in dimensions \(n=5,6\).

\section{Lower limits and super-polyharmonicity}\label{sec:auxiliary-structure}
Let \(m\geq2\), \(n>2m\), and let \(u\) be a positive solution of \eqref{eq:original}.
We prove that at a nonremovable singularity its lower limit is either zero or \(+\infty\).
This dichotomy yields Theorem~\ref{thm:positive-lower-limit} and describes the alternative
in which local super-polyharmonicity fails. Write \(p=\frac{n+2m}{n-2m}\), let
\(\av f(r)\) denote the average of \(f\) over \(\partial\mathbf B_r\), and set
\[
 c_{n,k}=\frac{\Gamma(\frac{n}{2}-k)}{4^k\pi^{\frac{n}{2}}\Gamma(k)},\,\  1\leq k\leq m.
\]
Thus \(c_{n,k}|x|^{2k-n}\) is the fundamental solution of \((-\Delta)^k\).
For \(0<R<1\), define
\[
 J_j(R)=\int_{\mathbf B_R}|y|^{2m-n-2j}u(y)^p\,dy,\,\  j=0,1.
\]
Whether \(J_j(R)\) is finite does not depend on the choice of \(R\).

We first separate the Riesz potential of \(u^p\) from a smooth polyharmonic remainder.
The representation also applies to the lower-order iterates without assuming their signs.

\begin{lemma}\label{lem:representation}
For each \(R\in(0,1)\), \(u\in L^1(\mathbf B_R)\cap L^p(\mathbf B_R)\), the equation holds
in \(\mathbf B_R\) in the distributional sense, and, for \(0<|x|<R\),
\begin{equation}\label{eq:representation}
 u(x)=P_R(x)+h_R(x),
 \,\ P_R(x)=c_{n,m}\int_{\mathbf B_R}|x-y|^{2m-n}u(y)^p\,dy,
\end{equation}
where \(h_R\) is smooth and \((-\Delta)^m h_R=0\) in \(\mathbf B_R\).
For \(j=1,\ldots,m-1\) and \(0<|x|<R\), one also has
\begin{equation}\label{eq:iterated}
 (-\Delta)^j u(x)=c_{n,m-j}\int_{\mathbf B_R}|x-y|^{2m-n-2j}u(y)^p\,dy
                 +(-\Delta)^j h_R(x).
\end{equation}
\end{lemma}

\begin{proof}
Since \(p=\frac{n+2m}{n-2m}>\frac{n}{n-2m}\), Lemma~2.1 of \cite{JX}, after scaling to
\(\mathbf B_R\), gives \(u^p\in L^1(\mathbf B_R)\) and a Green representation without
sign assumptions on the lower iterates. Splitting its Green kernel into
\(c_{n,m}|x-y|^{2m-n}\) and a smooth correction gives \eqref{eq:representation}, with
\(h_R\) smooth and polyharmonic across zero. In particular, the distributional equation
contains no term supported at zero, and \(u\in L^1(\mathbf B_R)\) follows from
\(u\in L^p(\mathbf B_R)\). Differentiating the potential in distributions gives
\eqref{eq:iterated}. Its kernels are locally integrable because \(n-2m+2j\leq n-2\),
and the identities hold pointwise away from zero.
\end{proof}

The remainder and its iterated Laplacians are bounded near zero, so divergence of the
weighted integrals forces the corresponding potentials to tend to \(+\infty\).
When \(J_0(R)\) is finite, spherical averaging instead determines the lower limit of \(u\).

\begin{proposition}\label{thm:dichotomy}
Suppose that the singularity of \(u\) at zero is nonremovable. If \(J_0(R)=\infty\), then
\[
 \lim_{x\to0}(-\Delta)^j u(x)=+\infty,\,\  j=0,\ldots,m-1.
\]
If \(J_0(R)<\infty\), then
\begin{equation}\label{eq:zero-branch}
 \liminf_{x\to0}u(x)=0,\,\  0<\av u(r)\leq Cr^2\quad(0<r<r_0).
\end{equation}
Here \(C,r_0>0\) may depend on \(u,m,n\). In addition, \(J_1(R)=\infty\) implies
\((-\Delta)^j u(x)\to+\infty\) for every \(j=1,\ldots,m-1\).
\end{proposition}

\begin{proof}
For a fixed \(j\in\{0,\ldots,m-1\}\), suppose that
\(\int_{\mathbf B_R}|y|^{2m-n-2j}u(y)^p\,dy=\infty\). On each annulus
\(\varepsilon<|y|<\frac{R}{2}\),
\[
 \int_{\varepsilon<|y|<\frac{R}{2}}|x-y|^{2m-n-2j}u(y)^p\,dy
 \to
 \int_{\varepsilon<|y|<\frac{R}{2}}|y|^{2m-n-2j}u(y)^p\,dy
 \,\  (x\to0).
\]
The convergence follows from uniform convergence of the kernel there.
The right-hand side tends to infinity as \(\varepsilon\to0^+\).
The remaining potential is nonnegative, while \((-\Delta)^j h_R\) is bounded near zero.
Thus \eqref{eq:representation}--\eqref{eq:iterated} give \((-\Delta)^j u(x)\to+\infty\).
Since \(R<1\), divergence of \(J_0\) implies divergence of all these weighted integrals
for \(j\geq0\), and divergence of \(J_1\) does so for \(j\geq1\).
This proves both divergence assertions.

Suppose that \(J_0(R)<\infty\). For \(0<\sigma\leq n-2\), the distributional
superharmonicity of \(|x-y|^{-\sigma}\) gives
\[
 \frac1{|\mathbb S^{n-1}|}\int_{\mathbb S^{n-1}}|r\theta-y|^{-\sigma}\,d\theta
 \leq |y|^{-\sigma},\,\  y\ne0.
\]
One may apply the mean value inequality to smooth regularizations; the average tends
to \(|y|^{-\sigma}\) as \(r\to0^+\). With \(\sigma=n-2m\), dominated convergence
therefore yields
\[
 \av{P_R}(r)\leq c_{n,m}J_0(R),\,\  \av{P_R}(r)\to c_{n,m}J_0(R).
\]
Fatou's lemma along any sequence \(x\to0\) gives \(\liminf P_R(x)\geq c_{n,m}J_0(R)\).
The reverse bound for the lower limit of \(u\) follows because its minimum on each sphere
is at most its average. Hence
\[
 L:=\liminf_{x\to0}u(x)=\lim_{r\to0^+}\av u(r)
   =c_{n,m}J_0(R)+h_R(0)\in[0,\infty).
\]
We claim that \(L=0\). For \(0<s<R\), the same representation gives
\[
 \begin{aligned}
 h_s(x)&=u(x)-c_{n,m}\int_{\mathbf B_s}|x-y|^{2m-n}u(y)^p\,dy,\,\  0<|x|<s,\\
 h_s(0)&=L-c_{n,m}J_0(s).
 \end{aligned}
\]
If \(L>0\), finiteness of \(J_0(R)\) allows us to choose \(s\) so that \(h_s(0)>0\),
and then \(r<s\) so that \(h_s>0\) on \(\mathbf B_r\). For \(x\in\mathbf B_r\),
\[
 h_r(x)=h_s(x)+c_{n,m}\int_{\mathbf B_s\setminus\mathbf B_r}|x-y|^{2m-n}u(y)^p\,dy>0.
\]
Set \(a=c_{n,m}^{\frac{1}{p-1}}r^{\frac{n-2m}{2}}\), \(v(x)=au(rx)\), and
\(\widetilde h(x)=ah_r(rx)\). A change of variables gives
\[
 v(x)=\int_{\mathbf B_1}|x-y|^{2m-n}v(y)^p\,dy+\widetilde h(x),\,\  0<|x|<1.
\]
Here \(v\in L^p(\mathbf B_1)\), while \(\widetilde h\) is positive, smooth and polyharmonic
in \(\mathbf B_1\). The integer-order case of Theorem~1.3 of \cite{JX} therefore applies.
Since the singularity is nonremovable, it yields \(u(x)\geq c_0|x|^{-\frac{n-2m}{2}}\)
for \(0<|x|<r_0\). Consequently,
\[
 J_0(R)\geq c_0^p|\mathbb S^{n-1}|\int_0^{r_0}t^{-1-\frac{n-2m}{2}}\,dt=\infty.
\]
This is a contradiction. Thus \(L=0\). Taylor expansion of the smooth remainder, whose
linear term has zero spherical average, now gives
\[
 0<\av u(r)\leq c_{n,m}J_0(R)+\av{h_R}(r)
 =c_{n,m}J_0(R)+h_R(0)+O(r^2)=O(r^2).
\]
This proves \eqref{eq:zero-branch}.
\end{proof}

\begin{proof}[Proof of Theorem~\ref{thm:positive-lower-limit}]
Since \(\liminf_{x\to0}u(x)>0\), Proposition~\ref{thm:dichotomy} excludes the finite
\(J_0(R)\) case and gives \((-\Delta)^j u(x)\to+\infty\) for \(j=1,\ldots,m-1\).
There are finitely many such indices, so all these functions are positive in one common
punctured neighborhood of zero.
\end{proof}

We next record the consequences for radial solutions and for singularities where the
lower-order signs fail.

\begin{corollary}\label{cor:necessary}
Suppose that the singularity of \(u\) at zero is nonremovable.
\begin{enumerate}
\renewcommand{\labelenumi}{(\roman{enumi})}
\item If \(u\) is radial, then it satisfies \eqref{eq:target} and \eqref{eq:B}.
\item If \eqref{eq:target} fails, then \(u\) is nonradial and
\[
 J_1(R)<\infty,\,\  \liminf_{x\to0}u(x)=0,\,\  \limsup_{x\to0}u(x)=\infty.
\]
Its spherical average satisfies \(0<\av u(r)\leq Cr^2\) for \(0<r<r_0\), with
\(C,r_0>0\). Moreover, for every small \(s>0\),
\[
 u\notin L^{\frac{2n}{n-2m}}(\mathbf B_s),
 \,\ u(x)\ne O(|x|^{2m-n-2})\quad\text{as }x\to0.
\]
\end{enumerate}
\end{corollary}

\begin{proof}
We first give the removability criterion used in both parts. If \(u\) is bounded near
zero, Lemma~\ref{lem:representation} supplies the distributional equation across zero.
Interior \(W^{2m,q}\) estimates, Sobolev embedding and Schauder estimates then give a
nonnegative \(C^{2m}\) extension. Its value at zero is allowed to vanish.

Suppose more generally that \(u\in L^{p+1}(\mathbf B_s)\), where \(p+1=\frac{2n}{n-2m}\).
Choose \(0<R_1<R_2<s\) and set
\(V=c_{n,m}u^{p-1}\in L^{\frac{n}{2m}}(\mathbf B_{R_2})\). The representation implies
\[
 0\leq u(x)\leq\int_{\mathbf B_{R_2}}\frac{V(y)u(y)}{|x-y|^{n-2m}}\,dy
              +|h_{R_2}(x)|,\,\  0<|x|<R_1,
\]
where \(h_{R_2}\) is bounded on \(\mathbf B_{R_1}\). Apply Corollary~1.1 of \cite{Li}
with potential order \(2m\), initial exponent \(p+1>\frac{n}{n-2m}\), and target exponent
\(pq\), where \(q>\frac{n}{2m}\); this choice gives \(pq>p+1\). It follows that \(u^p\in L^q\)
on a smaller ball. H\"older's inequality in \eqref{eq:representation}, applied on that
ball, makes \(u\) bounded near zero, so the preceding argument applies. This reasoning
requires no sign condition on the polyharmonic remainder.

If \(u\) is radial and \(J_0(R)<\infty\), Proposition~\ref{thm:dichotomy} gives
\(u(r)=\av u(r)\leq Cr^2\) for small \(r\). The bounded extension just established would
make the singularity removable. Hence \(J_0(R)=\infty\), and the divergent case of
Proposition~\ref{thm:dichotomy} proves (i).

For (ii), the \(J_1\)-divergence assertion in Proposition~\ref{thm:dichotomy} shows that
\(J_1(R)<\infty\), and hence \(J_0(R)<\infty\). The finite branch gives the zero lower limit
and the estimate for \(\av u\), while (i) excludes radial symmetry.
Nonremovability and the criterion above imply both \(\limsup_{x\to0}u(x)=\infty\) and
failure of the stated integrability. Finally, if \(u(x)\leq C|x|^{2m-n-2}\) near zero, then
\[
 \int_{\mathbf B_s}u^{\frac{2n}{n-2m}}\,dy
 =\int_{\mathbf B_s}u\,u^p\,dy
 \leq C\int_{\mathbf B_s}|y|^{2m-n-2}u(y)^p\,dy=CJ_1(s)<\infty,
\]
again a contradiction. This also recovers the growth criterion in Proposition~2.3 of \cite{JX}.
\end{proof}

\section{The ansatz and its properties}\label{sec:exact-counterexample}
\subsection{Construction of the ansatz}\label{subsec:ce-ansatz}
Fix an integer \(n\geq5\) and a large even integer \(\ell\).
The geometric parameters \(\gamma,\alpha\), with \(\alpha>0\), are chosen as in
Lemma~\ref{lem:ce-certificates} and fixed before the small parameter \(b>0\).
The constant in the bubble formula is
\[
 A_n=[n(n-4)(n^2-4)]^{\frac{n-4}{8}}.
\]
After fixing \(n,\ell,\gamma,\alpha\), choose
\(1<R<\gamma^{-1}\) and work in \(\Omega=\mathbf B_R\).
The resulting solution is restricted to \(\mathbf B_1\). Write
\(x=(\bar x,x'')\in\mathbb R^2\times\mathbb R^{n-2}\). For
\(\lambda>0\) and \(\mathbf x\in\mathbb R^n\), set
\[
 \mathcal B_{\lambda,\mathbf x}(x)
 =A_n\left(\frac\lambda{\lambda^2+|x-\mathbf x|^2}\right)^{\frac{n-4}{2}},
 \,\  \Delta^2 \mathcal B_{\lambda,\mathbf x}=\mathcal B_{\lambda,\mathbf x}^{\frac{n+4}{n-4}}.
\]
Let
\[
 \bar\xi_i=\left(\cos\frac{2\pi(i-1)}{\ell},
                   \sin\frac{2\pi(i-1)}{\ell}\right),
 \,\ \xi_i=(\bar\xi_i,\mathbf0),\,\  i=1,\ldots,\ell.
\]
For \(k\geq-3\), write \(\mathbf x_{k,i}=r_k\xi_i\).
The interior circles \(k=0,1,\ldots\) each carry \(\ell\) bubbles
of width \(\lambda_k\).

For \(b>0\), write
\(\lambda=(\lambda_k)_{k\geq-3}\) and
\(r=(r_k)_{k\geq-3}\), and choose
\begin{equation}\label{eq:ce-reference-parameters}
 r_k\ \hbox{near }\gamma^k,
 \,\ \lambda_k\ \hbox{near }\lambda_k^0:=
 \left(\frac{b\gamma^{k(n-2+\alpha)}}{\ell A_n}\right)^{\frac2{n-4}},
 \,\  k\geq-3.
\end{equation}
We write \(\lambda^0=(\lambda_k^0)_{k\geq-3}\).
The parameters \((\lambda,r)\) will satisfy \eqref{eq:ce-parameter-space};
in particular, \(r_k\to0\) as \(k\to\infty\).

Let \(I\) be the identity operator, and let
\[
 J_2[f](x)=f(0)+Df(0)\cdot x+\frac12D^2f(0)[x,x].
\]
Then \(\Delta^2((I-J_2)[f])=\Delta^2f\). 
We will apply \(I-J_2\) to the circle sums in
Definitions~\ref{def:ce-interior-sum} and \ref{def:ce-exterior-sum}
below, which makes the exterior sum \(h_{r,\lambda}\)
well defined.
\par\noindent
\begin{minipage}{\linewidth}
\begin{definition}\label{def:ce-interior-sum}
The contribution of the interior circles is
\[
 \mathbf h_{r,\lambda}(x)
 =\sum_{k=0}^\infty(I-J_2)
       \left[\sum_{i=1}^{\ell}\mathcal B_{\lambda_k,\mathbf x_{k,i}}\right](x).
\]
\end{definition}
Figure~\ref{fig:ce-centers} illustrates the circle arrangement in
Definition~\ref{def:ce-interior-sum}.

\begin{figure}[H]
 \centering
 \includegraphics[width=\linewidth]{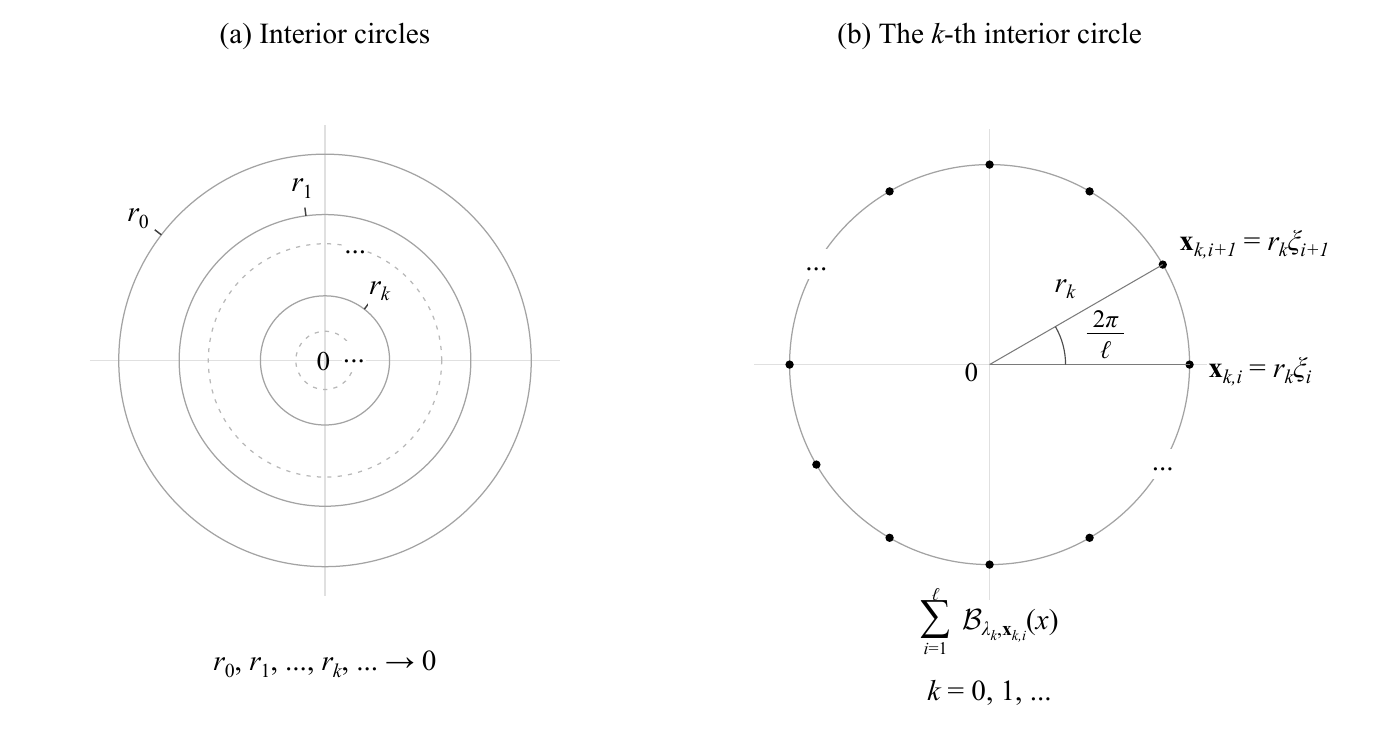}
 \caption{Schematic interior circles in the plane \(x''=0\), approaching the origin.
 The enlarged panel shows representative points
 \(\mathbf x_{k,i}=r_k\xi_i\). Each circle has \(\ell\) equally spaced
 vertices, with \(\ell\) large and even.}
 \label{fig:ce-centers}
\end{figure}
\end{minipage}
\par\smallskip

Set
\[
\mathcal K_a(x)=|x-a|^{4-n},
\,\  a\ne0,\quad x\ne a.
\]
For \(\mathbf x\in\mathbb R^n\setminus\overline\Omega\), \(\mathcal K_{\mathbf x}\) is smooth and
\(
 \Delta^2\mathcal K_{\mathbf x}=0
\) in \(\Omega\).

On the exterior circles, we use \(\mathcal K_{\mathbf x}\) in the directions
\(\xi_i\). We number them outward from 1. For \(k=1,2,3\),
we use the points
\(\mathbf x_{-k,i}=r_{-k}\xi_i\) with coefficient
\(A_n\lambda_{-k}^{\frac{n-4}{2}}\) at each point. For \(k\geq4\),
we use the points \(\gamma^{-k}\xi_i\) with coefficient
\(\frac{b}{\ell}\gamma^{-k(n-2+\alpha)}\). All exterior points lie outside \(\overline\Omega\)
for \((\lambda,r)\) satisfying \eqref{eq:ce-parameter-space}.
\par\noindent
\begin{minipage}{\linewidth}
\begin{definition}\label{def:ce-exterior-sum}
The contribution of the exterior circles is
\[
\begin{aligned}
 h_{r,\lambda}(x)={}&A_n\sum_{k=1}^{3}\lambda_{-k}^{\frac{n-4}{2}}
 (I-J_2)\left[\sum_{i=1}^{\ell}\mathcal K_{\mathbf x_{-k,i}}\right](x)\\
 &+\frac b\ell\sum_{k=4}^{\infty}\gamma^{-k(n-2+\alpha)}
 (I-J_2)\left[\sum_{i=1}^{\ell}\mathcal K_{\gamma^{-k}\xi_i}\right](x).
\end{aligned}
\]
\end{definition}
Figure~\ref{fig:ce-exterior} illustrates the circle arrangement in
Definition~\ref{def:ce-exterior-sum}.

\begin{figure}[H]
 \centering
 \includegraphics[width=\linewidth]{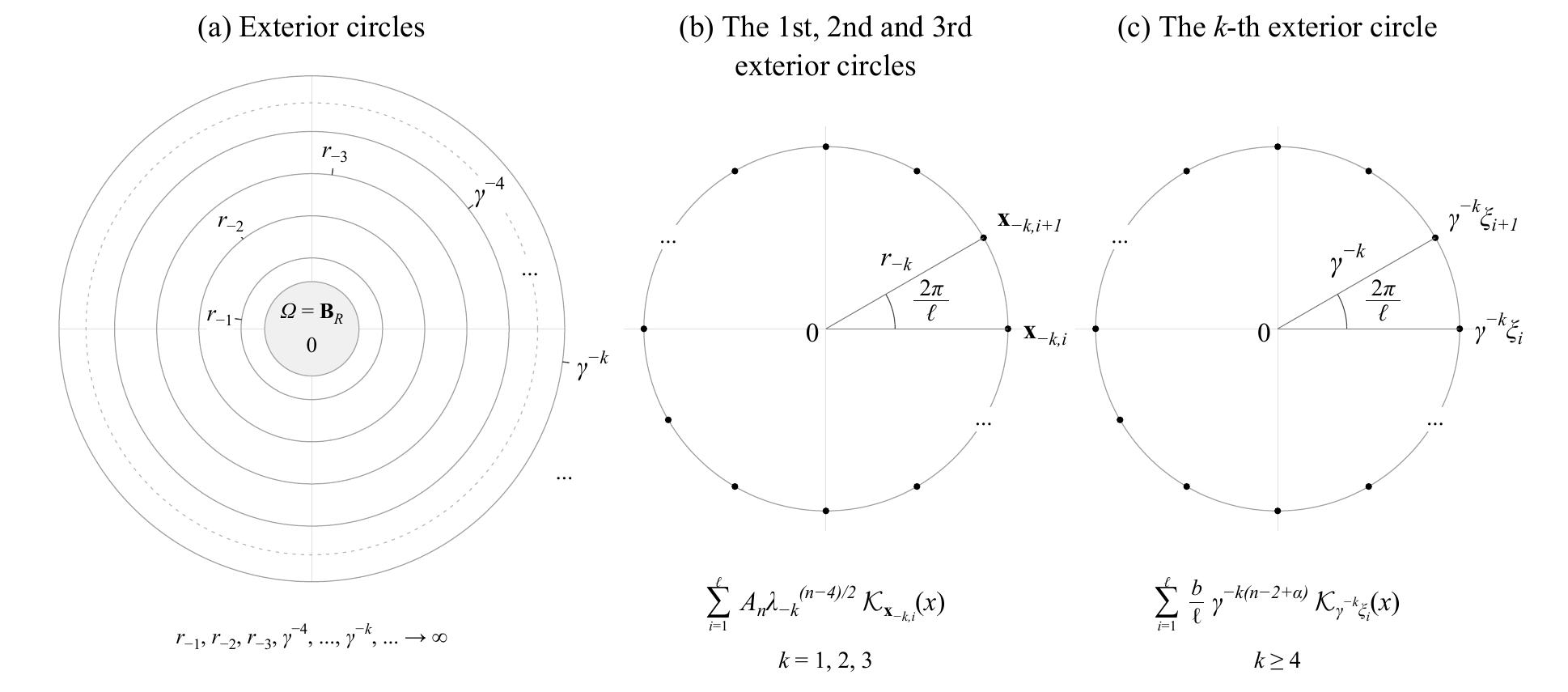}
 \caption{Schematic exterior circles in the plane \(x''=0\) and
 \(\mathcal K_{\mathbf x}(x)\), with \(\mathbf x\) at the vertices of one circle. In panel (b), \(k=1,2,3\);
 in panel (c), \(k\geq4\).
 Each circle has the same large even number \(\ell\) of
 equally spaced points; only representative points are drawn in the
 enlarged panels.}
 \label{fig:ce-exterior}
\end{figure}
\end{minipage}

We consider the ansatz
\[
 U_{\lambda,r}(x)=\mathbf h_{r,\lambda}(x) + h_{r,\lambda}(x) + b^{\frac{n-4}{n-3}}|x''|^2,
\]
and seek a solution of \eqref{eq:biharmonic-problem} in the form
\begin{equation}\label{eq:ce-solution-form}
 u=U_{\lambda,r}+\phi.
\end{equation}
The quadratic term \(b^{\frac{n-4}{n-3}}|x''|^2\) gives the lower bound
\eqref{eq:ce-profile-positive}, which combines with
\eqref{eq:ce-nonlinear-bound} to yield \(u>0\).
Combined with the estimates in the proof of Proposition~\ref{prop:ce-existence}
for \(\mathbf h_{r,\lambda},h_{r,\lambda}\) and \(\phi\), the quadratic term also gives
\(u(se_n)\to0\) as \(s\to0^+\), which leads to
\(\liminf_{x\to0}u(x)=0\).
The interior sum \(\mathbf h_{r,\lambda}\) gives
\(u(\mathbf x_{k,1})=A_n\lambda_k^{-\frac{n-4}{2}}(1+o(1))\to\infty\)
as \(k\to\infty\), which leads to \(\limsup_{x\to0}u(x)=\infty\).
The exterior sum \(h_{r,\lambda}\) is taken over infinitely many circles
so that we can choose suitable \(\gamma,\alpha\) in Section~\ref{subsec:ce-point-balance}.
We include the six exterior parameters \(\lambda_{-j},r_{-j}\), \(j=1,2,3\),
in \((\lambda,r)\) to solve \eqref{eq:ce-parameter-equations} (see Section~\ref{subsec:ce-parameter-choice}).

In Sections~\ref{sec:ce-pde} and~\ref{sec:ce-reduction},
we will show that there exist suitable \(\lambda,r\) and \(\phi\)
such that the function \(u\) given by \eqref{eq:ce-solution-form}
solves \eqref{eq:biharmonic-problem}.

\subsection{The choice of \(\gamma\) and \(\alpha\)}\label{subsec:ce-point-balance}
By symmetry, we may take \(i=1\) in our discussion. Define
\begin{equation}\label{eq:ce-background-decomposition}
 H_k(x):=U_{\lambda,r}(x)
          -\mathcal B_{\lambda_k,\mathbf x_{k,1}}(x),\,\  k\geq0.
\end{equation}
Fix \(k_0\geq0\). The leading terms in the equations for \(\lambda_{k_0},r_{k_0}\)
are determined by \(H_{k_0}(\mathbf x_{k_0,1})\) and
\(\partial_{\xi_1}H_{k_0}(\mathbf x_{k_0,1})\); see
Section~\ref{sec:ce-reduction}. We first choose \(\gamma,\alpha\)
so that \(H_{k_0}(\mathbf x_{k_0,1})\) and
\(\partial_{\xi_1}H_{k_0}(\mathbf x_{k_0,1})\) vanish at leading order.

For \(r>0\), define
\[
 \Phi_r(x)=(I-J_2)
 \left[\frac1\ell\sum_{i=1}^{\ell}\mathcal K_{r\xi_i}\right](x).
\]
By scaling, \(\Phi_r(x)=r^{4-n}\Phi_1(\frac{x}{r})\).

For \(0<\gamma<1\) and \(0<\alpha<2\), the reference parameters are
\(\lambda_k=\lambda_k^0\) and \(r_k=\gamma^k\), \(k\geq-3\).
Define
\begin{equation}\label{eq:ce-normalized-background}
\begin{aligned}
 H(s):={}&\sum_{k\in\mathbb Z\setminus\{0\}}
 \gamma^{k(n-2+\alpha)}\Phi_{\gamma^k}(s\xi_1)\\
 &+\frac1\ell\sum_{i=2}^{\ell}\mathcal K_{\xi_i}(s\xi_1)
       -1-\frac{(n-4)^2}{4}s^2.
\end{aligned}
\end{equation}
The series defining \(H\) in \eqref{eq:ce-normalized-background}
and its derivatives converge locally
uniformly for \(\gamma<s<\gamma^{-1}\), so
\(H\in C^\infty((\gamma,\gamma^{-1}))\).
The following lemma shows that \(b\gamma^{k_0(2+\alpha)}H(s)\)
is the leading term in the expansion of
\(H_{k_0}(\gamma^{k_0}s\xi_1)\) as \(b\rightarrow 0\).

\begin{lemma}\label{lem:ce-background-comparison}
For \(\lambda_k=\lambda_k^0\) and \(r_k=\gamma^k\),
\(k\geq-3\), and any fixed compact subset
\(K\subset(\gamma,R)\),
\begin{equation}\label{eq:ce-background-comparison}
\begin{split}
 H_{k_0}(\gamma^{k_0}s\xi_1)
 &=b\gamma^{k_0(2+\alpha)}H(s)
     +O\!\left(b^{\frac{n}{n-4}}\gamma^{4k_0}\right),\\
 \partial_{\xi_1}H_{k_0}(\gamma^{k_0}s\xi_1)
 &=b\gamma^{k_0(1+\alpha)}H'(s)
     +O\!\left(b^{\frac{n}{n-4}}\gamma^{3k_0}\right).
\end{split}
\end{equation}
Here \(\partial_{\xi_1}\) is the directional derivative along
\(\xi_1\). The estimates hold as \(b\rightarrow 0\), uniformly for
\(k_0\geq0\) and \(s\in K\); their constants may depend on
\(n,\ell,\gamma,\alpha,K\), but not on \(b,k_0,s\).
\end{lemma}
\begin{proof}
Let \(x=(\bar x,x'')=\gamma^{k_0}s\xi_1\). Then
\[
 H_{k_0}(\gamma^{k_0}s\xi_1)-b\gamma^{k_0(2+\alpha)}H(s)
 =\sum_{k=0}^{\infty}e_{k,k_0}(s),
\]
where
\[
 e_{k,k_0}(s):=
 \begin{cases}
  \displaystyle(I-J_2)\left[\sum_{i=1}^{\ell}d_{k,i}\right]
             (\gamma^{k_0}s\xi_1),&k\ne k_0,\\[2mm]
  \begin{aligned}
   &\sum_{i=2}^{\ell}d_{k_0,i}(\gamma^{k_0}s\xi_1) -J_2\left[\sum_{i=1}^{\ell}d_{k_0,i}\right]
             (\gamma^{k_0}s\xi_1),
  \end{aligned}&k=k_0,
 \end{cases}
\]
and, for \(k\geq0\),
\[
 d_{k,i}(x):=\mathcal B_{\lambda_k,r_k\xi_i}(x)
       -A_n\lambda_k^{\frac{n-4}{2}}\mathcal K_{r_k\xi_i}(x).
\]
For \(x\ne r_k\xi_i\), the mean value theorem gives
\begin{equation}\label{eq:ce-bubble-difference-estimate}
\begin{aligned}
 d_{k,i}(x)
 &=-\frac{n-4}{2}A_n\lambda_k^{\frac{n}{2}}
   \int_0^1\bigl(|x-r_k\xi_i|^2+\sigma\lambda_k^2\bigr)^{-\frac{n-2}{2}}
   \,d\sigma,\\
 |D^q d_{k,i}(x)|
 &\leq C\lambda_k^{\frac{n}{2}}|x-r_k\xi_i|^{2-n-q},
 \,\  0\leq q\leq4.
\end{aligned}
\end{equation}
For \(s\in K\),
\begin{equation}\label{eq:ce-circle-distance}
 |\gamma^{k_0}s\xi_1-r_k\xi_i|\geq c_K\max\{\gamma^{k_0},r_k\},
 \,\  (k,i)\ne(k_0,1).
\end{equation}
Since \(\ell\) is even,
\[
 \sum_{i=1}^{\ell}d_{k,i}(-x)=\sum_{i=1}^{\ell}d_{k,i}(x).
\]
For \(k<k_0\), Taylor expansion gives
\[
\begin{aligned}
 (I-J_2)\left[\sum_{i=1}^{\ell}d_{k,i}\right](\gamma^{k_0}s\xi_1)
 &=\frac{\gamma^{4k_0}s^4}{6}\sum_{i=1}^{\ell}\int_0^1(1-t)^3 \partial_{\xi_1}^{\,4}d_{k,i}(t\gamma^{k_0}s\xi_1)\,dt,\\
 |t\gamma^{k_0}s\xi_1-r_k\xi_i|
 &\geq(1-R\gamma)r_k,\,\  0\leq t\leq1.
\end{aligned}
\]
Together with \eqref{eq:ce-bubble-difference-estimate} and
\eqref{eq:ce-circle-distance} for \(k\geq k_0\), this gives
\[
 |e_{k,k_0}(s)|\leq C
 \begin{cases}
  \lambda_k^{\frac{n}{2}}r_k^{-n-2}\gamma^{4k_0},&k\leq k_0,\\
  \lambda_k^{\frac{n}{2}}r_k^{-n}\gamma^{2k_0},&k>k_0.
 \end{cases}
\]
Using \(\lambda_k=\lambda_k^0\) from
\eqref{eq:ce-reference-parameters}, we obtain
\[
\begin{aligned}
 \left|\sum_{k=0}^{\infty}e_{k,k_0}(s)\right|
 &\leq C\left[
  \gamma^{4k_0}\sum_{k=0}^{k_0}\lambda_k^{\frac{n}{2}}r_k^{-n-2}
  +\gamma^{2k_0}\sum_{k=k_0+1}^{\infty}\lambda_k^{\frac{n}{2}}r_k^{-n}\right]\\
 &\leq C b^{\frac{n}{n-4}}\left[
  \gamma^{4k_0}\sum_{k=0}^{k_0}\gamma^{\frac{k(8+n\alpha)}{n-4}}
  +\gamma^{2k_0}\sum_{k=k_0+1}^{\infty}
                  \gamma^{\frac{kn(2+\alpha)}{n-4}}\right]\\
 &\leq C b^{\frac{n}{n-4}}\gamma^{4k_0},
\end{aligned}
\]
since \(\frac{8+n\alpha}{n-4}>0\) and
\(2+\frac{n(2+\alpha)}{n-4}>4\). This gives the first estimate of
\eqref{eq:ce-background-comparison}. The second estimate follows by a
similar argument, differentiating with respect to \(s\) and using
\(\partial_{\xi_1}=\gamma^{-k_0}\partial_s\) along
\(x=\gamma^{k_0}s\xi_1\).
\end{proof}

At \(s=1\), \eqref{eq:ce-background-comparison} gives the leading value
and radial derivative of \(H_{k_0}\) at \(\gamma^{k_0}\xi_1\). Relative to their respective
scales \(b\gamma^{k_0(2+\alpha)}\) and
\(b\gamma^{k_0(1+\alpha)}\), both errors are
\(O(b^{\frac{4}{n-4}}\gamma^{k_0(2-\alpha)})\).
Thus we impose the two conditions
\begin{equation}\label{eq:ce-balance}
 H(1)=H'(1)=0
\end{equation}
so that \(b\gamma^{k_0(2+\alpha)}H(1)\) and
\(b\gamma^{k_0(1+\alpha)}H'(1)\) vanish. When the reference
parameters are perturbed as in \eqref{eq:ce-parameter-space}, the
additional errors in the value and radial derivative at \(\mathbf x_{k_0,1}\) are
controlled in Lemma~\ref{lem:ce-profile}.

We now choose \(\gamma,\alpha\) so that \eqref{eq:ce-balance} holds.
\begin{lemma}\label{lem:ce-certificates}
For each integer \(n\geq5\), there exist \(\ell_0(n)\) and
\(\kappa_n^*>2\pi\) such that
for every even \(\ell\geq\ell_0(n)\) one can choose
\[
 e^{-\frac{\kappa_n^*}{\ell}}\leq\gamma\leq e^{-\frac{2\pi}{\ell}},
 \,\ 0<\alpha<1+\frac3n,
\]
so that \eqref{eq:ce-balance} holds.
\end{lemma}
\begin{proof}
Lemma~\ref{lem:large-ell-root} gives the stated parameter range and
\eqref{eq:ce-balance}.
\end{proof}

\subsection{Estimates for the ansatz}
Take \(\ell_0(n)\) large enough.
For an even \(\ell\geq\ell_0(n)\), fix \(\gamma,\alpha\)
as in Lemma~\ref{lem:ce-certificates}.
For any fixed \(\rho\in(\gamma^{2-\alpha},1)\), assume that
\begin{equation}\label{eq:ce-parameter-space}
 \left|\left(\frac{\lambda_k}{\lambda_k^0}\right)^{\frac{n-4}{2}}-1\right|
 +\left|\frac{r_k}{\gamma^k}-1\right|
 \leq\varepsilon\rho^k,\,\  k\geq-3.
\end{equation}
The value of \(\rho\) used to solve the parameter equations
from \eqref{eq:ce-reduced-parameter-equation} is chosen in
Subsection~\ref{app:parameter-surjectivity}.
Fix \(\varepsilon_0>0\) small enough. The estimates in this
subsection hold uniformly for every \(0\leq\varepsilon\leq\varepsilon_0\) and all
small \(b>0\), with an upper bound for \(b\) depending on
\(\varepsilon_0\).
The numbers \(n,\ell,\gamma,\alpha,R,\rho,\varepsilon_0\) are chosen
before \(b\). The positive constants \(C\) in this subsection may
depend on these numbers,
but are independent of \(b,k\) and \((\lambda,r)\) satisfying
\eqref{eq:ce-parameter-space}; they may be increased to include lower bounds written
with \(C^{-1}\).
Equation~\eqref{eq:ce-parameter-space} gives
\[
\begin{gathered}
 C^{-1} b^{\frac{2}{n-4}}r_k^{\frac{2(n-2+\alpha)}{n-4}}
 \leq\lambda_k\leq C b^{\frac{2}{n-4}}r_k^{\frac{2(n-2+\alpha)}{n-4}},\\
 \frac{\lambda_k}{r_k}\leq C b^{\frac{2}{n-4}}r_k^{\frac{n+2\alpha}{n-4}}.
\end{gathered}
\]
To prove \(U_{\lambda,r}>0\), we first estimate \(\mathcal P_n(r)\),
the sum along the ray with polar angle \(\frac{\pi}{\ell}\) in the plane
\(x''=0\), defined for \(r>0\) by
\begin{equation}\label{eq:ce-midpoint-potential}
 \mathcal P_n(r)=\sum_{k\in\mathbb Z}\gamma^{k(n-2+\alpha)}
 \Phi_{\gamma^k}\left(r\cos\frac\pi\ell,r\sin\frac\pi\ell,\mathbf0\right).
\end{equation}
\begin{lemma}\label{lem:ce-midpoint-positive}
For the selected \(\gamma,\alpha\) and every large even
\(\ell\), there is \(c_{n,\ell}>0\) such that
\[
 \mathcal P_n(r)\geq c_{n,\ell}r^{2+\alpha},\,\  r>0.
\]
\end{lemma}
\begin{proof}
Lemma~\ref{lem:large-ell-lattice-positive} and
\eqref{eq:large-ell-joint-cancellation} give
\[
 \mathcal P_n(r)\geq c_{n,\ell}r^{2+\alpha},
 \,\  \gamma\leq r\leq1.
\]
The identity
\(\mathcal P_n(\gamma r)=\gamma^{2+\alpha}\mathcal P_n(r)\)
gives \(\mathcal P_n(r)\geq c_{n,\ell}r^{2+\alpha}\) for all \(r>0\).
\end{proof}

To locate the angular minimum of \(\Phi_{\gamma^k}(\bar x,0)\), we use the
following averaging identity.
\begin{lemma}\label{lem:ce-vertices}
For every polynomial \(P\) on \(\mathbb R^2\) of degree less than \(\ell\),
\[
 \frac1\ell\sum_{i=1}^{\ell}P(\bar\xi_i)
 =\frac1{2\pi}\int_0^{2\pi}P(\cos\theta,\sin\theta)\,d\theta.
\]
\end{lemma}
\begin{proof}
Every nonconstant sine or cosine term of order less than \(\ell\) has
zero average over \(\{\bar\xi_i\}_{i=1}^{\ell}\) and over the circle.
\end{proof}
The sum
\(b\sum_{k\in\mathbb Z}\gamma^{k(n-2+\alpha)}\Phi_{\gamma^k}(x)\)
is understood away from the points \(\gamma^k\xi_i\).
We now extend the lower bound for \(\mathcal P_n\) to all angles in the plane
\(x''=0\).
\begin{lemma}\label{lem:ce-trace}
For every \(\bar x\ne0\) with \(\bar x\ne\gamma^k\bar\xi_i\)
for \(k\in\mathbb Z\), \(i=1,\ldots,\ell\),
\begin{equation}\label{eq:ce-trace-lower}
 \sum_{k\in\mathbb Z}\gamma^{k(n-2+\alpha)}
                  \Phi_{\gamma^k}(\bar x,0)
       >c_{n,\ell}|\bar x|^{2+\alpha}.
\end{equation}
\end{lemma}
\begin{proof}
For circle \(k\), apply Hermite interpolation to
\(z\mapsto(1+s^2-2sz)^{\frac{4-n}{2}}\), with \(s=\frac{|\bar x|}{\gamma^k}\),
at \(z=\cos\frac{(2j-1)\pi}{\ell}\), \(j=1,\ldots,\frac{\ell}{2}\).
The remainder is nonnegative, and the average of the Hermite
interpolation polynomial of degree less than \(\ell\)
over \(\{\bar\xi_i\}_{i=1}^{\ell}\) is independent of rotation by
Lemma~\ref{lem:ce-vertices}.
Thus, for fixed \(|\bar x|\), \(\Phi_{\gamma^k}(\bar x,0)\) is minimized
at polar angle \(\frac{\pi}{\ell}\), by \eqref{eq:large-ell-hermite-remainder}.
Summing \(\gamma^{k(n-2+\alpha)}\Phi_{\gamma^k}(\bar x,0)\)
over \(k\in\mathbb Z\), using absolute convergence and
Lemma~\ref{lem:ce-midpoint-positive}, gives
\eqref{eq:ce-trace-lower} after reducing \(c_{n,\ell}\), if necessary.
\end{proof}

We now use the lower bound \eqref{eq:ce-trace-lower} to estimate
\(U_{\lambda,r}\) and \(H_{k_0}\) from \eqref{eq:ce-background-decomposition}.
\begin{lemma}\label{lem:ce-profile}
For small fixed \(\varepsilon_0>0\) and then
small \(b>0\), the following estimates hold uniformly for parameters
satisfying \eqref{eq:ce-parameter-space} with
\(0\leq\varepsilon\leq\varepsilon_0\).

The ansatz has the positive lower bound
\begin{equation}\label{eq:ce-profile-positive}
 U_{\lambda,r}(x)\geq C^{-1} b|\bar x|^{2+\alpha}+
                    \frac{b^{\frac{n-4}{n-3}}}4|x''|^2>0,
 \,\  x\in\Omega\setminus\{0\}.
\end{equation}

For every \(k_0\geq0\), the function \(H_{k_0}\) in
\eqref{eq:ce-background-decomposition} satisfies
\begin{equation}\label{eq:ce-background-taylor}
\begin{split}
 |H_{k_0}(\mathbf x_{k_0,1})|
 &\leq C(\varepsilon br_{k_0}^{2+\alpha}+b^{\frac{n}{n-4}}r_{k_0}^4),\\
 |DH_{k_0}(\mathbf x_{k_0,1})|
 &\leq C(\varepsilon br_{k_0}^{1+\alpha}+b^{\frac{n}{n-4}}r_{k_0}^3),\\
 \sup_{x\in\mathbf B_{4\eta r_{k_0}}(\mathbf x_{k_0,1})}
 |D^2H_{k_0}(x)|&\leq C b^{\frac{n-4}{n-3}}.
\end{split}
\end{equation}
Here \(\eta>0\) is fixed small enough. The corresponding
estimates for every \(i\) follow by symmetry. The support radii of
\(\chi_{k,i}\) in Subsection~\ref{subsec:ce-formulation}
are chosen so that \eqref{eq:ce-background-taylor} holds on their supports.

Outside \(\bigcup_{k\geq0}\bigcup_{i=1}^{\ell}
\mathbf B_{\eta r_k}(\mathbf x_{k,i})\), the ansatz satisfies
\begin{equation}\label{eq:ce-profile-upper}
 |U_{\lambda,r}(x)|\leq C b^{\frac{n-4}{n-3}}|x|^2,
 \,\  x\in\Omega\setminus\{0\}.
\end{equation}
\end{lemma}
\begin{proof}
To prove \eqref{eq:ce-profile-positive}, it suffices, as in the
proof of Lemma~\ref{lem:ce-trace}, to consider
\[
 \bar x=|\bar x|\left(\cos\frac\pi\ell,\sin\frac\pi\ell\right).
\]
For \(\bar x\ne0\), Taylor expansion and summation using
\eqref{eq:ce-parameter-space} give
\begin{equation}\label{eq:ce-profile-plane-comparison}
 \left|(\mathbf h_{r,\lambda}+h_{r,\lambda})(\bar x,0)
       -b\mathcal P_n(|\bar x|)\right|
 \leq Cb\left(b^{\frac4{n-4}}+\varepsilon\right)|\bar x|^{2+\alpha}
\end{equation}
and
\begin{equation}\label{eq:ce-profile-quadratic-bound}
 \left|\frac12\partial_{x_n}^2
       (\mathbf h_{r,\lambda}+h_{r,\lambda})(\bar x,0)\right|\leq Cb.
\end{equation}
Here \(\frac12\partial_{x_n}^2(\mathbf h_{r,\lambda}+h_{r,\lambda})(\bar x,0)\)
is the coefficient of \(|x''|^2\). Each summand in
\(\mathbf h_{r,\lambda}+h_{r,\lambda}\) is convex as a function of
\(|x''|^2\). Thus Lemma~\ref{lem:ce-midpoint-positive} together with
\eqref{eq:ce-profile-plane-comparison} and
\eqref{eq:ce-profile-quadratic-bound} implies
\[
\begin{aligned}
U_{\lambda,r}(x)
 &\geq(\mathbf h_{r,\lambda}+h_{r,\lambda})(\bar x,0) +\left[b^{\frac{n-4}{n-3}}
       +\frac12\partial_{x_n}^2
       (\mathbf h_{r,\lambda}+h_{r,\lambda})(\bar x,0)\right]|x''|^2\\
 &\geq b\left[c_{n,\ell}-C\left(b^{\frac4{n-4}}+\varepsilon\right)\right]
       |\bar x|^{2+\alpha}
       +\left(b^{\frac{n-4}{n-3}}-Cb\right)|x''|^2\\
 &\geq\frac{c_{n,\ell}}2b|\bar x|^{2+\alpha}
       +\frac14b^{\frac{n-4}{n-3}}|x''|^2,
\end{aligned}
\]
for small fixed \(\varepsilon_0\) and then small \(b\).
At \(\bar x=0\), the definition of \(J_2\) and convexity give
\[
 U_{\lambda,r}(0,x'')\geq b^{\frac{n-4}{n-3}}|x''|^2.
\]
This proves \eqref{eq:ce-profile-positive}.

For \eqref{eq:ce-background-taylor}, fix \(k_0\geq0\) and use the
definition of \(H_{k_0}\) in \eqref{eq:ce-background-decomposition}.
For the parameters in Lemma~\ref{lem:ce-background-comparison},
\eqref{eq:ce-background-comparison} and \(H(1)=H'(1)=0\) give
\[
 |H_{k_0}(\gamma^{k_0}\xi_1)|
 +\gamma^{k_0}|DH_{k_0}(\gamma^{k_0}\xi_1)|
 \leq Cb^{\frac{n}{n-4}}\gamma^{4k_0}.
\]
Only the radial component of \(DH_{k_0}(\gamma^{k_0}\xi_1)\) can be nonzero,
by symmetry.

Varying \((\lambda,r)\) subject to
\eqref{eq:ce-parameter-space}, while moving the
evaluation point to \(\mathbf x_{k_0,1}\), changes \(H_{k_0}\) by at most
\(C\varepsilon br_{k_0}^{2+\alpha}\) in the value and
\(C\varepsilon br_{k_0}^{1+\alpha}\) in the gradient.
Indeed, after rescaling \(x=r_{k_0}y\), the points
\(\frac{r_j\xi_i}{r_{k_0}}\), \(j\geq-3\), with \((j,i)\ne(k_0,1)\), stay a fixed distance
from \(\xi_1\). Differentiate \(H_{k_0}(r_{k_0}\xi_1)\)
with respect to \(\lambda_j,r_j\), \(j\geq-3\), and split
the sums defining \(H_{k_0}\) into
\(r_j\geq2r_{k_0}\), \(r_j\leq \frac{r_{k_0}}{2}\), and
\(\frac{r_{k_0}}{2}<r_j<2r_{k_0}\). For \(j\leq-4\), use the fixed radius
\(\gamma^j\) in the same split. This gives the bound for
\(H_{k_0}(\mathbf x_{k_0,1})\).
One spatial derivative gives the factor \(r_{k_0}^{-1}\) in the
bound for \(DH_{k_0}(\mathbf x_{k_0,1})\).
The added term \(b^{\frac{n-4}{n-3}}|x''|^2\) has zero value and gradient
at \(\mathbf x_{k_0,1}\). Combining the bounds
\(C\varepsilon br_{k_0}^{2+\alpha}\) and
\(C\varepsilon br_{k_0}^{1+\alpha}\) with
\eqref{eq:ce-background-comparison}, and using
\(C^{-1}\gamma^{k_0}\leq r_{k_0}\leq C\gamma^{k_0}\),
proves the first two estimates in
\eqref{eq:ce-background-taylor}.

For \(x\in\mathbf B_{4\eta r_{k_0}}(\mathbf x_{k_0,1})\), split the
sum defining \(D^2H_{k_0}\) at radii \(\frac{r_{k_0}}{2}\) and \(2r_{k_0}\).
The points \(\mathbf x_{j,i}\) with \((j,i)\ne(k_0,1)\) and
\(\frac{r_{k_0}}{2}<r_j<2r_{k_0}\) are separated from
\(\mathbf B_{4\eta r_{k_0}}(\mathbf x_{k_0,1})\).
In the sums defining \(\mathbf h_{r,\lambda}\) and \(h_{r,\lambda}\),
subtracting \(J_2\) leaves fourth-order Taylor remainders for
\(r_j\geq2r_{k_0}\) and for \(\gamma^j\geq2r_{k_0}\), \(j\leq-4\). Summing the geometric series, using \(0<\alpha<2\), gives
\[
 \left|D^2\bigl(H_{k_0}-b^{\frac{n-4}{n-3}}|x''|^2\bigr)(x)\right|
 \leq C b r_{k_0}^{\alpha}.
\]
The Hessian of \(b^{\frac{n-4}{n-3}}|x''|^2\) is bounded by
\(Cb^{\frac{n-4}{n-3}}\). Since \(r_{k_0}\leq C\) and
\(b\leq b^{\frac{n-4}{n-3}}\), this proves the third estimate in \eqref{eq:ce-background-taylor}.

To prove \eqref{eq:ce-profile-upper}, let
\(x\notin\bigcup_{k,i}\mathbf B_{\eta r_k}(\mathbf x_{k,i})\).
Split the sums defining \(\mathbf h_{r,\lambda}+h_{r,\lambda}\)
into \(r_k\geq2|x|\), \(r_k\leq\frac{|x|}{2}\), and
\(\frac{|x|}{2}<r_k<2|x|\), using \(\gamma^k\) in place of \(r_k\) for
\(k\leq-4\). For \(\frac{|x|}{2}<r_k<2|x|\), use
\(|x-\mathbf x_{k,i}|\geq\eta r_k\).
For \(r_k\geq2|x|\), the subtraction of \(J_2\) leaves a
fourth-order remainder. For \(r_k\leq\frac{|x|}{2}\), estimate
\(\mathcal B_{\lambda_k,\mathbf x_{k,i}}\) and its constant and quadratic
Taylor terms separately, and do the same for the terms of
\(h_{r,\lambda}\). Summation, using \(0<\alpha<2\), yields
\[
 |\mathbf h_{r,\lambda}(x)+h_{r,\lambda}(x)|
 \leq Cb|x|^{2+\alpha}.
\]
Adding \(b^{\frac{n-4}{n-3}}|x''|^2\) and using the boundedness of
\(\Omega\) and \(b\leq b^{\frac{n-4}{n-3}}\) gives \eqref{eq:ce-profile-upper}.
\end{proof}

\section{The constrained problem}\label{sec:ce-pde}
\subsection{The equation for \texorpdfstring{\(\phi\)}{phi}}\label{subsec:ce-formulation}
Fix \((\lambda,r)\) satisfying \eqref{eq:ce-parameter-space}, and let
\(U_{\lambda,r}\) be the ansatz constructed in
Section~\ref{sec:exact-counterexample}. We seek a solution of
\eqref{eq:biharmonic-problem} of the form \(u=U_{\lambda,r}+\phi\).
For \eqref{eq:biharmonic-problem}, define the linearized operator
\(\mathcal L_{\lambda,r}\) at \(U_{\lambda,r}\), the error \(\mathcal E\),
and the nonlinear term \(\mathcal N(\phi)\) by
\[
\begin{aligned}
 \mathcal L_{\lambda,r}\phi&=\Delta^2\phi-\frac{n+4}{n-4}
                       U_{\lambda,r}^{\frac8{n-4}}\phi,
 \,\ \mathcal E=U_{\lambda,r}^{\frac{n+4}{n-4}}-\Delta^2U_{\lambda,r},\\
 \mathcal N(\phi)&=(U_{\lambda,r}+\phi)^{\frac{n+4}{n-4}}
                  -U_{\lambda,r}^{\frac{n+4}{n-4}}
                  -\frac{n+4}{n-4}U_{\lambda,r}^{\frac8{n-4}}\phi.
\end{aligned}
\]
Thus solving \eqref{eq:biharmonic-problem} is equivalent to finding
\(\phi\) such that
\[
 \mathcal L_{\lambda,r}\phi=\mathcal E+\mathcal N(\phi),
 \,\ U_{\lambda,r}+\phi>0\quad\text{in }\Omega\setminus\{0\}.
\]

In this section, \(\sum_{k,i}\) denotes
\(\sum_{k=0}^{\infty}\sum_{i=1}^{\ell}\). Constants may depend
on the fixed \(n,\ell,\gamma,\alpha,R,\rho,\eta,\varepsilon_0\),
but are independent of \(k\), small \(b>0\), and
\((\lambda,r)\) satisfying \eqref{eq:ce-parameter-space} with
\(0\leq\varepsilon\leq\varepsilon_0\).
The notation \(o_b(1)\) denotes a quantity tending to zero as
\(b\rightarrow0\), uniformly in \(k,\lambda,r\).
Put \(\tau=\frac{2(n-4)}{n-2}\), and set
\begin{equation}\label{eq:ce-linear-norms}
\begin{aligned}
 \|\phi\|_*&=\mathop{\rm ess\,sup}_{x\in\Omega\setminus\{0\}}
 \left(\sum_{k,i}
 \frac{br_k^{2+\alpha}\lambda_k^\tau}
 {(\lambda_k+|x-\mathbf x_{k,i}|)^\tau}\right)^{-1}|\phi(x)|,\\
 \|f\|_{**}&=\mathop{\rm ess\,sup}_{x\in\Omega\setminus\{0\}}
 \left(\sum_{k,i}
 \frac{br_k^{2+\alpha}\lambda_k^\tau}
 {(\lambda_k+|x-\mathbf x_{k,i}|)^{\tau+4}}\right)^{-1}|f(x)|.
\end{aligned}
\end{equation}
As in \cite[Section~2]{WY}, the weights in the definitions of
\(\|\cdot\|_*\) and \(\|\cdot\|_{**}\) in
\eqref{eq:ce-linear-norms} are chosen together from the decay of
the fundamental solution.

With \(x=(\bar x,x'')\in\mathbb R^2\times\mathbb R^{n-2}\), define
\begin{equation}\label{eq:ce-linear-space}
\begin{aligned}
&\mathcal H=\Bigl\{\phi\in H^2(\Omega):\;
 \phi\text{ is even in }x_2,\text{ radially symmetric in }x'',\\
 &\quad\phi(r\cos\theta,r\sin\theta,x'')
 =\phi\left(r\cos\left(\theta+\frac{2\pi i}{\ell}\right),
 r\sin\left(\theta+\frac{2\pi i}{\ell}\right),x''\right),\quad i=1,\ldots,\ell\Bigr\}.
\end{aligned}
\end{equation}
Let
\[
 Z_{k,i,1}=\lambda_k^{\frac{n-2}{2}}\partial_{\lambda_k}
                  \mathcal B_{\lambda_k,\mathbf x_{k,i}},
 \,\ Z_{k,i,2}=\lambda_k^{\frac{n-2}{2}}\partial_{r_k}
                  \mathcal B_{\lambda_k,\mathbf x_{k,i}},
\]
where \(x,r_k\) are fixed in \(\partial_{\lambda_k}\), and
\(x,\lambda_k\) are fixed in \(\partial_{r_k}\).
Choose a fixed smooth radial function defining \(\chi_{k,i}\) by
translation to \(\mathbf x_{k,i}\) and scaling by \(r_k\), with
\(0\leq\chi_{k,i}\leq1\),
\[
 \chi_{k,i}=1\quad\text{on }\mathbf B_{2\eta r_k}(\mathbf x_{k,i}),
 \,\  \operatorname{supp}\chi_{k,i}
 \subset\mathbf B_{4\eta r_k}(\mathbf x_{k,i}).
\]
Choose \(\eta>0\) small enough that
\(\mathbf B_{4\eta r_k}(\mathbf x_{k,i})\) are disjoint and
contained in \(\Omega\) for every \((\lambda,r)\) satisfying
\eqref{eq:ce-parameter-space}. We take \(b\) small enough that
\(\frac{\lambda_k}{r_k}\ll\eta\).

Let \(G_n\) be the fundamental solution of \(\Delta^2\) in
\(\mathbb R^n\):
\[
 G_n(x)=\frac{|x|^{4-n}}{2(n-4)(n-2)|\mathbb S^{n-1}|},
 \,\  \Delta^2G_n=\delta_0.
\]
For \(f\) with \(\|f\|_{**}<\infty\), write
\[
\begin{aligned}
 (G_n*f)(x)&:=\int_\Omega G_n(x-y)f(y)\,dy,\\
 J_2[G_n*f](x)&:=\int_\Omega
 \left[G_n(-y)+DG_n(-y)\cdot x
       +\frac12D^2G_n(-y)[x,x]\right]f(y)\,dy,
\end{aligned}
\]
and define
\begin{equation}\label{eq:ce-normalized-potential}
 \mathcal Tf(x)=(I-J_2)[G_n*f](x),
\end{equation}
so that \(\Delta^2\mathcal Tf=f\) in \(\Omega\). Write \(\mathbf c=(c_{k,1},c_{k,2})_{k\geq0}\) and
\(\|\mathbf c\|_\infty=\sup_{k\geq0,\ \nu=1,2}|c_{k,\nu}|\).
We seek \(\phi\) and \(\mathbf c\) with
\(\|\phi\|_*<\infty\) and \(\|\mathbf c\|_\infty<\infty\)
satisfying
\begin{equation}\label{eq:ce-constrained-correction}
\left\{\begin{aligned}
 &\mathcal L_{\lambda,r}\phi=\mathcal E+\mathcal N(\phi)
       +\frac{n+4}{n-4}\sum_{k,i,\nu}c_{k,\nu}br_k^{2+\alpha}\chi_{k,i}
       \mathcal B_{\lambda_k,\mathbf x_{k,i}}^{\frac8{n-4}}Z_{k,i,\nu}
              &&\text{in }\Omega\setminus\{0\},\\
 &\phi\in\mathcal H,\,\  \phi=\mathcal T(\Delta^2\phi),\\
 &\left\langle
   \chi_{k,i}\mathcal B_{\lambda_k,\mathbf x_{k,i}}^{\frac8{n-4}}Z_{k,i,\nu},\phi
  \right\rangle=0 &&(k\geq0,\ i=1,\ldots,\ell,\ \nu=1,2).
\end{aligned}\right.
\end{equation}
Here \(\langle v,w\rangle=\int_\Omega vw\).

In Subsection~\ref{subsec:ce-linear}, we introduce the linear
constrained problem \eqref{eq:ce-linear-problem} for \(\mathcal L_{\lambda,r}\)
and establish its solvability. In Subsection~\ref{subsec:ce-nonlinear},
we estimate the nonlinear term \(\mathcal N(\phi)\) and the error \(\mathcal E\), and then use the contraction mapping principle
to solve the nonlinear constrained problem
\eqref{eq:ce-constrained-correction}. In Section~\ref{sec:ce-reduction},
we choose \(\lambda,r\) so that \(\mathbf c\) vanishes.

\subsection{The linear constrained problem}\label{subsec:ce-linear}
For a measurable \(f\) with the symmetries in
\eqref{eq:ce-linear-space} and finite \(\|f\|_{**}\), consider
\begin{equation}\label{eq:ce-linear-problem}
\left\{\begin{aligned}
 &\mathcal L_{\lambda,r}\phi=f+\frac{n+4}{n-4}\sum_{k,i,\nu}c_{k,\nu}br_k^{2+\alpha}\chi_{k,i}
       \mathcal B_{\lambda_k,\mathbf x_{k,i}}^{\frac8{n-4}}Z_{k,i,\nu}
              &&\text{in }\Omega\setminus\{0\},\\
 &\phi\in\mathcal H,\,\  \phi=\mathcal T(\Delta^2\phi),\\
 &\left\langle
   \chi_{k,i}\mathcal B_{\lambda_k,\mathbf x_{k,i}}^{\frac8{n-4}}Z_{k,i,\nu},\phi
  \right\rangle=0 &&(k\geq0,\ i=1,\ldots,\ell,\ \nu=1,2).
\end{aligned}\right.
\end{equation}

\begin{lemma}\label{lem:ce-linear-estimate}
For small \(b>0\), any solution of the linear constrained problem
\eqref{eq:ce-linear-problem} with \(\|\phi\|_*<\infty\) and
\(\|\mathbf c\|_\infty<\infty\) satisfies
\begin{equation}\label{eq:ce-linear-bound}
 \|\phi\|_*+\|\mathbf c\|_\infty+\|\Delta^2\phi\|_{**}
 \leq C\|f\|_{**}.
\end{equation}
Moreover,
\begin{equation}\label{eq:ce-linear-pointwise}
 |\phi(x)|\leq Cb|x|^{2+\alpha}\|f\|_{**}.
\end{equation}
The constants are independent of \(k,b,\lambda,r\).
\end{lemma}
\begin{proof}
Since \(0<\tau<n-4\), the convolution inequality in
\cite[Lemma~A.2]{GLN}, with \(m=2\), gives
\begin{equation}\label{eq:ce-bubble-convolution}
 \int_{\mathbb R^n}\frac{|x-y|^{4-n}\lambda_k^\tau}
 {(\lambda_k+|y-\mathbf x_{k,i}|)^{\tau+4}}\,dy
 \leq\frac{C\lambda_k^\tau}
 {(\lambda_k+|x-\mathbf x_{k,i}|)^\tau}.
\end{equation}
By symmetry it suffices to take \(i=1\). We first prove
\begin{equation}\label{eq:ce-linear-local-bound}
\sup_{\mathbf B_{4\eta r_k}(\mathbf x_{k,1})}|\phi|
       +br_k^{2+\alpha}\sum_{\nu=1}^2|c_{k,\nu}| \leq C\sup_{\mathbf B_{4\eta r_k}(\mathbf x_{k,1})}
 \left|\phi-G_n*\bigl[\chi_{k,1}(\Delta^2\phi-f)\bigr]\right|.
\end{equation}
We argue by contradiction, as in \cite[Lemma~2.1]{WY}.
If \eqref{eq:ce-linear-local-bound} fails, take a sequence with
\(b\rightarrow0\). Divide \(\phi,f,\mathbf c\) by the left-hand side
of \eqref{eq:ce-linear-local-bound}, so that it equals one and the
right-hand side tends to zero. Rescale
\(x=\mathbf x_{k,1}+\lambda_k y\).
The estimate \eqref{eq:ce-background-taylor} gives
\[
 \sup_x\int_{\mathbf B_{4\eta r_k}(\mathbf x_{k,1})}
 |x-y|^{4-n}
 \left|U_{\lambda,r}^{\frac8{n-4}}(y)
       -\mathcal B_{\lambda_k,\mathbf x_{k,1}}^{\frac8{n-4}}(y)\right|\,dy
 =o_b(1).
\]
By \eqref{eq:ce-bubble-convolution}, the convolution
\(G_n*[\chi_{k,1}(\Delta^2\phi-f)]\), after the rescaling
\(x=\mathbf x_{k,1}+\lambda_k y\) and normalization in
\eqref{eq:ce-linear-local-bound}, gives the bound
\(C(1+|y|)^{-\tau}+o(1)\) for the rescaled \(\phi\) on
\(|y|<\frac{4\eta r_k}{\lambda_k}\).
Local compactness yields a limit \(\psi\), and the coefficients
\(br_k^{2+\alpha}c_{k,\nu}\) converge to \(a_\nu\).
The limiting equation is
\begin{equation}\label{eq:ce-rescaled-linear-equation}
 \Delta^2\psi-\frac{n+4}{n-4}\mathcal B_{1,0}^{\frac8{n-4}}\psi
 =\frac{n+4}{n-4}\sum_{\nu=1}^2a_\nu
                       \mathcal B_{1,0}^{\frac8{n-4}}Z_\nu,
 \,\  |\psi(y)|\leq\frac C{(1+|y|)^\tau},
\end{equation}
where
\[
 Z_1(y)=\frac{(n-4)A_n(|y|^2-1)}{2(1+|y|^2)^{\frac{n-2}{2}}},
 \,\ Z_2(y)=\frac{(n-4)A_n y_1}{(1+|y|^2)^{\frac{n-2}{2}}}.
\]
Passing to the limit in
\(G_n*[\chi_{k,1}(\Delta^2\phi-f)]\), with
\(x=\mathbf x_{k,1}+\lambda_k y\), gives
\begin{equation}\label{eq:ce-rescaled-linear-potential}
 \psi(y)=\frac{n+4}{n-4}\int_{\mathbb R^n}G_n(y-z)
 \mathcal B_{1,0}^{\frac8{n-4}}(z)
 \left(\psi(z)+\sum_{\nu=1}^2a_\nu Z_\nu(z)\right)\,dz.
\end{equation}
Since \(\mathcal B_{1,0}^{\frac{8}{n-4}}(z)\leq C(1+|z|)^{-8}\),
repeated use of convolution inequalities, followed by differentiation
under the integral sign in \eqref{eq:ce-rescaled-linear-potential}, gives
\[
 |D^j\psi(y)|\leq\frac C{(1+|y|)^{n-4+j}},\,\  j=0,1,2.
\]
Approximation by compactly supported functions then gives
\(\psi\in D_0^{2,2}(\mathbb R^n)\).

For \(L\geq1\) and \(\nu=1,2\), scaling gives
\[
\begin{aligned}
&\int_{\Omega\setminus\mathbf B_{L\lambda_k}(\mathbf x_{k,1})}
 \chi_{k,1}\mathcal B_{\lambda_k,\mathbf x_{k,1}}^{\frac8{n-4}}
 |Z_{k,1,\nu}\phi| \leq C\lambda_k^{n-4}
 \sup_{\mathbf B_{4\eta r_k}(\mathbf x_{k,1})}|\phi|
 \int_{|y|>L}\mathcal B_{1,0}^{\frac8{n-4}}|Z_\nu(y)|\,dy\\
 &\,\ \leq CL^{-4}\lambda_k^{n-4}
 \sup_{\mathbf B_{4\eta r_k}(\mathbf x_{k,1})}|\phi|.
\end{aligned}
\]
Thus we can pass to the limit in the orthogonality conditions
in \eqref{eq:ce-linear-problem}, after the rescaling
\(x=\mathbf x_{k,1}+\lambda_k y\) and normalization in
\eqref{eq:ce-linear-local-bound}, to obtain
\[
 \int_{\mathbb R^n}\mathcal B_{1,0}^{\frac8{n-4}}Z_\nu\psi=0,
 \,\  \nu=1,2.
\]
Testing \eqref{eq:ce-rescaled-linear-equation} with
\(Z_1,Z_2\) gives \(a_1=a_2=0\).
By symmetry and the non-degeneracy result of Lu and Wei \cite{LW00},
\(\psi\in\operatorname{span}\{Z_1,Z_2\}\).
The conditions
\(\int_{\mathbb R^n}\mathcal B_{1,0}^{\frac8{n-4}}Z_\nu\psi=0\),
\(\nu=1,2\), imply \(\psi=0\).
Uniform convergence on compact sets to \(\psi=0\) and the bound
\(C(1+|y|)^{-\tau}+o(1)\) for \(|y|<\frac{4\eta r_k}{\lambda_k}\)
contradict the normalization of the left-hand side of
\eqref{eq:ce-linear-local-bound}.

We now use \eqref{eq:ce-linear-local-bound} to prove
\eqref{eq:ce-linear-bound} and \eqref{eq:ce-linear-pointwise}.
The sum in the definition of \(\|\phi\|_*\) satisfies
\begin{equation}\label{eq:ce-weight-bounds}
 C^{-1}b^{1+\frac4{n-2}}
 \leq\sum_{k,i}\frac{br_k^{2+\alpha}\lambda_k^\tau}
       {(\lambda_k+|x-\mathbf x_{k,i}|)^\tau}\leq Cb.
\end{equation}
In particular, \(\|\phi\|_\infty\leq Cb\|\phi\|_*\).
Since \(\lambda_k+|x-\mathbf x_{k,i}|\) is uniformly bounded on
\(\Omega\), \eqref{eq:ce-weight-bounds} also gives
\begin{equation}\label{eq:ce-source-weight-lower}
 \sum_{k,i}\frac{br_k^{2+\alpha}\lambda_k^\tau}
 {(\lambda_k+|x-\mathbf x_{k,i}|)^{\tau+4}}
 \geq C^{-1}b^{1+\frac4{n-2}}.
\end{equation}

For \(\mathcal Tf\), we have
\begin{equation}\label{eq:ce-free-potential-bound}
 \|\mathcal Tf\|_*\leq C\|f\|_{**},
 \,\ |\mathcal Tf(x)|\leq Cb|x|^{2+\alpha}\|f\|_{**}.
\end{equation}
Indeed, decompose
\[
 f=\left(1-\sum_{k,i}\chi_{k,i}\right)f
       +\sum_{k,i}\chi_{k,i}f.
\]
Since the supports of \(\chi_{k,i}\) are disjoint,
\eqref{eq:ce-parameter-space} and \eqref{eq:ce-linear-norms} give
\[
\begin{aligned}
 \left|\left(1-\sum_{k,i}\chi_{k,i}(x)\right)f(x)\right|
 &\leq\|f\|_{**}\sum_{k,i}(1-\chi_{k,i}(x))
 \frac{br_k^{2+\alpha}\lambda_k^\tau}
 {(\lambda_k+|x-\mathbf x_{k,i}|)^{\tau+4}}\\
 &\leq Cb^{1+\frac4{n-2}}\|f\|_{**}.
\end{aligned}
\]
By symmetry, the odd-order terms in the Taylor expansion of
\(G_n(x-y)\) at \(x=0\) integrate to zero against \(f(y)\). Taylor's theorem and
\eqref{eq:ce-bubble-convolution} give
\begin{equation}\label{eq:ce-source-potential-parts}
\begin{aligned}
 &\left|\mathcal T\left[\left(1-\sum_{k,i}\chi_{k,i}\right)f\right](x)\right|
 \leq Cb^{1+\frac4{n-2}}\|f\|_{**}|x|^4
       \left(1+\log\frac{2R}{|x|}\right),\\
 &\left|\sum_{k,i}\mathcal T(\chi_{k,i}f)(x)\right|
 \leq Cb\|f\|_{**}|x|^{2+\alpha}.
\end{aligned}
\end{equation}
For the second estimate in \eqref{eq:ce-source-potential-parts},
split the sum into three sums over \(r_k\geq2|x|\),
\(r_k\leq\frac{|x|}{2}\), and \(\frac{|x|}{2}<r_k<2|x|\).
Taylor's theorem applies when \(r_k\geq2|x|\).
For \(r_k\leq\frac{|x|}{2}\), estimate \(G_n*(\chi_{k,i}f)\)
and \(J_2[G_n*(\chi_{k,i}f)]\) separately.
For \(\frac{|x|}{2}<r_k<2|x|\), use \eqref{eq:ce-bubble-convolution}.
Since \(\alpha<2\), \eqref{eq:ce-source-potential-parts} proves
the pointwise estimate in \eqref{eq:ce-free-potential-bound}.

For the \(*\)-norm estimate in \eqref{eq:ce-free-potential-bound},
\eqref{eq:ce-bubble-convolution} gives
\(\|G_n*f\|_*\leq C\|f\|_{**}\). Moreover,
\[
 \left|\int_\Omega D^jG_n(-y)f(y)\,dy\right|
 \leq Cb^{1+\frac4{n-2}}\|f\|_{**},\,\  j=0,1,2.
\]
Thus \eqref{eq:ce-weight-bounds} gives
\(\|J_2(G_n*f)\|_*\leq C\|f\|_{**}\), and we use
\(\mathcal Tf=G_n*f-J_2(G_n*f)\).

To estimate \(\phi\), write
\[
 \phi=\mathcal Tf+\mathcal T(\Delta^2\phi-f).
\]
We first bound \(\Delta^2\phi-f\) using the linear equation
\eqref{eq:ce-linear-problem}. Condition
\eqref{eq:ce-parameter-space} gives
\(\lambda_k^{\frac{n-4}{2}}\leq Cbr_k^{n-2+\alpha}\).
Direct integration, using \eqref{eq:ce-background-taylor}, gives
\[
 \int_\Omega |x|^{-n}U_{\lambda,r}^{\frac8{n-4}}\,dx
 +\sum_{k,i}\int_\Omega |x|^{-n}\chi_{k,i}
       \mathcal B_{\lambda_k,\mathbf x_{k,i}}^{\frac8{n-4}}\,dx
 =o_b(1).
\]
By \eqref{eq:ce-linear-problem},
\(br_k^{2+\alpha}\leq Cb\), \(|Z_{k,i,\nu}|\leq C\), and
\(\|\phi\|_\infty\leq Cb\|\phi\|_*\), we obtain
\begin{equation}\label{eq:ce-linear-source-moment}
\begin{aligned}
&\int_\Omega |x|^{-n}|\Delta^2\phi-f|\,dx \leq C\|\phi\|_\infty
       \int_\Omega |x|^{-n}U_{\lambda,r}^{\frac8{n-4}}\,dx +Cb\|\mathbf c\|_\infty
       \sum_{k,i}\int_\Omega |x|^{-n}\chi_{k,i}
       \mathcal B_{\lambda_k,\mathbf x_{k,i}}^{\frac8{n-4}}\,dx\\
 &\quad\leq b\,o_b(1)(\|\phi\|_*+\|\mathbf c\|_\infty).
\end{aligned}
\end{equation}

We now use \eqref{eq:ce-linear-source-moment} to estimate the remaining term
\(\mathcal T(\Delta^2\phi-f)\). The weighted integral controls the
normalized kernel away from its singularity; near the evaluation point,
we also need a pointwise bound outside the bubble cores. Write
\[
 g=\Delta^2\phi-f,\,\  N=\|\phi\|_*+\|\mathbf c\|_\infty,
 \,\  q=\frac8{n-4}.
\]
By \eqref{eq:ce-linear-source-moment}, there is a quantity
\(\varepsilon_b\to0\), uniform in the admissible parameters, such that
\[
 A_g:=\int_\Omega |y|^{-n}|g(y)|\,dy\leq b\varepsilon_bN.
\]
The symmetries in \eqref{eq:ce-linear-space} and the evenness of \(\ell\)
give \(g(-y)=g(y)\) a.e. We extend \(g\) by zero outside
\(\Omega=\mathbf B_R\), preserving this symmetry.

Take the fixed \(\eta\) small enough, in particular \(\eta<\frac{1}{16}\),
and let \(D_{1/2}=\bigcup_{m,i}\mathbf B_{\frac{\eta r_m}{2}}(\mathbf x_{m,i})\).
We claim that, with \(\delta_b=b^{\frac{8}{n-3}}+b^{\frac{8}{n-4}}\),
\[
 |g(y)|\leq Cb\delta_bN\,\ \text{for a.e. }y\in\Omega\setminus D_{1/2}.
\]
Indeed, \eqref{eq:ce-linear-problem}, \eqref{eq:ce-weight-bounds},
\(|Z_{m,i,\nu}|\leq C\), and \(br_m^{2+\alpha}\leq Cb\) imply
\[
 |g(y)|\leq CbN\left(U_{\lambda,r}(y)^q+
 \sum_{m,i}\chi_{m,i}(y)\mathcal B_{\lambda_m,\mathbf x_{m,i}}(y)^q\right).
\]
The bound \eqref{eq:ce-profile-upper} also holds outside \(D_{1/2}\).
To check the additional annuli
\(\frac{\eta r_m}{2}\leq|y-\mathbf x_{m,i}|<\eta r_m\), use
\eqref{eq:ce-parameter-space} to obtain
\(\mathcal B_{\lambda_m,\mathbf x_{m,i}}(y)\leq Cbr_m^{2+\alpha}\).
For \(H_{m,i}=U_{\lambda,r}-\mathcal B_{\lambda_m,\mathbf x_{m,i}}\),
Taylor's theorem and \eqref{eq:ce-background-taylor} give
\[
 |H_{m,i}(y)|\leq C\left(\varepsilon br_m^{2+\alpha}
 +b^{\frac{n}{n-4}}r_m^4+b^{\frac{n-4}{n-3}}r_m^2\right)
 \leq Cb^{\frac{n-4}{n-3}}r_m^2.
\]
Since \(C^{-1}r_m\leq |y|\leq Cr_m\) there, it follows that
\(U_{\lambda,r}(y)\leq Cb^{\frac{n-4}{n-3}}|y|^2\) outside \(D_{1/2}\),
and hence \(U_{\lambda,r}(y)^q\leq Cb^{\frac{8}{n-3}}\).
Moreover, if \(\chi_{m,i}(y)\ne0\) and \(y\notin D_{1/2}\), then
\[
 \mathcal B_{\lambda_m,\mathbf x_{m,i}}(y)^q
 \leq C\lambda_m^4r_m^{-8}
 \leq Cb^{\frac{8}{n-4}}r_m^{\frac{8(2+\alpha)}{n-4}}
 \leq Cb^{\frac{8}{n-4}}.
\]
Here the radii are uniformly bounded above. The supports of the
\(\chi_{m,i}\) are disjoint, so these estimates prove the claimed
pointwise bound.

Fix \(x\ne0\), put \(s=|x|\) and \(\vartheta=\frac{\eta}{16}\), and consider
the two configurations needed below. If
\(x\in\mathbf B_{4\eta r_k}(\mathbf x_{k,1})\), set \(\psi=\chi_{k,1}\);
if \(x\notin\bigcup_{m,i}\mathbf B_{\eta r_m}(\mathbf x_{m,i})\), set
\(\psi=0\). In the first configuration, \(C^{-1}r_k\leq s\leq Cr_k\) and
\(\operatorname{supp}\psi\subset\{|y|\leq2s\}\); the support condition
is automatic in the second. In either configuration,
\[
 |(1-\psi(y))g(y)|\leq Cb\delta_bN
 \,\ \text{for a.e. }y\in\mathbf B_{\vartheta s}(x)\cap\Omega.
\]
To verify this, note first that \(x\) lies outside every
\(\mathbf B_{\eta r_m}(\mathbf x_{m,i})\) except possibly the selected one.
In the first configuration this follows from the disjointness of the
\(\mathbf B_{4\eta r_m}(\mathbf x_{m,i})\).
If \(|y-x|<\vartheta s\) and \(|y-\mathbf x_{m,i}|<\frac{\eta r_m}{2}\), then
\[
 r_m>\frac{(1-\vartheta)s}{1+\frac{\eta}{2}}>\frac{s}{2},
 \,\ |x-\mathbf x_{m,i}|<(2\vartheta+\frac{\eta}{2})r_m<\eta r_m.
\]
Thus such a half-sized core must be the selected one, on which
\(\psi=1\). Everywhere else the preceding pointwise bound applies.

We can now estimate the normalized potential directly. Put
\[
 P_2(x,y)=G_n(-y)+DG_n(-y)\cdot x+\tfrac12D^2G_n(-y)[x,x].
\]
The definition \eqref{eq:ce-normalized-potential} and the support condition
on \(\psi\) give
\[
\begin{aligned}
 \mathcal Tg(x)-G_n*(\psi g)(x)
 &=\int_{|y|\leq2s}G_n(x-y)(1-\psi(y))g(y)\,dy\\
 &\quad-\int_{|y|\leq2s}P_2(x,y)g(y)\,dy\\
 &\quad+\int_{|y|>2s}[G_n(x-y)-P_2(x,y)]g(y)\,dy.
\end{aligned}
\]
All integrals use the zero extension of \(g\). For \(0<|y|\leq2s\),
\(|P_2(x,y)|\leq Cs^4|y|^{-n}\), so the second integral is bounded by
\(Cs^4A_g\). For \(|y|>2s\), Taylor expansion through degree three gives
\[
 G_n(x-y)-P_2(x,y)=\tfrac16D^3G_n(-y)[x,x,x]+R_4(x,y),
 \,\  |R_4(x,y)|\leq Cs^4|y|^{-n}.
\]
The cubic term has zero integral because it is odd in \(y\), whereas
\(g\) and the integration region are centrally symmetric. The third
integral is therefore also bounded by \(Cs^4A_g\). In the first integral,
the part with \(|x-y|<\vartheta s\) is bounded by
\[
 Cb\delta_bN\int_{|z|<\vartheta s}|z|^{4-n}\,dz
 \leq Cb\delta_bNs^4.
\]
On the remaining part, \(|x-y|^{4-n}\leq C_\vartheta s^{4-n}\) and
\(\int_{|y|\leq2s}|g(y)|\,dy\leq(2s)^nA_g\). Consequently, in both
configurations,
\[
 |\mathcal Tg(x)-G_n*(\psi g)(x)|
 \leq Cs^4(A_g+b\delta_bN)\leq b\omega_bs^4N,
 \,\ \omega_b=C(\varepsilon_b+\delta_b)\to0.
\]
The convergence is uniform in \(k\) and all admissible parameters.
Using \(\phi=\mathcal Tf+\mathcal Tg\) and the pointwise estimate in
\eqref{eq:ce-free-potential-bound}, the first configuration, with
\(C^{-1}r_k\leq s\leq Cr_k\), yields
\begin{equation}\label{eq:ce-linear-local-remainder}
\sup_{\mathbf B_{4\eta r_k}(\mathbf x_{k,1})}
 \left|\phi-G_n*\bigl[\chi_{k,1}(\Delta^2\phi-f)\bigr]\right| \leq Cbr_k^{2+\alpha}\|f\|_{**}
       +b\,o_b(1)r_k^4(\|\phi\|_*+\|\mathbf c\|_\infty).
\end{equation}
In the second configuration \(\psi=0\), and the same estimate gives
\begin{equation}\label{eq:ce-linear-away-bound}
 |\phi(x)|\leq Cb|x|^{2+\alpha}\|f\|_{**}
       +b\,o_b(1)|x|^4(\|\phi\|_*+\|\mathbf c\|_\infty).
\end{equation}

Here \(o_b(1)\) may be taken to be \(\omega_b\) in both estimates.
Combining \eqref{eq:ce-linear-local-bound},
\eqref{eq:ce-linear-local-remainder} and \eqref{eq:ce-linear-away-bound},
and using \(\alpha<2\), we obtain
\begin{equation}\label{eq:ce-linear-origin-bound}
 \frac{|\phi(x)|}{b|x|^{2+\alpha}}+\|\mathbf c\|_\infty
 \leq C\left[\|f\|_{**}
       +o_b(1)(\|\phi\|_*+\|\mathbf c\|_\infty)\right].
\end{equation}
For \(\tau<4\), we have
\[
\begin{aligned}
 br_k^{2+\alpha}\mathcal B_{\lambda_k,\mathbf x_{k,i}}^q
 &\leq C\frac{br_k^{2+\alpha}\lambda_k^4}
 {(\lambda_k+|x-\mathbf x_{k,i}|)^8}\\
 &\leq C\frac{br_k^{2+\alpha}\lambda_k^\tau}
 {(\lambda_k+|x-\mathbf x_{k,i}|)^{\tau+4}}.
\end{aligned}
\]
Moreover, for \(x\in\mathbf B_{\eta r_k}(\mathbf x_{k,i})\),
we have \(|x|\leq Cr_k\) and, by \eqref{eq:ce-background-taylor},
\(|H_{k,i}(x)|\leq Cb^{\frac{n-4}{n-3}}r_k^2\). Hence
\[
 b|x|^{2+\alpha}U_{\lambda,r}^q
 \leq C\frac{br_k^{2+\alpha}\lambda_k^\tau}
 {(\lambda_k+|x-\mathbf x_{k,i}|)^{\tau+4}}+Cb^{1+\frac{8}{n-3}}.
\]
For \(x\in\Omega\setminus\bigl(\{0\}\cup
\bigcup_{k,i}\mathbf B_{\eta r_k}(\mathbf x_{k,i})\bigr)\),
\eqref{eq:ce-profile-upper} gives
\[
 b|x|^{2+\alpha}U_{\lambda,r}^q
 \leq Cb^{1+\frac{8}{n-3}}|x|^{2+\alpha+2q}
 \leq Cb^{1+\frac{8}{n-3}}.
\]
Combining \(\frac{8}{n-3}>\frac{4}{n-2}\) with
\eqref{eq:ce-source-weight-lower}, we have, for \(0<b<1\),
\[
 Cb^{1+\frac{8}{n-3}}\leq Cb^{1+\frac{4}{n-2}}
 \leq C\sum_{k,i}\frac{br_k^{2+\alpha}\lambda_k^\tau}
 {(\lambda_k+|x-\mathbf x_{k,i}|)^{\tau+4}}.
\]
Then
\[
 \left|\sum_{k,i,\nu}c_{k,\nu}br_k^{2+\alpha}\chi_{k,i}
       \mathcal B_{\lambda_k,\mathbf x_{k,i}}^q Z_{k,i,\nu}\right|
 \leq C\|\mathbf c\|_\infty
       \sum_{k,i}\frac{br_k^{2+\alpha}\lambda_k^\tau}
       {(\lambda_k+|x-\mathbf x_{k,i}|)^{\tau+4}},
\]
since \(|\chi_{k,i}Z_{k,i,\nu}|\leq C\).
Thus \eqref{eq:ce-linear-problem} and \eqref{eq:ce-linear-origin-bound} give
\[
\begin{aligned}
 \|g\|_{**}
 &\leq C\left(\sup_{x\in\Omega\setminus\{0\}}
          \frac{|\phi(x)|}{b|x|^{2+\alpha}}+\|\mathbf c\|_\infty\right)
 \leq C\bigl(\|f\|_{**}+o_b(1)N\bigr),\\
 \|\Delta^2\phi\|_{**}
 &\leq \|f\|_{**}+\|g\|_{**}
 \leq C\bigl(\|f\|_{**}+o_b(1)N\bigr).
\end{aligned}
\]
Apply \eqref{eq:ce-free-potential-bound} to
\(\phi=\mathcal T(\Delta^2\phi)\) and absorb
\(o_b(1)(\|\phi\|_*+\|\mathbf c\|_\infty)\).
This proves \eqref{eq:ce-linear-bound} and
\eqref{eq:ce-linear-pointwise}.
The proof also applies uniformly when the multiplier terms and
orthogonality conditions in \eqref{eq:ce-linear-problem} are
restricted to \(0\leq k\leq M\), and
\(U_{\lambda,r}^{\frac{8}{n-4}}\) is replaced by the coefficient
\(U_{\lambda,r}^{\frac{8}{n-4}}(1-\sum_{k>M,i}\chi_{k,i})\) in
\eqref{eq:ce-truncated-coefficient}. Indeed, this replacement preserves symmetry and
decreases the upper bounds used in \eqref{eq:ce-linear-source-moment}
and in the pointwise estimate above. On the core
\(\mathbf B_{2\eta r_m}(\mathbf x_{m,i})\) of an omitted bubble,
\(g=0\); taking \(\psi=0\) there gives the estimate of the form
\eqref{eq:ce-linear-away-bound} also in its surrounding ball.
All constants are thus independent of the truncation index.
\end{proof}

\begin{proposition}\label{lem:ce-inverse}
For small \(b>0\), every measurable \(f\) with the
symmetries in \eqref{eq:ce-linear-space} and finite \(\|f\|_{**}\)
has a unique solution \((\phi,\mathbf c)\) of the linear constrained
problem \eqref{eq:ce-linear-problem} with
\(\phi\in\mathcal H\), \(\|\phi\|_*<\infty\) and
\(\|\mathbf c\|_\infty<\infty\).
This solution satisfies \eqref{eq:ce-linear-bound} and
\eqref{eq:ce-linear-pointwise}.
\end{proposition}
\begin{proof}
The Hardy--Littlewood--Sobolev inequality and
\eqref{eq:ce-linear-norms} give
\begin{equation}\label{eq:ce-potential-energy-bound}
\begin{aligned}
 \|\mathcal Tf\|_{H^2(\Omega)}
 &\leq C\left(\|f\|_{L^{\frac{2n}{n+4}}(\Omega)}
       +\int_\Omega |y|^{2-n}|f(y)|\,dy\right)\\
 &\leq Cb\|f\|_{**}
       \left(1+\sum_k r_k^{\frac{n}{2}+\alpha}+\sum_k r_k^\alpha\right)
 \leq Cb\|f\|_{**}.
\end{aligned}
\end{equation}

In \eqref{eq:ce-linear-problem}, retain the terms containing
\(c_{k,\nu}\) and the orthogonality conditions for
\(0\leq k\leq M\). Replace \(U_{\lambda,r}^{\frac{8}{n-4}}\) in
\(\mathcal L_{\lambda,r}\) by
\begin{equation}\label{eq:ce-truncated-coefficient}
 U_{\lambda,r}^{\frac8{n-4}}
 \left(1-\sum_{k>M,i}\chi_{k,i}\right).
\end{equation}
For fixed \(b>0\), the space of continuous functions on
\(\Omega\setminus\{0\}\) with the symmetries in \eqref{eq:ce-linear-space} and finite
\(\|\cdot\|_*\) is complete. The coefficient in
\eqref{eq:ce-truncated-coefficient} is bounded and tends to zero
at the origin. The operator
\[
 \phi\mapsto\mathcal T\left[
 \frac{n+4}{n-4}U_{\lambda,r}^{\frac8{n-4}}
 \left(1-\sum_{k>M,i}\chi_{k,i}\right)\phi
 \right]
\]
maps \(\{\phi:\|\phi\|_*\leq1\}\) to a uniformly bounded and
equicontinuous family on \(\overline\Omega\), with each function
vanishing at the origin. By the Arzel\`a--Ascoli theorem, the
operator is compact in \(\|\cdot\|_*\).
Adding the terms containing \(c_{k,\nu}\) and imposing the
orthogonality conditions in \eqref{eq:ce-linear-problem} for
\(0\leq k\leq M\) gives a Fredholm operator of index zero.
Lemma~\ref{lem:ce-linear-estimate} gives injectivity and hence
solvability, as in \cite[Proposition~2.2]{WY}.

Let \((\phi_M,\mathbf c^{(M)})\) denote the solution with
\(0\leq k\leq M\). The bounds \eqref{eq:ce-linear-bound} and
\eqref{eq:ce-linear-pointwise} are independent of \(M\).
Local elliptic estimates
and a diagonal subsequence give locally uniform convergence of
\(\phi_M\) and convergence of each \(c_{k,\nu}^{(M)}\).
Passing to the limit gives the equation and every orthogonality
condition in \eqref{eq:ce-linear-problem}, together with
\eqref{eq:ce-linear-bound} and \eqref{eq:ce-linear-pointwise}.
Moreover, \eqref{eq:ce-linear-bound} gives
\begin{equation}\label{eq:ce-truncated-source-majorant}
 |\Delta^2\phi_M(x)|\leq C\|f\|_{**}
 \sum_{k,i}\frac{br_k^{2+\alpha}\lambda_k^\tau}
 {(\lambda_k+|x-\mathbf x_{k,i}|)^{\tau+4}}.
\end{equation}
The sum on the right-hand side of
\eqref{eq:ce-truncated-source-majorant} belongs to
\(L^{\frac{2n}{n+4}}(\Omega)\) and is integrable
against \(|x|^{2-n}\), as in the proof of
\eqref{eq:ce-potential-energy-bound}.
Dominated convergence and the first inequality in
\eqref{eq:ce-potential-energy-bound} therefore give
\[
 \phi_M=\mathcal T(\Delta^2\phi_M)
 \to\mathcal T(\Delta^2\phi)
 \quad\text{in }H^2(\Omega).
\]
Thus \(\phi=\mathcal T(\Delta^2\phi)\in\mathcal H\), and
\eqref{eq:ce-linear-bound} and \eqref{eq:ce-potential-energy-bound} give
\[
 \|\phi\|_{H^2(\Omega)}\leq Cb\|\Delta^2\phi\|_{**}
                        \leq Cb\|f\|_{**}.
\]
Uniqueness follows from Lemma~\ref{lem:ce-linear-estimate} applied
to the linear constrained problem \eqref{eq:ce-linear-problem}.
\end{proof}

For fixed \(\lambda,r\), write \(\mathcal A(f)=\phi\) for the solution
in Proposition~\ref{lem:ce-inverse}. Then \(\mathcal A\) is a linear operator and
\(\|\mathcal A(f)\|_*\leq C\|f\|_{**}\).

\subsection{The nonlinear constrained problem}\label{subsec:ce-nonlinear}
We now solve the nonlinear constrained problem \eqref{eq:ce-constrained-correction}.
\begin{proposition}\label{lem:ce-projected}
Fix \(\varepsilon_0>0\) small enough. For all
small \(b>0\), uniformly for parameters in
\eqref{eq:ce-parameter-space} with \(0\leq\varepsilon\leq\varepsilon_0\),
there is a unique solution of \eqref{eq:ce-constrained-correction} satisfying
\[
 U_{\lambda,r}+\phi>0,
 \,\ \|\phi\|_*+\|\mathbf c\|_\infty\leq4C\varepsilon_0\leq1,
\]
where \(C\) is fixed large enough. It satisfies
\begin{equation}\label{eq:ce-nonlinear-bound}
\begin{aligned}
 \|\phi\|_*+\|\mathbf c\|_\infty
 &\leq C\left(\varepsilon+b^{\frac{3n-2}{(n-3)(n-2)}}\right),\\
 |\phi(x)|&\leq C\left(\varepsilon+b^{\frac{3n-2}{(n-3)(n-2)}}\right)
                         b|x|^{2+\alpha}.
\end{aligned}
\end{equation}
Moreover, \(U_{\lambda,r}+\phi\in C^\infty(\Omega\setminus\{0\})\).
\end{proposition}

We first estimate \(\mathcal N(\phi)\) and \(\mathcal E\).
\begin{lemma}\label{lem:ce-nonlinear-estimates}
Suppose that \(\|\phi_j\|_*\leq1\) and
\(U_{\lambda,r}+\phi_j>0\) for \(j=1,2\). Then
\[
 \|\mathcal N(\phi_1)-\mathcal N(\phi_2)\|_{**}
 \leq o_b(1)\|\phi_1-\phi_2\|_*.
\]
\end{lemma}
\begin{proof}
Elementary inequalities and the fundamental theorem of calculus
give, for every \(n\geq5\),
\begin{equation}\label{eq:ce-nonlinear-pointwise}
|\mathcal N(\phi_1)-\mathcal N(\phi_2)|
 \leq C\left[(|\phi_1|+|\phi_2|)^{\frac8{n-4}} +(1+U_{\lambda,r})^{\frac{\max\{12-n,0\}}{n-4}}
                         (|\phi_1|+|\phi_2|)\right]|\phi_1-\phi_2|.
\end{equation}
Since \(\|\phi_j\|_\infty\leq Cb\) and
\(\|v\|_{**}\leq C\|v\|_*\), we have
\[
 \left\|(|\phi_1|+|\phi_2|)^{\frac8{n-4}}
                  |\phi_1-\phi_2|\right\|_{**}
 \leq Cb^{\frac8{n-4}}\|\phi_1-\phi_2\|_*.
\]
For \(5\leq n<12\), \eqref{eq:ce-profile-upper} gives
\[
 \{x\in\Omega\setminus\{0\}:U_{\lambda,r}(x)>1\}
 \subset\bigcup_{k,i}\mathbf B_{\eta r_k}(\mathbf x_{k,i}).
\]
By symmetry, take \(x\in\mathbf B_{\eta r_k}(\mathbf x_{k,1})\)
with \(U_{\lambda,r}(x)>1\), and write
\(x=\mathbf x_{k,1}+\lambda_k y\).
Then \(U_{\lambda,r}(x)\leq C\mathcal B_{\lambda_k,\mathbf x_{k,1}}(x)\)
and \(1+|y|\leq C\lambda_k^{-\frac{1}{2}}\).
Using \(\|\phi_j\|_\infty\leq Cb\) and
\(\|\phi_1-\phi_2\|_\infty\leq Cb\|\phi_1-\phi_2\|_*\), we obtain
\begin{equation}\label{eq:ce-nonlinear-local-bound}
\begin{aligned}
&\frac{(\lambda_k+|x-\mathbf x_{k,1}|)^{\tau+4}}
        {br_k^{2+\alpha}\lambda_k^\tau}
       (1+U_{\lambda,r}(x))^{\frac{12-n}{n-4}} \times(|\phi_1(x)|+|\phi_2(x)|)
                   |\phi_1(x)-\phi_2(x)|\\
 &\quad\leq\frac{Cb\lambda_k^{\frac{n-4}{2}}}{r_k^{2+\alpha}}
       (1+|y|)^{\tau+n-8}\|\phi_1-\phi_2\|_*\\
 &\quad\leq\frac{Cb}{r_k^{2+\alpha}}
       \left(\lambda_k^{\frac{n-4}{2}}+\lambda_k^{\frac{4-\tau}{2}}\right)
       \|\phi_1-\phi_2\|_* \leq Cb\|\phi_1-\phi_2\|_*.
\end{aligned}
\end{equation}
The last inequality in \eqref{eq:ce-nonlinear-local-bound}
follows from \eqref{eq:ce-parameter-space} and
\(0<\alpha<1+\frac{3}{n}\). For \(n\geq12\) or \(U_{\lambda,r}(x)\leq1\),
\((1+U_{\lambda,r}(x))^{\frac{\max\{12-n,0\}}{n-4}}\) is bounded,
and \(\|v\|_{**}\leq C\|v\|_*\) applies directly.
Combining \eqref{eq:ce-nonlinear-pointwise} and
\eqref{eq:ce-nonlinear-local-bound} with
\(\|\phi_j\|_\infty\leq Cb\) and \(\|v\|_{**}\leq C\|v\|_*\) gives
\[
 \|\mathcal N(\phi_1)-\mathcal N(\phi_2)\|_{**}
 \leq C\left(b+b^{\frac8{n-4}}\right)\|\phi_1-\phi_2\|_*,
\]
which proves the lemma.
\end{proof}

\begin{lemma}\label{lem:ce-residual}
For \((\lambda,r)\) satisfying \eqref{eq:ce-parameter-space},
\[
 \|\mathcal E\|_{**}
 \leq C\left(\varepsilon+b^{\frac{3n-2}{(n-3)(n-2)}}\right).
\]
\end{lemma}
\begin{proof}
The local \(C^4\) convergence of the series defining
\(\mathbf h_{r,\lambda}\) and \(h_{r,\lambda}\) gives
\[
 \Delta^2U_{\lambda,r}
 =\sum_{k=0}^\infty\sum_{i=1}^{\ell}
 \mathcal B_{\lambda_k,\mathbf x_{k,i}}^{\frac{n+4}{n-4}},
 \,\  \text{in }\Omega\setminus\{0\}.
\]
By symmetry, take \(i=1\). Recall from
\eqref{eq:ce-background-decomposition} that
\(U_{\lambda,r}=\mathcal B_{\lambda_k,\mathbf x_{k,1}}+H_k\).
On \(\mathbf B_{\eta r_k}(\mathbf x_{k,1})\), we have
\begin{equation}\label{eq:ce-residual-local-decomposition}
 |\mathcal E|\leq C\mathcal B_{\lambda_k,\mathbf x_{k,1}}^{\frac8{n-4}}|H_k|
       +C|H_k|^{\frac{n+4}{n-4}}
       +\sum_{(m,i)\ne(k,1)}
          \mathcal B_{\lambda_m,\mathbf x_{m,i}}^{\frac{n+4}{n-4}}.
\end{equation}
Taylor's theorem, \eqref{eq:ce-background-taylor} and \(\tau<2\) give
\begin{equation}\label{eq:ce-residual-local-bound}
 \mathcal B_{\lambda_k,\mathbf x_{k,1}}^{\frac8{n-4}}|H_k|
 \leq C\left(\varepsilon+b^{\frac{3n-8}{(n-4)(n-3)}}\right)
 \frac{br_k^{2+\alpha}\lambda_k^\tau}
 {(\lambda_k+|x-\mathbf x_{k,1}|)^{\tau+4}}.
\end{equation}
On \(\mathbf B_{\eta r_k}(\mathbf x_{k,1})\),
\(|H_k|\leq Cb^{\frac{n-4}{n-3}}r_k^2\) gives
\(|H_k|^{\frac{n+4}{n-4}}\leq Cb^{\frac{n+4}{n-3}}\), and
\(\sum_{(m,i)\ne(k,1)}\mathcal B_{\lambda_m,\mathbf x_{m,i}}^{\frac{n+4}{n-4}}
\leq Cb^{\frac{n+4}{n-4}}\).
Outside \(\bigcup_{k,i}\mathbf B_{\eta r_k}(\mathbf x_{k,i})\),
\(|\mathcal E(x)|\leq C(b^{\frac{n+4}{n-3}}+b^{\frac{n+4}{n-4}})\).
Thus \eqref{eq:ce-residual-local-decomposition} and
\eqref{eq:ce-residual-local-bound} give, for \(x\in\Omega\setminus\{0\}\),
\begin{equation}\label{eq:ce-residual-pointwise-bound}
 |\mathcal E(x)|\leq
 C\left(\varepsilon+b^{\frac{3n-8}{(n-4)(n-3)}}\right)
 \sum_{k,i}\frac{br_k^{2+\alpha}\lambda_k^\tau}
 {(\lambda_k+|x-\mathbf x_{k,i}|)^{\tau+4}}
 +C\left(b^{\frac{n+4}{n-3}}+b^{\frac{n+4}{n-4}}\right).
\end{equation}
Combining \eqref{eq:ce-source-weight-lower} and
\eqref{eq:ce-residual-pointwise-bound} with the inequalities
\[
 \frac7{n-3}-\frac4{n-2}=\frac{3n-2}{(n-3)(n-2)}
 <\min\left\{\frac{3n-8}{(n-4)(n-3)},
                   \frac8{n-4}-\frac4{n-2}\right\}
\]
gives the asserted estimate.
\end{proof}

\begin{proof}[Proof of Proposition~\ref{lem:ce-projected}]
We solve
\[
 \phi=\mathcal A(\mathcal E+\mathcal N(\phi))
\]
by contraction in the set of continuous functions on
\(\Omega\setminus\{0\}\) with the symmetries in
\eqref{eq:ce-linear-space} satisfying
\begin{equation}\label{eq:ce-nonlinear-closed-set}
 \|\phi\|_*\leq2C\|\mathcal E\|_{**},
 \,\ |\phi(x)|\leq\tfrac12U_{\lambda,r}(x).
\end{equation}
Here \(C\) is the constant in \eqref{eq:ce-linear-bound} and
\eqref{eq:ce-linear-pointwise}. The set defined by
\eqref{eq:ce-nonlinear-closed-set} is complete under
\(\|\cdot\|_*\), and \(2C\|\mathcal E\|_{**}\leq1\) by
Lemma~\ref{lem:ce-residual}, for \(\varepsilon_0\) and then \(b\)
small enough.

For every \(\phi\) satisfying \eqref{eq:ce-nonlinear-closed-set},
Lemma~\ref{lem:ce-nonlinear-estimates} and \(\mathcal N(0)=0\) give
\[
 \|\mathcal E+\mathcal N(\phi)\|_{**}
 \leq\|\mathcal E\|_{**}+o_b(1)\|\phi\|_*
 \leq2\|\mathcal E\|_{**}.
\]
Since \(\Omega\) is bounded and \(0<b\leq 1\), we have
\(b|x''|^{2+\alpha}\leq Cb^{\frac{n-4}{n-3}}|x''|^2\), so
\eqref{eq:ce-profile-positive} gives
\begin{equation}\label{eq:ce-contraction-profile-lower}
 U_{\lambda,r}(x)\geq C^{-1}b|x|^{2+\alpha}.
\end{equation}
For fixed small \(\varepsilon_0>0\) and \(b>0\) small enough,
Proposition~\ref{lem:ce-inverse}, Lemma~\ref{lem:ce-residual} and
\eqref{eq:ce-contraction-profile-lower} give
\[
\begin{aligned}
 \|\mathcal A(\mathcal E+\mathcal N(\phi))\|_*
 &\leq2C\|\mathcal E\|_{**},\\
 |\mathcal A(\mathcal E+\mathcal N(\phi))(x)|
 &\leq2C\|\mathcal E\|_{**}b|x|^{2+\alpha}
       \leq\tfrac12U_{\lambda,r}(x).
\end{aligned}
\]
Thus \(\phi\mapsto\mathcal A(\mathcal E+\mathcal N(\phi))\)
maps the set defined by \eqref{eq:ce-nonlinear-closed-set} into itself.
For any two functions \(\phi_1,\phi_2\) satisfying
\eqref{eq:ce-nonlinear-closed-set},
\[
 \|\mathcal A(\mathcal E+\mathcal N(\phi_1))
       -\mathcal A(\mathcal E+\mathcal N(\phi_2))\|_*
 \leq Co_b(1)\|\phi_1-\phi_2\|_*
 \leq\tfrac12\|\phi_1-\phi_2\|_*.
\]
The contraction mapping principle gives a unique fixed point.
Proposition~\ref{lem:ce-inverse} gives
\(\phi\in\mathcal H\), \(\phi=\mathcal T(\Delta^2\phi)\), and
the coefficients \(\mathbf c\). Equations~\eqref{eq:ce-linear-bound}
and \eqref{eq:ce-linear-pointwise}, with
\(f=\mathcal E+\mathcal N(\phi)\), give \eqref{eq:ce-nonlinear-bound}
by Lemma~\ref{lem:ce-residual}, since \(Co_b(1)\|\phi\|_*\)
may be absorbed.
Uniqueness follows from Lemma~\ref{lem:ce-nonlinear-estimates}
and \eqref{eq:ce-linear-bound}.
Since \(U_{\lambda,r}+\phi\geq \frac{U_{\lambda,r}}{2}>0\), interior elliptic
regularity and bootstrapping give the asserted smoothness.
\end{proof}

\section{Solving the original equation \eqref{eq:biharmonic-problem}}\label{sec:ce-reduction}
\subsection{Local Pohozaev identities and vanishing of \texorpdfstring{\(c_{k,1},c_{k,2}\)}{c(k,1), c(k,2)}}
We will choose \((\lambda,r)\) satisfying \eqref{eq:ce-parameter-space},
with \(0\leq\varepsilon\leq\varepsilon_0\), so that
\(c_{k,1}=c_{k,2}=0\) for all \(k\geq0\).
First, we use local Pohozaev identities to derive a condition
equivalent to \(c_{k,1}=c_{k,2}=0\).
Let \(u=U_{\lambda,r}+\phi\) be the function given by
Proposition~\ref{lem:ce-projected}. Define
\(\mathcal R_k=(\mathcal R_{k,1},\mathcal R_{k,2})^T\) and
\(\mathcal R=(\mathcal R_k)_{k\geq0}\) by
\begin{equation}\label{eq:ce-reduced-map}
\begin{aligned}
 \mathcal R_{k,1}(\lambda,r)
 &=\frac{2}{(n-4)b\gamma^{k(2+\alpha)}
        \int_{\mathbb R^n}\mathcal B_{\lambda_k,\mathbf{x}_{k,1}}^{\frac{n+4}{n-4}}}
 \times\left\{\begin{aligned}
 &\int\zeta_{k,1}\Delta^2u
     \left[(x-\mathbf x_{k,1})\cdot\nabla T_k+\frac{n-4}{2}T_k\right]\\
 &+\frac{n-4}{2n}\int
       ((x-\mathbf{x}_{k,1})\cdot\nabla\zeta_{k,1})u^{\frac{2n}{n-4}}
 \end{aligned}\right\},\\
 \mathcal R_{k,2}(\lambda,r)
 &=\frac{r_k}{b\gamma^{k(2+\alpha)}
        \int_{\mathbb R^n}\mathcal B_{\lambda_k,\mathbf{x}_{k,1}}^{\frac{n+4}{n-4}}}
 \times\biggl[\int\zeta_{k,1}\Delta^2u\,\partial_{x_1}T_k
    +\frac{n-4}{2n}\int(\partial_{x_1}\zeta_{k,1})u^{\frac{2n}{n-4}}\biggr],
\end{aligned}
\end{equation}
where
\[
 T_k=u-G_n*(\zeta_{k,1}\Delta^2u).
\]
The functions \(\zeta_{k,i}\) are obtained from smooth radial functions,
fixed independently of \(b,\lambda,r\), by translation to
\(\mathbf x_{k,i}\) and scaling by \(r_k\), so that
\(\zeta_{k,i}=1\) on \(\operatorname{supp}\chi_{k,i}\).
The supports \(\operatorname{supp}\zeta_{k,i}\) are slightly larger than
\(\operatorname{supp}\chi_{k,i}\), pairwise disjoint, and have closures
in \(\mathbf B_{4\eta r_k}(\mathbf x_{k,i})\subset\Omega\).
Extend \(\zeta_{k,1}\Delta^2u\) by zero outside \(\Omega\).
The form of \(\mathcal R_k\) follows from the local Pohozaev
identities derived in the proof of the following lemma.

\begin{lemma}
\label{lem:ce-local-identities}
For small \(b>0\), uniformly for \((\lambda,r)\) satisfying
\eqref{eq:ce-parameter-space} and all \(k\geq0\), the equation
\begin{equation}\label{eq:ce-parameter-equations}
 \mathcal R_k(\lambda,r)=0,\,\  k\geq0,
\end{equation}
is equivalent to \(c_{k,1}=c_{k,2}=0\).
\end{lemma}
\begin{proof}
Fix \(k\geq0\) and set \(f=\zeta_{k,1}\Delta^2u\). Antisymmetry of
\(DG_n\) and homogeneity of \(G_n\) of degree \(4-n\) give
\begin{equation}\label{eq:ce-potential-identities}
\int f\,\partial_{x_1}(G_n*f)=0,\,\   \int f\left[(x-\mathbf x_{k,1})\cdot\nabla(G_n*f)
                   +\frac{n-4}{2}(G_n*f)\right]=0.
\end{equation}
The identities in \eqref{eq:ce-potential-identities} follow by
interchanging the variables in the double integrals,
which are absolutely convergent because \(f\) has compact support
and \(G_n,DG_n\) are locally integrable.
Using \(u=G_n*f+T_k\) and
\eqref{eq:ce-constrained-correction}, integration by parts gives
\begin{equation}\label{eq:ce-local-pohozaev}
\begin{aligned}
&\frac{n+4}{n-4}br_k^{2+\alpha}\sum_{\nu=1}^2c_{k,\nu} \times\int\chi_{k,1}\mathcal B_{\lambda_k,\mathbf x_{k,1}}^{\frac8{n-4}}
        Z_{k,1,\nu}\left[(x-\mathbf x_{k,1})\cdot\nabla u+\frac{n-4}{2}u\right]\\
 &\,\ =\int f\left[(x-\mathbf x_{k,1})\cdot\nabla T_k+\frac{n-4}{2}T_k\right] +\frac{n-4}{2n}\int
          ((x-\mathbf x_{k,1})\cdot\nabla\zeta_{k,1})u^{\frac{2n}{n-4}},\\
 &\frac{n+4}{n-4}br_k^{2+\alpha}\sum_{\nu=1}^2c_{k,\nu}
   \int\chi_{k,1}\mathcal B_{\lambda_k,\mathbf x_{k,1}}^{\frac8{n-4}}
        Z_{k,1,\nu}\partial_{x_1}u\\
 &\,\ =\int f\partial_{x_1}T_k
   +\frac{n-4}{2n}\int(\partial_{x_1}\zeta_{k,1})u^{\frac{2n}{n-4}}.
\end{aligned}
\end{equation}
Equations~\eqref{eq:ce-local-pohozaev} are the local Pohozaev
identities for dilation and radial translation. The terms containing \(c_{j,\nu}\chi_{j,i}\) with
\((j,i)\ne(k,1)\) vanish because
\(\operatorname{supp}\chi_{j,i}\cap\operatorname{supp}\zeta_{k,1}
=\varnothing\), and \(\zeta_{k,1}=1\) on the support of
\(\chi_{k,1}\). Define the coefficient matrix \(M_k\) by
\begin{equation}\label{eq:ce-coefficient-matrix}
\begin{aligned}
(M_k)_{1\nu}&=\frac{n+4}{n-4}br_k^{2+\alpha} \times\int\chi_{k,1}\mathcal B_{\lambda_k,\mathbf x_{k,1}}^{\frac8{n-4}}
        Z_{k,1,\nu}\left[(x-\mathbf x_{k,1})\cdot\nabla u+\frac{n-4}{2}u\right],\\
 (M_k)_{2\nu}&=\frac{n+4}{n-4}br_k^{2+\alpha}
   \int\chi_{k,1}\mathcal B_{\lambda_k,\mathbf x_{k,1}}^{\frac8{n-4}}
        Z_{k,1,\nu}\partial_{x_1}u,
 \,\  \nu=1,2.
\end{aligned}
\end{equation}
Equations~\eqref{eq:ce-reduced-map}, \eqref{eq:ce-local-pohozaev}
and \eqref{eq:ce-coefficient-matrix} give
\begin{equation}\label{eq:ce-coefficient-system}
 \mathcal R_k=0
 \quad\Longleftrightarrow\quad
 M_k\begin{pmatrix}c_{k,1}\\c_{k,2}\end{pmatrix}=0.
\end{equation}
It remains to prove that \(M_k\) is invertible.

Use \(u=\mathcal B_{\lambda_k,\mathbf{x}_{k,1}}+H_k+\phi\).
On \(\operatorname{supp}\chi_{k,1}\), \eqref{eq:ce-background-taylor} and
\eqref{eq:ce-nonlinear-bound} give
\[
 |H_k+\phi|
 \leq C\bigl(b^{\frac{n-4}{n-3}}r_k^2+br_k^{2+\alpha}\bigr).
\]
Integrating by parts transfers the derivatives of \(H_k+\phi\)
to \(\chi_{k,1}\mathcal B_{\lambda_k,\mathbf x_{k,1}}^{\frac{8}{n-4}}
Z_{k,1,\nu}\), for \(\nu=1,2\). By scaling, the contribution of
\(H_k+\phi\) to \((M_k)_{1\nu}\) is
\(br_k^{2+\alpha}\lambda_k^{\frac{n-4}{2}}o_b(1)\).
For \((M_k)_{2\nu}\), it is
\(br_k^{2+\alpha}\lambda_k^{\frac{n-6}{2}}o_b(1)\).
Here we use
\[
 C\bigl(b^{\frac{n-4}{n-3}}r_k^2+br_k^{2+\alpha}\bigr)\lambda_k^{\frac{n-4}{2}}
 =o_b(1),
\]
since \(\lambda_k^{\frac{n-4}{2}}\leq Cbr_k^{n-2+\alpha}\) and
\(r_k\leq e^{\varepsilon_0}\).
By symmetry, the contribution of \(\mathcal B_{\lambda_k,\mathbf x_{k,1}}\)
to \((M_k)_{12}\) and \((M_k)_{21}\) is zero.
Since \(\frac{\lambda_k}{r_k}\to0\) uniformly, scaling in
\eqref{eq:ce-coefficient-matrix} gives
\begin{equation}\label{eq:ce-coefficient-determinant}
\det M_k
 =\left(\frac{n+4}{n-4}\right)^2
     b^2r_k^{4+2\alpha}\lambda_k^{n-5} \times\prod_{\nu=1}^2
 \left(\int_{\mathbb R^n}\mathcal B_{1,0}^{\frac8{n-4}}Z_\nu^2\right)
 (1+o_b(1)).
\end{equation}
Both integrals in \eqref{eq:ce-coefficient-determinant} are positive,
so \(\det M_k>0\) for small \(b\), uniformly in \(k,\lambda,r\).
Thus \(M_k\) is invertible, and \eqref{eq:ce-coefficient-system} proves
\(\mathcal R_k=0\) if and only if \(c_{k,1}=c_{k,2}=0\).
\end{proof}

It remains to choose \(\lambda,r\) so that
\eqref{eq:ce-parameter-equations} holds for every \(k\geq0\).
In Subsection~\ref{subsec:ce-linearization}, we consider the linearization of
\(\mathcal R(\lambda,r)\) as \(b\to0\) at \(\lambda_k=\lambda_k^0\),
\(r_k=\gamma^k\), \(k\geq-3\), and define the linear operator \(L\)
and the remainder \(E_b\).
In Subsection~\ref{app:parameter-surjectivity}, we construct a
bounded right inverse of \(L\), allowing
\(\lambda_{-j},r_{-j}\), \(j=1,2,3\), to vary.
Subsection~\ref{subsec:ce-regular-parts} gives the estimates,
continuity and compactness of \(E_b\).
Finally, in Subsection~\ref{subsec:ce-parameter-choice}, we apply the
Schauder fixed point theorem to choose \(\lambda,r\) with
\(\mathcal R(\lambda,r)=0\).
This gives \(c_{k,1}=c_{k,2}=0\) for every \(k\geq0\).

\subsection{The linear operator \texorpdfstring{\(L\)}{L} and the remainder \texorpdfstring{\(E_b\)}{E(b)}}\label{subsec:ce-linearization}
We first determine \(\lim_{b\to0}\mathcal R(\lambda,r)\).
The bounds for \(\alpha\) in Lemma~\ref{lem:ce-certificates} give
\[
 2-\alpha\geq1-\frac3n>0,
 \,\ n-(n-6)\alpha\geq\min\left\{n,3+\frac{18}n\right\}>0.
\]
Set
\[
 \beta=\min\left\{2-\alpha,\frac{n-(n-6)\alpha}{n-4}\right\}>0.
\]
For the parameter equations from \eqref{eq:ce-reduced-parameter-equation},
we further require \(\gamma^\beta<\rho<1\).

We first estimate \(T_k\) and the local integrals of \(\Delta^2u\) in
\eqref{eq:ce-reduced-map}. Write
\[
\begin{aligned}
 T_k
 &=U_{\lambda,r}-G_n*(\zeta_{k,1}\Delta^2U_{\lambda,r})
       +T_{\phi,k},\\
 T_{\phi,k}
 &=\phi-G_n*(\zeta_{k,1}\Delta^2\phi).
\end{aligned}
\]
\begin{lemma}\label{lem:ce-regular-source}
For \(j=0,1,2\),
\begin{equation}\label{eq:ce-regular-correction}
 \sup_{x\in\operatorname{supp}\zeta_{k,1}}|D^jT_{\phi,k}(x)|
 \leq b\,o_b(1)r_k^{4-j}(1+|\log r_k|).
\end{equation}
On the support of \(\zeta_{k,1}\), \(T_k\) satisfies
\begin{equation}\label{eq:ce-regular-derivatives}
 \sup|DT_k|\leq Cb^{\frac{n-4}{n-3}}r_k,
 \,\ \sup|D^2T_k|\leq Cb^{\frac{n-4}{n-3}}.
\end{equation}
The local integrals of \(\Delta^2u\) satisfy
\begin{equation}\label{eq:ce-local-source-integrals}
 \int\zeta_{k,1}|x-\mathbf{x}_{k,1}|^j|\Delta^2u|
 \leq C\left[\lambda_k^{\frac{n-4}{2}+j}
       +b^{\frac{n+4}{n-3}}r_k^{n+\frac{2(n+4)}{n-4}+j}\right],
 \,\  j=0,1.
\end{equation}
Moreover,
\begin{equation}\label{eq:ce-relative-source-integral}
 \frac{\left|\int\zeta_{k,1}\Delta^2u
       -\int_{\mathbb R^n}\mathcal B_{\lambda_k,\mathbf{x}_{k,1}}^{\frac{n+4}{n-4}}\right|}
      {\int_{\mathbb R^n}\mathcal B_{\lambda_k,\mathbf{x}_{k,1}}^{\frac{n+4}{n-4}}}
 \leq o_b(1)r_k^{2-\alpha}.
\end{equation}
All bounds are uniform in \(k\geq0\).
\end{lemma}
\begin{proof}
Write \(p=\frac{n+4}{n-4}\), \(\sigma=\frac{n-4}{n-3}\), and abbreviate
\(B_{k,i}=\mathcal B_{\lambda_k,\mathbf x_{k,i}}\). We first construct a
decomposition
\begin{equation}\label{eq:ce-correction-source-decomposition}
 \Delta^2\phi=f_0+\sum_{k,i}f_{k,i},
 \,\ |f_0(x)|\leq C\left(b^{\frac{n+4}{n-3}}
 |x|^{\frac{2(n+4)}{n-4}}+b^{\frac{n+4}{n-4}}\right),
\end{equation}
where \(f_{k,i}\) is supported in \(\{\zeta_{k,i}=1\}\), \(f_0\) and each
\(F_k:=\sum_i f_{k,i}\) have the symmetries in \eqref{eq:ce-linear-space}, and
\begin{equation}\label{eq:ce-local-source-bound}
 |f_{k,i}(\mathbf x_{k,i}+\lambda_k y)|
 \leq \frac{Cbr_k^{2+\alpha}}{\lambda_k^4(1+|y|)^6}.
\end{equation}
In particular, \(f_0\) and \(F_k\) are even under \(x\mapsto-x\), since
\(\ell\) is even.

Choose \(0\leq\theta_{k,i}\leq1\) from one fixed smooth radial cutoff,
equal to one on \(\mathbf B_{\eta r_k}(\mathbf x_{k,i})\) and supported in
\(\mathbf B_{2\eta r_k}(\mathbf x_{k,i})\subset\{\chi_{k,i}=1\}\). Set
\[
\begin{aligned}
 w_{k,i}&=u-B_{k,i},&
 D_{k,i}&=u^p-B_{k,i}^p,&
 \omega_{k,i}&=\left[1+\left(\frac{|w_{k,i}|}{B_{k,i}}\right)^{p-1}\right]^{-1},\\
 P_{k,i}&=p\sum_{\nu=1}^2c_{k,\nu}br_k^{2+\alpha}\chi_{k,i}B_{k,i}^{p-1}Z_{k,i,\nu},\\
 f_{k,i}&=\theta_{k,i}\omega_{k,i}D_{k,i}+P_{k,i},&
 f_0&=\Delta^2\phi-\sum_{k,i}f_{k,i}.
\end{aligned}
\]
The constrained equation gives
\(\Delta^2\phi=u^p-\sum_{k,i}B_{k,i}^p+\sum_{k,i}P_{k,i}\).
For every \(p>1\),
\[
 |\omega_{k,i}D_{k,i}|\leq CB_{k,i}^{p-1}|w_{k,i}|,
 \,\ |(1-\omega_{k,i})D_{k,i}|\leq C|w_{k,i}|^p.
\]
Indeed, \(|D_{k,i}|\leq C(B_{k,i}^{p-1}|w_{k,i}|+|w_{k,i}|^p)\), and
multiplication by \(\omega_{k,i}\) or \(1-\omega_{k,i}\) reduces the expression in
parentheses to its first or second term. Consequently,
\[
 f_0=\left(1-\sum_{k,i}\theta_{k,i}\right)u^p
 -\sum_{k,i}(1-\theta_{k,i})B_{k,i}^p
 +\sum_{k,i}\theta_{k,i}(1-\omega_{k,i})D_{k,i}.
\]
The first term is supported outside the balls of radius \(\eta r_k\), where
\eqref{eq:ce-profile-upper} and \eqref{eq:ce-nonlinear-bound} give
\(|u|\leq Cb^\sigma|x|^2\). The second term is bounded by \(Cb^p\), by the
bubble tail estimates. On \(\operatorname{supp}\theta_{k,i}\),
\eqref{eq:ce-background-taylor} and \eqref{eq:ce-nonlinear-bound} give
\(|w_{k,i}|\leq Cb^\sigma r_k^2\), and
\(C^{-1}r_k\leq |x|\leq Cr_k\). The supports are disjoint, so these bounds
prove \eqref{eq:ce-correction-source-decomposition}. All the cutoffs use the
same radial profiles, and each symmetry in \eqref{eq:ce-linear-space}
permutes \(B_{k,i},w_{k,i},\omega_{k,i}\) and \(P_{k,i}\) within the same circle;
the coefficients \(c_{k,\nu}\) are independent of \(i\). Thus this explicit
decomposition has the claimed symmetries.

On \(\operatorname{supp}\theta_{k,i}\), the same Taylor estimates give
\[
 |w_{k,i}(\mathbf x_{k,i}+\lambda_k y)|\leq Cbr_k^{2+\alpha}(1+|y|)^2.
\]
For the quadratic term, use
\[
 \frac{b^{\sigma-1}\lambda_k^2}{r_k^{2+\alpha}}
 \leq Cb^{\frac{3n-8}{(n-4)(n-3)}}
 r_k^{\frac{2n+(8-n)\alpha}{n-4}}\leq C;
\]
the exponent of \(r_k\) is positive for \(0<\alpha<2\). Multiplying by
\(B_{k,i}^{p-1}\leq C\lambda_k^{-4}(1+|y|)^{-8}\), and using
\(|Z_{k,i,\nu}|\leq C(1+|y|)^{4-n}\) for \(P_{k,i}\), proves
\eqref{eq:ce-local-source-bound}. Since \(\tau<\min\{2,n-4\}\), we may weaken its
decay to \((1+|y|)^{-\tau-4}\). Integration yields
\(\int|f_{k,i}|\leq Cbr_k^{n-2+\alpha-\tau}\lambda_k^\tau\), hence
\begin{equation}\label{eq:ce-local-source-mass}
 \int|f_{k,i}|\leq Cb^{\frac{n+2}{n-2}}
 r_k^{n+\alpha+\frac{4(1+\alpha)}{n-2}},
 \,\ \int_\Omega|x|^{2-n}|\Delta^2\phi|<\infty.
\end{equation}
The weighted integrability follows by summing these masses, since
\begin{equation}\label{eq:ce-exponent-margins}
 0<\alpha<2,\,\  \mu:=\alpha+\frac{4(1+\alpha)}{n-2}>0.
\end{equation}
We now prove \eqref{eq:ce-regular-correction} on the entire support of
\(\zeta_{k,1}\). The reference geometry and small \(\varepsilon_0\) give
\[
 |\mathbf x_{m,i}-\mathbf x_{k,1}|\geq d_0\max\{r_m,r_k\},
 \quad (m,i)\ne(k,1),
 \,\ d_0=\min\left\{1-\frac{\gamma(1+\varepsilon_0)}{1-\varepsilon_0},
 2\sin\frac{\pi}{\ell}\right\}>0.
\]
We choose \(\eta\) smaller if necessary so that \(8\eta<\frac{d_0}{2}\).

Then
\(\operatorname{supp}\zeta_{k,1}\) is at distance at least
\((\frac{d_0}{2})\max\{r_k,r_m\}\) from \(\operatorname{supp}f_{m,i}\) whenever
\((m,i)\ne(k,1)\). Since \(\zeta_{k,1}f_{k,1}=f_{k,1}\), the exact identity is
\[
 T_{\phi,k}=\mathcal T f_0-G_n*(\zeta_{k,1}f_0)
 +\sum_{m\geq0}\mathcal T F_m-G_n*f_{k,1}.
\]
Put \(M_m=\sum_i\int|f_{m,i}|\leq Cb^{\frac{n+2}{n-2}}r_m^{n+\mu}\).
For \(x\in\operatorname{supp}\zeta_{k,1}\), \(j=0,1,2\), and a fixed
\(C_0\) large enough, Taylor expansion of the normalized kernel gives
\[
 |D^j\mathcal T F_m(x)|\leq CM_m
 \begin{cases}
 r_k^{4-j}r_m^{-n},&r_m\geq C_0r_k,\\
 r_k^{2-j}r_m^{2-n},&r_m\leq \frac{r_k}{C_0}.
 \end{cases}
\]
In the first case the cubic Taylor term integrates to zero because \(F_m\)
is even. In the second case the bound includes every term subtracted by
\(J_2\). For \(\frac{r_k}{C_0}<r_m<C_0r_k\), the preceding separation bounds the
other sources by \(Cr_k^{4-n-j}M_m\). The singular convolution of
\(f_{k,1}\) cancels exactly, leaving \(-J_2[G_n*f_{k,1}]\), which has the same
bound. Summing the three ranges gives
\[
\begin{aligned}
 \left|D^j\left(\sum_m\mathcal T F_m-G_n*f_{k,1}\right)(x)\right|
 &\leq Cb^{\frac{n+2}{n-2}}\left[
 r_k^{4-j}\sum_{r_m\geq C_0r_k}r_m^\mu
 +r_k^{2-j}\sum_{r_m\leq \frac{r_k}{C_0}}r_m^{2+\mu}+r_k^{4+\mu-j}
 \right]\\
 &\leq Cb^{\frac{n+2}{n-2}}r_k^{4-j}.
\end{aligned}
\]
Here the comparable range contains a uniformly bounded number of circles.

Finally, \eqref{eq:ce-correction-source-decomposition} gives
\(\|f_0\|_\infty\leq A_b:=C(b^{\frac{n+4}{n-3}}+b^{\frac{n+4}{n-4}})\).
For \(s=|x|\) with \(C^{-1}r_k\leq s\leq Cr_k\), split \(\mathcal T f_0\)
into \(|y|\leq2s\) and \(|y|>2s\). Integrating \(|D^jG_n(x-y)|\) and the
Taylor terms on the first region gives \(CA_b s^{4-j}\). On the second,
evenness of \(f_0\) cancels the cubic term, leaving
\(CA_b s^{4-j}(1+|\log s|)\), from the radial integral of \(|y|^{-n}\). Also,
\(\sup_{\operatorname{supp}\zeta_{k,1}}
 |D^j G_n*(\zeta_{k,1}f_0)|\leq CA_b r_k^{4-j}\), by direct kernel integration
over a ball of radius \(Cr_k\). Together these estimates yield
\[
 \sup_{\operatorname{supp}\zeta_{k,1}}|D^jT_{\phi,k}|
 \leq Cb\left(b^{\frac4{n-2}}+b^{\frac7{n-3}}+b^{\frac8{n-4}}\right)
 r_k^{4-j}(1+|\log r_k|),\,\  j=0,1,2,
\]
which proves \eqref{eq:ce-regular-correction} with all constants uniform in
\(k,\lambda,r\).

The integral representation of
\(\mathcal B_{\lambda_k,\mathbf x_{k,1}}\) on \(\mathbb R^n\) gives
\[
T_k={}H_k+T_{\phi,k} +G_n*\left[(1-\zeta_{k,1})
       \mathcal B_{\lambda_k,\mathbf{x}_{k,1}}^{\frac{n+4}{n-4}}
       -\zeta_{k,1}\sum_{(m,i)\ne(k,1)}
       \mathcal B_{\lambda_m,\mathbf{x}_{m,i}}^{\frac{n+4}{n-4}}\right].
\]
Use \eqref{eq:ce-background-taylor} for
\(DH_k\) and \(D^2H_k\). For \(j=1,2\),
\(D^jG_n*[(1-\zeta_{k,1})
\mathcal B_{\lambda_k,\mathbf x_{k,1}}^{\frac{n+4}{n-4}}]\)
is bounded by \(C\lambda_k^{\frac{n+4}{2}}r_k^{-n-j}\), whereas
\(D^jG_n*[\zeta_{k,1}\sum_{(m,i)\ne(k,1)}
\mathcal B_{\lambda_m,\mathbf x_{m,i}}^{\frac{n+4}{n-4}}]\)
is bounded by \(Cb^{\frac{n+4}{n-4}}r_k^{4-j}\). Together with the
\(j=1,2\) estimates in \eqref{eq:ce-regular-correction}, these give
\eqref{eq:ce-regular-derivatives}.
The equation \eqref{eq:ce-constrained-correction} also gives,
on \(\operatorname{supp}\zeta_{k,1}\),
\[
 |\Delta^2u|\leq C\left(
 \mathcal B_{\lambda_k,\mathbf{x}_{k,1}}^{\frac{n+4}{n-4}}
 +b^{\frac{n+4}{n-3}}r_k^{\frac{2(n+4)}{n-4}}\right).
\]
Equation~\eqref{eq:ce-local-source-integrals} follows by the change of variables
\(x=\mathbf{x}_{k,1}+\lambda_k y\). For \(j=1\), the rescaled radial
integrand is \(O(s^{-4})\) as \(s\to\infty\), integrable for every
\(n\geq5\). The estimate does not require \(\Delta^2u\geq0\).

Since \(\Delta^2u=\Delta^2\phi+
\sum_{m,i}\mathcal B_{\lambda_m,\mathbf{x}_{m,i}}^{\frac{n+4}{n-4}}\),
\begin{equation}\label{eq:ce-local-mass-difference}
\begin{aligned}
 \int\zeta_{k,1}\Delta^2u
       -\int_{\mathbb R^n}\mathcal B_{\lambda_k,\mathbf{x}_{k,1}}^{\frac{n+4}{n-4}}
 &=\int\zeta_{k,1}\Delta^2\phi\\
 &\quad+\int\zeta_{k,1}\sum_{(m,i)\ne(k,1)}
              \mathcal B_{\lambda_m,\mathbf{x}_{m,i}}^{\frac{n+4}{n-4}}
 -\int_{\mathbb R^n}(1-\zeta_{k,1})
              \mathcal B_{\lambda_k,\mathbf{x}_{k,1}}^{\frac{n+4}{n-4}}.
\end{aligned}
\end{equation}
For the three integrals on the right-hand side of
\eqref{eq:ce-local-mass-difference}, use
\eqref{eq:ce-correction-source-decomposition} and
\eqref{eq:ce-local-source-bound} in the first, the bound
\(Cb^{\frac{n+4}{n-4}}r_k^n\) for the second, and
\(C\lambda_k^{\frac{n-4}{2}}\left(\frac{\lambda_k}{r_k}\right)^4\) for the third.
Divide by \(\int_{\mathbb R^n}\mathcal B_{\lambda_k,\mathbf{x}_{k,1}}^{\frac{n+4}{n-4}}\).
By \eqref{eq:ce-exponent-margins}, the resulting powers of \(r_k\)
are at least \(2-\alpha\), and the absolute value of the quotient is at most
\[
 C\left(b^{\frac4{n-2}}+b^{\frac7{n-3}}+b^{\frac8{n-4}}\right)
 r_k^{2-\alpha}
 \leq o_b(1)r_k^{2-\alpha}.
\]
The exponents \(\frac{4}{n-2}\), \(\frac{7}{n-3}\) and \(\frac{8}{n-4}\)
are positive for \(n\geq5\).
This proves \eqref{eq:ce-relative-source-integral}.
\end{proof}

For \(k\geq-3\), write
\[
 \mathbf p_k=\lambda_k-\lambda_k^0,
 \,\ \mathbf r_k=r_k-\gamma^k.
\]
We write \(\mathbf p=(\mathbf p_k)_{k\geq-3}\) and
\(\mathbf r=(\mathbf r_k)_{k\geq-3}\).
We use the norm
\begin{equation}\label{eq:ce-parameter-increment-norm}
 \|(\mathbf p,\mathbf r)\|_Y
 :=\sup_{k\geq-3}\rho^{-k}
 \left(\frac{n-4}{2}\frac{|\mathbf p_k|}{\lambda_k^0}
                   +\frac{|\mathbf r_k|}{\gamma^k}\right)
\end{equation}
on the space
\[
 Y=\left\{(\mathbf p,\mathbf r):
 \rho^{-k}\left(\frac{\mathbf p_k}{\lambda_k^0},
               \frac{\mathbf r_k}{\gamma^k}\right)\to(0,0)
                   \text{ as }k\to\infty\right\}.
\]
For small \(\varepsilon_0\), the mean value theorem shows that
\eqref{eq:ce-parameter-space} implies
\(\|(\mathbf p,\mathbf r)\|_Y\leq C\varepsilon\), whenever
\((\mathbf p,\mathbf r)\in Y\).
Conversely, a small \(\|(\mathbf p,\mathbf r)\|_Y\) gives
\eqref{eq:ce-parameter-space} with
\(\varepsilon=C\|(\mathbf p,\mathbf r)\|_Y\).
The constant \(C\) is independent of \(b\).
Let
\[
 \mathcal X_\rho=\{f=(f_k)_{k\geq0}:\rho^{-k}f_k\to0\},
 \,\ \|f\|_{\mathcal X_\rho}:=\sup_{k\geq0}\rho^{-k}
 \bigl(|f_{k,1}|+|f_{k,2}|\bigr).
\]
The components \(f_{k,1},f_{k,2}\) correspond to
\(\mathcal R_{k,1},\mathcal R_{k,2}\).

{\emergencystretch=2em
Equation~\eqref{eq:ce-relative-source-integral} bounds the difference between
\(\int\zeta_{k,1}\Delta^2u\) and
\(\int_{\mathbb R^n}\mathcal B_{\lambda_k,\mathbf{x}_{k,1}}^{\frac{n+4}{n-4}}\).

Evaluating \(\left[(x-\mathbf x_{k,1})\cdot\nabla T_k+\frac{n-4}{2}T_k\right]\) and
\(\partial_{x_1}T_k\) at \(\mathbf{x}_{k,1}\) in the first
integrals in \eqref{eq:ce-reduced-map} leaves the difference in
the numerator of \eqref{eq:ce-relative-source-integral}
multiplied by \(T_k(\mathbf{x}_{k,1})\) and
\(r_k\partial_{x_1}T_k(\mathbf{x}_{k,1})\), respectively.
At \(\mathbf x_{k,1}\), both \(b^{\frac{n-4}{n-3}}|x''|^2\) and its
\(x_1\)-derivative vanish. Thus \eqref{eq:ce-background-taylor} and
Lemma~\ref{lem:ce-regular-source} give
\par}
\begin{equation}\label{eq:ce-reduced-evaluation-bound}
 |T_k(\mathbf{x}_{k,1})|
 +r_k|\partial_{x_1}T_k(\mathbf{x}_{k,1})|
 \leq C\|(\mathbf p,\mathbf r)\|_Ybr_k^{2+\alpha}
       +b\,o_b(1)r_k^4(1+|\log r_k|).
\end{equation}
For the contribution of \(f_{k,i}\) to
\eqref{eq:ce-relative-source-integral}, the exponent of \(r_k\) satisfies
\[
 \frac{2(n+2\alpha)}{n-2}-(2-\alpha)
 =\alpha+\frac{4(1+\alpha)}{n-2}>0.
\]
For \(f_0\) in \eqref{eq:ce-correction-source-decomposition} and the two
bubble integrals in \eqref{eq:ce-local-mass-difference}, the powers
of \(r_k\) after division by
\(\int_{\mathbb R^n}\mathcal B_{\lambda_k,\mathbf x_{k,1}}^{\frac{n+4}{n-4}}\)
are at least \(2-\alpha\).
Multiplying \eqref{eq:ce-reduced-evaluation-bound} by the
error in \eqref{eq:ce-relative-source-integral} and dividing by
\(br_k^{2+\alpha}\) gives the bound
\begin{equation}\label{eq:ce-reduced-mass-error}
 o_b(1)(1+\|(\mathbf p,\mathbf r)\|_Y)r_k^{2-\alpha}(1+|\log r_k|).
\end{equation}

Equation~\eqref{eq:ce-regular-derivatives} and the \(j=1\) estimate in
\eqref{eq:ce-local-source-integrals} give the Taylor error in
\(\mathcal R_{k,1}\) and \(\mathcal R_{k,2}\), with the coefficients
in \eqref{eq:ce-reduced-map},
\begin{equation}\label{eq:ce-reduced-taylor-error}
 \frac{Cb^{\frac{n-4}{n-3}}r_k\lambda_k}{br_k^{2+\alpha}}
 \leq Cb^{\frac{n-2}{(n-4)(n-3)}}
          r_k^{\frac{n-(n-6)\alpha}{n-4}}.
\end{equation}
On the support of \(\nabla\zeta_{k,1}\),
\(u\leq Cb^{\frac{n-4}{n-3}}r_k^2\).
The integrals containing \(\nabla\zeta_{k,1}\) in
\eqref{eq:ce-reduced-map}, and the second term in
\eqref{eq:ce-local-source-integrals} with \(j=1\), contribute at most
\begin{equation}\label{eq:ce-reduced-cutoff-error}
 Cb^{\frac6{n-3}}r_k^{\frac{4n}{n-4}-2\alpha}.
\end{equation}
Since \(\frac{4n}{n-4}-2\alpha>2-\alpha\), the bound
\eqref{eq:ce-reduced-cutoff-error} is at most
\(o_b(1)r_k^{2-\alpha}\).

For \(T_{\phi,k}\), the \(j=0,1\) bounds in
\eqref{eq:ce-regular-correction} give
\begin{equation}\label{eq:ce-reduced-correction-error}
 \frac{|T_{\phi,k}(\mathbf{x}_{k,1})|
       +r_k|\partial_{x_1}T_{\phi,k}(\mathbf{x}_{k,1})|}
      {br_k^{2+\alpha}}
 \leq o_b(1)r_k^{2-\alpha}(1+|\log r_k|).
\end{equation}
The terms \(b^{\frac{n}{n-4}}r_k^4\) and \(b^{\frac{n}{n-4}}r_k^3\) in
\eqref{eq:ce-background-taylor} contribute at most
\(Cb^{\frac{4}{n-4}}r_k^{2-\alpha}\) after division by
\(br_k^{2+\alpha}\), with an additional factor \(r_k\) for the
radial derivative. The exponents \(2-\alpha\) and
\(\frac{n-(n-6)\alpha}{n-4}\) are at least \(\beta\).

Using \eqref{eq:ce-reduced-mass-error}--\eqref{eq:ce-reduced-correction-error}
in \eqref{eq:ce-reduced-map} and letting \(b\to0\), we obtain
\begin{equation}\label{eq:ce-point-parameter-sum}
\begin{aligned}
 \lim_{b\to0}\mathcal R_{k,1}(\lambda,r)
 &=\gamma^{-k(2+\alpha)}\Biggl[
  \sum_{j\geq-3}\gamma^{j(n-2+\alpha)}
  \left(\frac{\lambda_j}{\lambda_j^0}\right)^{\frac{n-4}{2}}\Phi_{r_j}(x)\\
 &\,\ \,\ {}
  +\sum_{j\leq-4}\gamma^{j(n-2+\alpha)}\Phi_{\gamma^j}(x)
  -\frac{\gamma^{k(n-2+\alpha)}}{\ell}
  \left(\frac{\lambda_k}{\lambda_k^0}\right)^{\frac{n-4}{2}}
  \mathcal K_{r_k\xi_1}(x)\Biggr]_{x=r_k\xi_1},\\[6pt]
 \lim_{b\to0}\mathcal R_{k,2}(\lambda,r)
 &=r_k\gamma^{-k(2+\alpha)}\partial_{x_1}\Biggl[
  \sum_{j\geq-3}\gamma^{j(n-2+\alpha)}
  \left(\frac{\lambda_j}{\lambda_j^0}\right)^{\frac{n-4}{2}}\Phi_{r_j}(x)\\
 &\,\ \,\ {}
  +\sum_{j\leq-4}\gamma^{j(n-2+\alpha)}\Phi_{\gamma^j}(x)
  -\frac{\gamma^{k(n-2+\alpha)}}{\ell}
  \left(\frac{\lambda_k}{\lambda_k^0}\right)^{\frac{n-4}{2}}
  \mathcal K_{r_k\xi_1}(x)\Biggr]\bigg|_{x=r_k\xi_1},
\end{aligned}
\end{equation}
with \(\frac{\lambda_j}{\lambda_j^0}\) and \(\frac{r_j}{\gamma^j}\) fixed.
For the terms with \(j=k\) in \eqref{eq:ce-point-parameter-sum}, we use
\begin{equation}\label{eq:ce-self-circle-regular-part}
 \Phi_{r_k}(x)-\frac1\ell\mathcal K_{r_k\xi_1}(x)
 =\frac1\ell\sum_{i=2}^{\ell}\mathcal K_{r_k\xi_i}(x)
  -J_2\left[\frac1\ell\sum_{i=1}^{\ell}\mathcal K_{r_k\xi_i}\right](x).
\end{equation}
The right-hand side of \eqref{eq:ce-self-circle-regular-part} is smooth
near \(r_k\xi_1\). In \eqref{eq:ce-point-parameter-sum}, the value
and the \(x_1\)-derivative of
\(\Phi_{r_k}-\frac{\mathcal K_{r_k\xi_1}}{\ell}\) at \(x=r_k\xi_1\)
are obtained from the right-hand side of
\eqref{eq:ce-self-circle-regular-part}.

Since \(r_k\) is uniformly comparable to \(\gamma^k\), combining
\eqref{eq:ce-reduced-mass-error}, \eqref{eq:ce-reduced-taylor-error},
\eqref{eq:ce-reduced-cutoff-error}, \eqref{eq:ce-reduced-correction-error}
and the bound \(Cb^{\frac{4}{n-4}}r_k^{2-\alpha}\) bounds the difference
between \(\mathcal R_k(\lambda,r)\) and the pair of right-hand sides in \eqref{eq:ce-point-parameter-sum} by
\begin{equation}\label{eq:ce-reduced-limit-error}
 o_b(1)(1+k)\gamma^{\beta k}.
\end{equation}

At \(\lambda=\lambda^0\), \(r_j=\gamma^j\), both right-hand sides of
\eqref{eq:ce-point-parameter-sum} vanish by \(H(1)=H'(1)=0\).
For \((\mathbf p,\mathbf r)\in Y\), set
\begin{equation}\label{eq:ce-parameter-path}
 \lambda_j(t)=\lambda_j^0+t\mathbf p_j,
 \,\ r_j(t)=\gamma^j+t\mathbf r_j,\,\  j\geq-3.
\end{equation}
We define \(L:Y\to\mathcal X_\rho\) by
\begin{equation}\label{eq:ce-limit-linearization}
 [L(\mathbf p,\mathbf r)]_k
 :=\left.\frac{d}{dt}\right|_{t=0}
 \left[\lim_{b\to0}\mathcal R_k(\lambda(t),r(t))\right],
 \,\  k\geq0.
\end{equation}
Here \(\frac{\mathbf p_j}{\lambda_j^0}\) and \(\frac{\mathbf r_j}{\gamma^j}\)
are kept fixed as \(b\to0\), and the limit is taken before
differentiation.

For \(\lambda,r\) satisfying \eqref{eq:ce-parameter-space}, with
\((\mathbf p,\mathbf r)\in Y\), define \(E_b\) by
\begin{equation}\label{eq:ce-reduced-parameter-equation}
 \mathcal R(\lambda,r)=L(\mathbf p,\mathbf r)+E_b(\lambda,r).
\end{equation}
Thus \(E_b\) contains both the error in taking \(b\to0\) and the
Taylor remainder of \(\lim_{b\to0}\mathcal R(\lambda,r)\).

In Subsection~\ref{app:parameter-surjectivity}, we choose \(\rho\)
in \((\gamma^\beta,1)\) and construct a right inverse of \(L\) whose bound,
with the norm \eqref{eq:ce-parameter-increment-norm}, is independent
of \(b\).

Subsection~\ref{subsec:ce-regular-parts} gives the estimates for
\(E_b(\lambda,r)\).

\subsection{Surjectivity of \texorpdfstring{\(L\)}{L}}\label{app:parameter-surjectivity}
To compute the coefficients of \(L\) in
\eqref{eq:ce-limit-linearization}, set
\[
\begin{aligned}
 \mathcal F(s)&:=\Phi_1(s\xi_1)\\
 &=\frac1\ell\sum_{i=1}^{\ell}\mathcal K_{\xi_i}(s\xi_1)
   -1-\frac{(n-4)^2}{4}s^2,\,\  s\ne\pm1.
\end{aligned}
\]
For \(j\geq-3\) with \(j\ne k\),
\eqref{eq:ce-reference-parameters} and scaling give
\begin{equation}\label{eq:ce-circle-parameter-vector}
\begin{aligned}
 &\gamma^{j(n-2+\alpha)-k(2+\alpha)}
 \left(\frac{\lambda_j}{\lambda_j^0}\right)^{\frac{n-4}{2}}
 \begin{pmatrix}
  \Phi_{r_j}(\mathbf x_{k,1})\\
  r_k\partial_{x_1}\Phi_{r_j}(\mathbf x_{k,1})
 \end{pmatrix}\\
 &\,\ =
 \gamma^{(j-k)(2+\alpha)}
 \left(\frac{\lambda_j}{\lambda_j^0}\right)^{\frac{n-4}{2}}
 \left(\frac{r_j}{\gamma^j}\right)^{4-n}
 \begin{pmatrix}\mathcal F(s)\\s\mathcal F'(s)\end{pmatrix},
 \,\  s=\frac{r_k}{r_j}.
\end{aligned}
\end{equation}
In \eqref{eq:ce-circle-parameter-vector}, the derivatives
\(\frac{2\lambda_j^0}{n-4}\partial_{\lambda_j}\) and
\(\gamma^j\partial_{r_j}\), evaluated at
\(\lambda=\lambda^0\), \(r_j=\gamma^j\) and \(r_k=\gamma^k\),
holding \(r_k\) fixed during differentiation, give the columns of
\[
 C_{j-k}=\gamma^{(2+\alpha)(j-k)}
 \left.
 \begin{pmatrix}
 \mathcal F(s)&(4-n)\mathcal F(s)-s\mathcal F'(s)\\
 s\mathcal F'(s)&(3-n)s\mathcal F'(s)-s^2\mathcal F''(s)
 \end{pmatrix}\right|_{s=\gamma^{k-j}}.
\]
The columns of \(C_0\) are obtained by applying
\(\frac{2\lambda_k^0}{n-4}\partial_{\lambda_k}\) and
\(\gamma^k\partial_{r_k}\) to \eqref{eq:ce-point-parameter-sum}
at \(\lambda=\lambda^0\), \(r_j=\gamma^j\).
These derivatives include the summand containing
\(\Phi_{r_k}-\frac{\mathcal K_{r_k\xi_1}}{\ell}\) in
\eqref{eq:ce-self-circle-regular-part}, the evaluation point
\(x=r_k\xi_1\) in every sum in \eqref{eq:ce-point-parameter-sum},
and the factor \(r_k\) in \(\mathcal R_{k,2}\).
By scaling and \(H(1)=H'(1)=0\),
\[
 C_0=-\sum_{m\ne0}C_m.
\]
Differentiating \eqref{eq:ce-point-parameter-sum} along
\eqref{eq:ce-parameter-path}, we therefore obtain
\begin{equation}\label{eq:ce-sequence-operator}
 [L(\mathbf p,\mathbf r)]_k
 =\sum_{j\geq-3}C_{j-k}
 \begin{pmatrix}
 \displaystyle\frac{n-4}{2}\frac{\mathbf p_j}{\lambda_j^0}\\[1mm]
 \displaystyle\frac{\mathbf r_j}{\gamma^j}
 \end{pmatrix},\,\  k\geq0.
\end{equation}
The matrices \(C_m\) are independent of \(b\); the dependence of
\(L\) on \(b\) is through \(\lambda_j^0\).
The same factors \(\lambda_j^0\) occur in the norm
\eqref{eq:ce-parameter-increment-norm}.

To solve \(\mathcal R(\lambda,r)=0\), we construct a bounded right inverse
\(\mathcal S\) of \(L\). If \((\mathbf p,\mathbf r)=-\mathcal S E_b(\lambda,r)\),
then by \eqref{eq:ce-reduced-parameter-equation} and \(L\mathcal S=I\),
\[
 \mathcal R(\lambda,r)=-L\mathcal S E_b(\lambda,r)+E_b(\lambda,r)=0.
\]
\begin{lemma}\label{lem:ce-sequence}
There exists \(\rho\in(\gamma^\beta,1)\), independent of \(b\), such that
\(L:Y\to\mathcal X_\rho\) is onto and has a bounded real linear
right inverse \(\mathcal S:\mathcal X_\rho\to Y\) satisfying
\[
 L\mathcal S=I,
 \,\ \|\mathcal S f\|_Y\leq C\|f\|_{\mathcal X_\rho},
\]
where \(C\) is independent of \(b\).
\end{lemma}
\begin{proof}
First set \(\mathbf p_{-j}=\mathbf r_{-j}=0\), \(j=1,2,3\), and define
\[
 Y_+=\{(\mathbf p,\mathbf r)\in Y:
       \mathbf p_{-j}=\mathbf r_{-j}=0,\ j=1,2,3\}.
\]
The map
\begin{equation}\label{eq:ce-parameter-isometry}
 (\mathbf p,\mathbf r)\mapsto
 \left(
 \begin{pmatrix}
 \displaystyle\frac{n-4}{2}\frac{\mathbf p_k}{\lambda_k^0}\\[1mm]
 \displaystyle\frac{\mathbf r_k}{\gamma^k}
 \end{pmatrix}
 \right)_{k\geq-3}
\end{equation}
identifies \(Y\) isometrically with the space of two-component
sequences that tend to zero after multiplication by \(\rho^{-k}\).
Under \eqref{eq:ce-parameter-isometry}, \(L\) is the convolution
with the matrices \(C_j\).
The coefficients \(C_j\) in \eqref{eq:ce-sequence-operator}
give the Laurent series
\[
 \mathcal Q(z)=\sum_{j\in\mathbb Z}C_jz^j
             =\sum_{j\ne0}C_j(z^j-1)
             =\begin{pmatrix}
                \mathcal Q_{11}(z)&\mathcal Q_{12}(z)\\
                \mathcal Q_{21}(z)&\mathcal Q_{22}(z)
               \end{pmatrix},
\]
where
\[
\begin{aligned}
 \mathcal Q_{11}(z)
 &=\sum_{j\ne0}\gamma^{(2+\alpha)j}
       \mathcal F(\gamma^{-j})(z^j-1),\\
 \mathcal Q_{21}(z)
 &=\sum_{j\ne0}\gamma^{(1+\alpha)j}
       \mathcal F'(\gamma^{-j})(z^j-1),\\
 \mathcal Q_{12}(z)&=(4-n)\mathcal Q_{11}(z)-\mathcal Q_{21}(z),\\
 \mathcal Q_{22}(z)
 &=(3-n)\mathcal Q_{21}(z)
   -\sum_{j\ne0}\gamma^{\alpha j}\mathcal F''(\gamma^{-j})(z^j-1).
\end{aligned}
\]
The entries \(\mathcal Q_{ab}(z)\), \(a,b=1,2\), are analytic for
\(\gamma^{2-\alpha}<|z|<\gamma^{-\alpha}\). To use the Laurent
coefficients of \(\mathcal Q^{-1}\), we will choose
\[
 \gamma^\beta<\rho<1,
 \,\ \det\mathcal Q(z)\ne0\quad\text{for }|z|=\rho.
\]
The numbers \(\gamma,\alpha\) satisfy \eqref{eq:ce-balance} and
are fixed before \(b\).

The circle average in \(\mathcal F\) has scalar Taylor coefficients
\(c_m\). For the fixed even integer \(\ell\geq6\), the binomial
coefficients
\[
 u_j=\frac{(\frac{n-4}{2})_j}{j!},\,\  j\geq0,
\]
give
\begin{equation}\label{eq:ce-full-circle-coefficients}
 \frac1\ell\sum_{i=1}^{\ell}|s\xi_1-\xi_i|^{4-n}
 =\sum_{\substack{j,k\geq0\\\ell\mid(j-k)}}u_ju_ks^{j+k}
 =\sum_{\substack{m\geq0\\m\ {\rm even}}}c_ms^m,\,\  |s|<1.
\end{equation}
Here \((a)_j=a(a+1)\cdots(a+j-1)\), with \((a)_0=1\).
Averaging over the angular points retains exactly the pairs
whose indices differ by a multiple of \(\ell\). Thus
\[
 c_0=1,\,\  c_2=\frac{(n-4)^2}{4},
 \,\ \mathcal F(s)=\sum_{\substack{m\geq4\\m\ {\rm even}}}c_ms^m
 \quad (|s|<1).
\]
The pair \((j,j)\) gives \(u_j^2\), while summing over all pairs
with total degree \(2j\) gives
\[
 u_j^2\leq c_{2j}\leq\frac{(n-4)_{2j}}{(2j)!}.
\]
Consequently \(c_{2j}\geq C^{-1}(1+j)^{n-6}\) and
\(c_{2j}\leq C(1+j)^{n-5}\), where \(C\) may depend on \(n,\ell\).

Using \(\mathcal Q_{12}(z)=(4-n)\mathcal Q_{11}(z)-\mathcal Q_{21}(z)\),
we define the symmetric matrix
\begin{equation}\label{eq:app-symbol-factorization}
\begin{aligned}
 \widetilde{\mathcal Q}(z)
 &:=\frac1{z-1}
 \begin{pmatrix}1&0\\0&-1\end{pmatrix}
 \mathcal Q(z)
 \begin{pmatrix}1&n-4\\0&1\end{pmatrix}\\
 &=\frac1{z-1}
 \begin{pmatrix}
  \mathcal Q_{11}(z)&-\mathcal Q_{21}(z)\\
  -\mathcal Q_{21}(z)&(4-n)\mathcal Q_{21}(z)-\mathcal Q_{22}(z)
 \end{pmatrix},\\
 \det\mathcal Q(z)&=-(z-1)^2\det\widetilde{\mathcal Q}(z).
\end{aligned}
\end{equation}
Since \(\mathcal Q(1)=0\), we extend \(\widetilde{\mathcal Q}\) analytically
to \(z=1\).
For any
\(\gamma^{2-\alpha}<\rho<1\), expansion of the circle coefficients gives
\begin{equation}\label{eq:ce-pole-expansion}
\widetilde{\mathcal Q}(\rho z)
 =-\sum_{l\notin\{0,2\}}\frac{\mu_l}{z-D_l}
       \begin{pmatrix}1&-l\\-l&l^2\end{pmatrix}
  +\frac{\mu_0}{z-D_0}\begin{pmatrix}1&0\\0&0\end{pmatrix}
       +\frac{\mu_2}{z-D_2}\begin{pmatrix}1&-2\\-2&4\end{pmatrix}.
\end{equation}
Here \(D_l\) are the pole locations and \(\mu_l>0\) are their
coefficients. The index \(l\) labels terms in this expansion,
rather than circles. It ranges over even \(l\geq4\),
\(l=4-n-m\) with even \(m\geq0\), and \(l=0,2\).
The summands with \(l=0,2\) in \eqref{eq:ce-pole-expansion}
come from \(1\) and \(\frac{(n-4)^2s^2}{4}\) subtracted in
\(\mathcal F(s)\).

For even \(m\geq4\),
\[
 D_m=\frac{\gamma^{m-2-\alpha}}{\rho}<1,
 \,\ \mu_m=\frac{c_mD_m}{1-\rho D_m}.
\]
For even \(m\geq0\),
\[
 D_{4-n-m}=\frac{\gamma^{2-n-m-\alpha}}{\rho}>1,
 \,\ \mu_{4-n-m}=\frac{c_m}{\rho(1-\gamma^{n-2+m+\alpha})}.
\]
For \(l=0,2\), the pairs \((D_l,\mu_l)\) are
\[
(D_0,\mu_0)=
 \left(\frac1{\rho\gamma^{2+\alpha}},
       \frac1{\rho(1-\gamma^{2+\alpha})}\right),
 \,\ (D_2,\mu_2)=
 \left(\frac1{\rho\gamma^\alpha},
       \frac{(n-4)^2}{4\rho(1-\gamma^\alpha)}\right).
\]
The formulas use only \(0<\gamma<1\), \(0<\alpha<2\) and
\(\gamma^{2-\alpha}<\rho<1\). They give
\begin{equation}\label{eq:ce-pole-summability}
\begin{gathered}
 \sup_{D_l<1}D_l<1,\,\  \sup_{D_l>1}D_l^{-1}<1,\\
 \sum_{D_l<1}\mu_l(1+l^2)<\infty,
 \,\ \sum_{D_l>1}\mu_lD_l^{-1}(1+l^2)<\infty.
\end{gathered}
\end{equation}

\subsubsection{The determinant \texorpdfstring{\(\det\mathcal Q\)}{det Q} and the range \texorpdfstring{\(L(Y_+)\)}{L(Y+)}}
First take any \(\gamma^{2-\alpha}<\rho<1\). The coefficient bound
\(c_m\geq C^{-1}(1+m)^{-1}\), for even \(m\geq4\), implies
\begin{equation}\label{eq:ce-positive-partial-sums}
 \sum_{\substack{4\leq m\leq N\\m\ {\rm even}}}
       c_m(w_1-mw_2)^2\to\infty
 \quad\text{for every }w=(w_1,w_2)\in\mathbb R^2\setminus\{0\}.
\end{equation}
The partial sums in \eqref{eq:ce-positive-partial-sums}
increase with \(N\). Compactness of the unit circle
therefore gives \(N\) large enough that
\[
 \frac12\sum_{\substack{4\leq m\leq N\\m\ {\rm even}}}
       c_m(w_1-mw_2)^2
 -\sum_{l=0,2}\frac{\mu_l}{D_l}(w_1-lw_2)^2
 \geq1,\,\  |w|=1.
\]
Fix such an \(N\), and choose
\(0<\delta\leq\min\{D_m:4\leq m\leq N,\ m\text{ even}\}\).
For even \(m\) with \(4\leq m\leq N\),
\[
 \frac{\mu_m}{\delta+D_m}
 \geq\frac{c_m}{2(1-\rho D_m)}\geq\frac{c_m}{2}.
\]
For \(l=0,2\), \(\frac{\mu_l}{\delta+D_l}\leq\frac{\mu_l}{D_l}\).
The expansion of \(\widetilde{\mathcal Q}\) now gives
\[
 \widetilde{\mathcal Q}(-\rho\delta)\geq I_2.
\]
The two convergent sums in \eqref{eq:ce-pole-summability} imply uniform
convergence of the series for \(\widetilde{\mathcal Q}(\rho z)\) on
\(|z-z_0|\leq\frac{|z_0|}{2}\), for every \(z_0<0\).
Thus it extends analytically along the segment from \(-1\) to
\(-\delta\). If \(\det\mathcal Q\) vanished identically on the
annulus, \eqref{eq:app-symbol-factorization} would give
\(\det\widetilde{\mathcal Q}(-\rho\delta)=0\), a contradiction.

The zeros of \(\det\mathcal Q\) in each closed subannulus are
isolated and finite. Since \(\beta\leq2-\alpha\), we may choose
\[
 \gamma^\beta<\rho<1,
 \,\ \det\mathcal Q(z)\ne0\quad (|z|=\rho).
\]
We fix this \(\rho\). Both \(\mathcal Q(z)\) and
\(\mathcal Q(z)^{-1}\) are analytic in a neighborhood of
\(\{z:|z|=\rho\}\).
Writing \([z^m]\) for the coefficient of \(z^m\) in a Laurent
series, we have
\[
 \sum_{m\in\mathbb Z}\rho^m
 \left(\|C_m\|+\bigl\|[z^m]\mathcal Q(z)^{-1}\bigr\|\right)<\infty.
\]
The matrix norms are induced by \((a_1,a_2)\mapsto|a_1|+|a_2|\).
The Laurent coefficients \([z^m]\mathcal Q(z)^{-1}\) define a bounded map
\(B_+:\mathcal X_\rho\to Y_+\):
\[
\begin{aligned}
 (B_+g)_k&=
 \begin{pmatrix}
  \dfrac{2\lambda_k^0}{n-4}&0\\
  0&\gamma^k
 \end{pmatrix}
 \sum_{j\geq0}\bigl([z^{j-k}]\mathcal Q(z)^{-1}\bigr)g_j,
 \,\  k\geq0,\\
 (B_+g)_{-1}&=(B_+g)_{-2}=(B_+g)_{-3}=0.
\end{aligned}
\]
Indeed,
\[
 \|B_+g\|_Y
 \leq\left(\sum_{m\in\mathbb Z}\rho^m
             \bigl\|[z^m]\mathcal Q(z)^{-1}\bigr\|\right)
       \|g\|_{\mathcal X_\rho}.
\]
The factors \(\frac{2\lambda_k^0}{n-4}\) and \(\gamma^k\) in
\(B_+\) cancel with \(\frac{n-4}{2\lambda_k^0}\) and \(\gamma^{-k}\)
in \eqref{eq:ce-parameter-increment-norm}.
Truncating the sum over \(j\) gives \(B_+g\in Y_+\).

In the coordinates \eqref{eq:ce-parameter-isometry}, truncate the
Laurent series for \(\mathcal Q\) and \(\mathcal Q^{-1}\) to finitely many powers.
For input and output indices \(j,k\geq0\), the difference between
the product with intermediate index \(m\geq0\) and the product
with \(m\in\mathbb Z\) consists of finitely many terms with
\(m<0\) and has finite-dimensional range. The convergence of
\(\sum_m\rho^m\|C_m\|\) and
\(\sum_m\rho^m\|[z^m]\mathcal Q(z)^{-1}\|\)
allows passage to the operator norm limit. Since \(\mathcal Q\mathcal Q^{-1}=\mathcal Q^{-1}\mathcal Q=I_2\),
\[
 (L|_{Y_+})B_+-I_{\mathcal X_\rho},
 \,\  B_+(L|_{Y_+})-I_{Y_+}
 \quad\text{are compact}.
\]
{\emergencystretch=2em
Consequently \(L|_{Y_+}\) has finite-dimensional kernel,
closed range and finite-dimensional quotient
\(\mathcal X_\rho/L(Y_+)\).
\par}

To prove \(L(Y)=\mathcal X_\rho\), consider
a continuous linear functional on \(\mathcal X_\rho\). It has the form
\begin{equation}\label{eq:ce-dual-functional}
 g\mapsto\sum_{k\geq0}\rho^{-k}h_k^Tg_k,
 \,\ h_k\in\mathbb R^2,\,\  \sum_{k\geq0}\bigl(|h_{k,1}|+|h_{k,2}|\bigr)<\infty.
\end{equation}
The functional in \eqref{eq:ce-dual-functional} vanishes on
\(L(Y_+)\) precisely when
\[
 \sum_{k\geq0}\rho^{-k}C_{j-k}^Th_k=0,\,\  j\geq0.
\]
Define the vector power series
\[
 f(z)=(\rho z-1)\sum_{k\geq0}
       \begin{pmatrix}h_{k,1}\\-h_{k,2}\end{pmatrix}z^k.
\]
The definition of \(\widetilde{\mathcal Q}\) gives
\[
 \mathcal Q(\rho z)^T\sum_{k\geq0}h_kz^k
 =\begin{pmatrix}1&0\\4-n&1\end{pmatrix}
       \widetilde{\mathcal Q}(\rho z)f(z).
\]
Thus vanishing of the functional on \(L(Y_+)\) is equivalent to
\[
 [z^k]\bigl(\widetilde{\mathcal Q}(\rho z)f(z)\bigr)=0,
 \,\  k\geq0.
\]
Since \(\rho<1\), the correspondence between \((h_k)\) and \(f\)
is invertible among power series with absolutely summable
coefficient norms. The Laurent coefficients
\([z^m]\widetilde{\mathcal Q}(\rho z)\), \(m\in\mathbb Z\), are
absolutely summable. Consequently, each coefficient of
\(\widetilde{\mathcal Q}(\rho z)f(z)\) is an absolutely convergent sum.

\subsubsection{The coefficient equations for \texorpdfstring{\(f\)}{f}}
Write \(f(z)=\sum_{k\geq0}f_kz^k\). We show that the vector power
series satisfying
\begin{equation}\label{eq:ce-zero-constant-coefficient-equations}
 \sum_{k\geq0}\bigl(|f_{k,1}|+|f_{k,2}|\bigr)<\infty,
 \,\ [z^k]\bigl(\widetilde{\mathcal Q}(\rho z)f(z)\bigr)=0\ (k\geq0),
 \,\  f(0)=0
\end{equation}
form a space of dimension at most two. For each pole index \(l\)
and coefficient index \(k\), define
\[
 y_{l,k}(f)=
 \begin{cases}
  \displaystyle\sum_{j\geq0}D_l^j\sqrt{\mu_l}(1,-l)f_{k+j+1},
       &D_l<1,\\
  \displaystyle-\sum_{j=0}^kD_l^{-(k-j+1)}
                   \sqrt{\mu_l}(1,-l)f_j,&D_l>1.
 \end{cases}
\]
Comparing the nonnegative powers gives
\begin{equation}\label{eq:ce-coefficient-recurrence}
\begin{aligned}
 \sum_{l\notin\{0,2\}}\sqrt{\mu_l}\begin{pmatrix}1\\-l\end{pmatrix}y_{l,k}
 &=\sqrt{\mu_0}\begin{pmatrix}1\\0\end{pmatrix}y_{0,k}
   +\sqrt{\mu_2}\begin{pmatrix}1\\-2\end{pmatrix}y_{2,k},\\
 y_{l,k+1}
 &=D_l^{-1}\bigl(y_{l,k}-\sqrt{\mu_l}(1,-l)f_{k+1}\bigr).
\end{aligned}
\end{equation}
Write \(y_k(f)=(y_{l,k}(f))_l\). The terms with \(l=0,2\)
have negative signs in the form
\[
 [x,y]=\sum_{l\notin\{0,2\}}D_l\overline{x_l}y_l
       -D_0\overline{x_0}y_0-D_2\overline{x_2}y_2.
\]
We use \([\cdot,\cdot]\) on sequences with
\(\sum_lD_l|x_l|^2<\infty\), with norm
\((\sum_lD_l|x_l|^2)^{\frac{1}{2}}\).
Only the coordinates \(l=0,2\) have negative coefficients.
The definitions give
\begin{equation}\label{eq:ce-coefficient-sequence-bounds}
 |y_{l,k}|\leq C\sqrt{\mu_l}\sqrt{1+l^2}
 \begin{cases}
  \displaystyle\sup_{j\geq k+1}\bigl(|f_{j,1}|+|f_{j,2}|\bigr),&D_l<1,\\
  \displaystyle D_l^{-1}\sum_{j=0}^kD_l^{-j}\bigl(|f_{k-j,1}|+|f_{k-j,2}|\bigr),&D_l>1.
 \end{cases}
\end{equation}
Since \(f_j\to0\) and \(\sup_{D_l>1}D_l^{-1}<1\),
\[
 \sup_{j\geq k+1}\bigl(|f_{j,1}|+|f_{j,2}|\bigr)\to0,
 \,\ \sup_{D_l>1}\sum_{j=0}^kD_l^{-j}\bigl(|f_{k-j,1}|+|f_{k-j,2}|\bigr)\to0
 \quad\text{as }k\to\infty.
\]
Combining \eqref{eq:ce-pole-summability} with \eqref{eq:ce-coefficient-sequence-bounds} gives
\begin{equation}\label{eq:ce-coefficient-decay}
 \sup_{k\geq0}\left(\sum_lD_l|y_{l,k}|^2\right)^{\frac{1}{2}}
       \leq C\sup_{j\geq0}\bigl(|f_{j,1}|+|f_{j,2}|\bigr),
 \,\ \sum_lD_l|y_{l,k}|^2\to0.
\end{equation}
Equations~\eqref{eq:ce-pole-summability} and
\eqref{eq:ce-coefficient-sequence-bounds} also give
\begin{equation}\label{eq:ce-coefficient-pairing-bounds}
 \sum_l|y_{l,i}(f)y_{l,j}(g)|<\infty,
 \,\ \sum_l\sqrt{\mu_l}\sqrt{1+l^2}|y_{l,k}|<\infty
\end{equation}
for \(f,g\) satisfying \eqref{eq:ce-zero-constant-coefficient-equations}.
The bounds \eqref{eq:ce-coefficient-sequence-bounds}--\eqref{eq:ce-coefficient-pairing-bounds}
justify the coefficient comparisons, recurrences and pairings.

Expanding the recurrence in \eqref{eq:ce-coefficient-recurrence} and using its first identity gives
\[
 [y_{i+1}(f),y_j(g)]=[y_i(f),y_{j+1}(g)],
 \,\ [y_i(f),y_j(g)]=[y_0(f),y_{i+j}(g)].
\]
Combine the sequences \(y_j(f)\) using the coefficients of
\(\sqrt{1-z}\):
\[
 \widehat y(f)=\sum_{j\geq0}\frac{(-\frac{1}{2})_j}{j!}y_j(f).
\]
The series for \(\widehat y(f)\) converges absolutely in
\((\sum_lD_l|x_l|^2)^{\frac{1}{2}}\), since
\(\sum_{j\geq0}|\frac{(-\frac{1}{2})_j}{j!}|=2\). The identity
\((\sqrt{1-z})^2=1-z\) then gives
\[
 [\widehat y(f),\widehat y(g)]
 =[y_0(f),y_0(g)]-[y_0(f),y_1(g)].
\]
If \(f(0)=0\), then \(y_{l,0}(f)=0\) for \(D_l>1\), including \(l=0,2\).
The first identity in \eqref{eq:ce-coefficient-recurrence} at \(k=0\)
therefore gives
\[
 \sum_{D_l<1}\sqrt{\mu_l}\begin{pmatrix}1\\-l\end{pmatrix}y_{l,0}(f)=0.
\]
Combining the two identities in \eqref{eq:ce-coefficient-recurrence} at \(k=0\),
we obtain
\[
 [y_0(f),y_1(f)]=\sum_{D_l<1}|y_{l,0}(f)|^2.
\]
Hence
\[
 [\widehat y(f),\widehat y(f)]
 =-\sum_{D_l<1}(1-D_l)|y_{l,0}(f)|^2.
\]
For nonzero \(f\), \([\widehat y(f),\widehat y(f)]<0\). Indeed, equality to zero
would imply \((1,-m)f(D_m)=0\) for every even \(m\geq4\).
If \(f_p\ne0\) is the first nonzero Taylor coefficient, division by
\(D_m^p\) gives
\[
 (f_p)_1-m(f_p)_2+O(mD_m)=0.
\]
Since \(D_m=\frac{\gamma^{m-2-\alpha}}{\rho}\) tends exponentially to zero,
dividing \((f_p)_1-m(f_p)_2+O(mD_m)=0\) by \(m\) and taking
\(m\to\infty\) gives \((f_p)_2=0\). Substituting \((f_p)_2=0\)
in the same identity and taking \(m\to\infty\) gives \((f_p)_1=0\),
a contradiction.
It follows that
\[
 f\mapsto\bigl(\widehat y_0(f),\widehat y_2(f)\bigr)
 \quad\text{is injective}.
\]
The subscripts \(0,2\) refer to the two negative coordinates of the
form \([\cdot,\cdot]\), so the space of solutions of
\eqref{eq:ce-zero-constant-coefficient-equations} has dimension at most two.

\subsubsection{Variation of \texorpdfstring{\(\lambda_{-j},r_{-j}\)}{lambda(-j), r(-j)}, \(j=1,2,3\)}
For each \(j=1,2,3\), varying
\(\mathbf p_{-j},\mathbf r_{-j}\) adds any sequence
\((C_{-k-j}w)_{k\geq0}\), \(w\in\mathbb R^2\), to
\(L(\mathbf p,\mathbf r)\).
Since \(\rho>\gamma^{2-\alpha}\),
\((C_{-k-j}w)_{k\geq0}\in\mathcal X_\rho\).
Applying the functional in \eqref{eq:ce-dual-functional} gives
\[
\sum_{k\geq0}\rho^{-k}h_k^TC_{-k-j}w
 =\rho^jw^T\begin{pmatrix}1&0\\4-n&1\end{pmatrix}
  \times[z^{-j}]
        \bigl(\widetilde{\mathcal Q}(\rho z)f(z)\bigr).
\]
If the functional in \eqref{eq:ce-dual-functional} is nonzero and
vanishes on \(L(Y_+)\) and on
\((C_{-k-j}w)_{k\geq0}\) for all \(w\in\mathbb R^2\),
\(j=1,2,3\), then
\[
 [z^{-j}]\bigl(\widetilde{\mathcal Q}(\rho z)f(z)\bigr)=0,
 \,\  j=1,2,3.
\]
Thus \(zf,z^2f,z^3f\) satisfy
\eqref{eq:ce-zero-constant-coefficient-equations}.
These three functions are linearly independent, contradicting
the dimension bound for the solution space of
\eqref{eq:ce-zero-constant-coefficient-equations}.
Hence any continuous functional vanishing on \(L(Y_+)\) and on
\((C_{-k-j}w)_{k\geq0}\) for every \(w\in\mathbb R^2\),
\(j=1,2,3\), is zero.

The range \(L(Y_+)\) is closed and has finite codimension.
Adding the contributions of
\(\mathbf p_{-j},\mathbf r_{-j}\), \(j=1,2,3\), preserves closedness,
so the range of \(L\) is closed and dense, giving
\[
\begin{aligned}
 L(Y_+)+
 \left\{\left(\sum_{j=1}^3C_{-k-j}w_j\right)_{k\geq0}:
          w_1,w_2,w_3\in\mathbb R^2\right\}
 &=\mathcal X_\rho.
\end{aligned}
\]
Hence \(L\) is onto. Its kernel has finite dimension;
restricting \(L\) to a closed complement of its kernel gives the
bounded right inverse \(\mathcal S\).
Choose the complement of \(\ker L\) after applying
\eqref{eq:ce-parameter-isometry}. The matrices \(C_j\) and the weighted supremum norm of the
two-component sequences in \eqref{eq:ce-parameter-isometry}
are independent of \(b\). Transferring the right inverse back to \(Y\) therefore gives
\(\|\mathcal S f\|_Y\leq C\|f\|_{\mathcal X_\rho}\), with \(C\)
independent of \(b\).
\end{proof}

\subsection{Estimates for \texorpdfstring{\(E_b\)}{E(b)}}\label{subsec:ce-regular-parts}
For the fixed point argument in Subsection~\ref{subsec:ce-parameter-choice}, we
estimate \(E_b\) in \eqref{eq:ce-reduced-parameter-equation} and show that it has
relatively compact image in \(\mathcal X_\rho\).

\begin{lemma}\label{lem:ce-reduced}
For every fixed small \(b\), consider \(\lambda,r\) satisfying
\eqref{eq:ce-parameter-space} with \(\varepsilon=\varepsilon_0\)
and \((\mathbf p,\mathbf r)\in Y\).
Then \(E_b(\lambda,r)\) is continuous in the norm
\eqref{eq:ce-parameter-increment-norm} and has relatively compact
image in \(\mathcal X_\rho\). Moreover,
\begin{equation}\label{eq:ce-reduced-remainder-bound}
 \rho^{-k}\bigl(|E_b(\lambda,r)_{k,1}|+|E_b(\lambda,r)_{k,2}|\bigr)
 \leq C\bigl(\|(\mathbf p,\mathbf r)\|_Y^2+o_b(1)\bigr)(k+1)
 \left[\rho^k+\left(\frac{\gamma^\beta}{\rho}\right)^k\right],\,\  k\geq0.
\end{equation}
\end{lemma}
\begin{proof}
To prove \eqref{eq:ce-reduced-remainder-bound}, we combine
\eqref{eq:ce-reduced-limit-error} with Taylor's theorem applied
to both right-hand sides of
\eqref{eq:ce-point-parameter-sum} along
\eqref{eq:ce-parameter-path}, \(0\leq t\leq1\).
The condition \eqref{eq:ce-parameter-space} holds along
\eqref{eq:ce-parameter-path}:
\(\lambda_j(t)\) lies between \(\lambda_j^0\) and \(\lambda_j\),
and \(r_j(t)\) lies between \(\gamma^j\) and \(r_j\).
In \eqref{eq:ce-circle-parameter-vector}, the derivatives
\(\lambda_j^0\partial_{\lambda_j}\),
\(\gamma^j\partial_{r_j}\) and
\(\gamma^k\partial_{r_k}\), of total order at most two, are bounded by
\(C\gamma^{(2-\alpha)(k-j)}\) for \(j<k\) and
\(C\gamma^{\alpha(j-k)}\) for \(j>k\).
Consequently, the second derivative with respect to \(t\)
of the vector in \eqref{eq:ce-circle-parameter-vector}, evaluated
at \(\lambda(t),r(t)\), is bounded by
\[
 C\|(\mathbf p,\mathbf r)\|_Y^2(\rho^j+\rho^k)^2
 \begin{cases}
 \gamma^{(2-\alpha)(k-j)},&j<k,\\
 \gamma^{\alpha(j-k)},&j>k.
 \end{cases}
\]
The second \(t\)-derivatives of the \(j=k\) contribution to
\eqref{eq:ce-point-parameter-sum}, using
\eqref{eq:ce-self-circle-regular-part}, and of the sum indexed by
\(j\leq-4\) in \eqref{eq:ce-point-parameter-sum} are bounded by
\(C\|(\mathbf p,\mathbf r)\|_Y^2\rho^{2k}\).
Using \((\rho^j+\rho^k)^2\leq2(\rho^{2j}+\rho^{2k})\),
Taylor's theorem gives
\begin{equation}\label{eq:ce-parameter-quadratic-error}
\begin{aligned}
 &\sum_{\nu=1}^2\left|\lim_{b\to0}\mathcal R_{k,\nu}(\lambda,r)
       -[L(\mathbf p,\mathbf r)]_{k,\nu}\right|\\
 &\quad\leq C\|(\mathbf p,\mathbf r)\|_Y^2\left[
 \rho^{2k}
 +\sum_{j=0}^k\gamma^{(2-\alpha)(k-j)}\rho^{2j}
 +\sum_{j>k}\gamma^{\alpha(j-k)}\rho^{2j}
 +\gamma^{(2-\alpha)k}\right].
\end{aligned}
\end{equation}
The term \(C\|(\mathbf p,\mathbf r)\|_Y^2\gamma^{(2-\alpha)k}\) in
\eqref{eq:ce-parameter-quadratic-error} bounds the terms quadratic
in \(\mathbf p_{-j},\mathbf r_{-j}\), \(j=1,2,3\):
after \(I-J_2\), their values are \(O(r_k^4)\), and their radial
derivatives are \(O(r_k^3)\).
The sum \(\sum_{j>k}\gamma^{\alpha(j-k)}\rho^{2j}\)
in \eqref{eq:ce-parameter-quadratic-error} is at most
\(C\rho^{2k}\). For the sum over \(0\leq j\leq k\),
\begin{equation}\label{eq:ce-parameter-geometric-sum}
 \rho^{-k}\sum_{j=0}^k\gamma^{(2-\alpha)(k-j)}\rho^{2j}
 =\sum_{j=0}^k
     \left(\frac{\gamma^{2-\alpha}}\rho\right)^{k-j}\rho^j
 \leq(k+1)\left[\rho^k+\left(\frac{\gamma^\beta}{\rho}\right)^k\right].
\end{equation}
Equation~\eqref{eq:ce-parameter-geometric-sum} also covers
\(\rho^2=\gamma^{2-\alpha}\), when the sum before division by \(\rho^k\)
is exactly \((k+1)\rho^{2k}\).
After division by \(\rho^k\), the terms \(\rho^{2k}\) and
\(\gamma^{(2-\alpha)k}\) in \eqref{eq:ce-parameter-quadratic-error}
are also bounded by \(C(k+1)[\rho^k+(\frac{\gamma^\beta}{\rho})^k]\).
Adding \(o_b(1)(1+k)(\frac{\gamma^\beta}{\rho})^k\), from
\eqref{eq:ce-reduced-limit-error}, proves the bound in
Lemma~\ref{lem:ce-reduced}.

To prove continuity, let \((\mathbf p,\mathbf r)\) converge in \(Y\).
For each fixed \(k,i\), the functions
\(\mathcal B_{\lambda_k,\mathbf x_{k,i}}\), \(\chi_{k,i}\)
and \(\zeta_{k,i}\) converge in \(C^q(\overline\Omega)\) for every
integer \(q\geq0\), as do the terms containing \(\lambda_{-j},r_{-j}\),
\(j=1,2,3\), in \(h_{r,\lambda}\).
On each compact subset of \(\Omega\setminus\{0\}\), for large \(K\),
\eqref{eq:ce-parameter-space} and the definition of \(J_2\) give
\[
 \left|\sum_{k>K}(I-J_2)
 \left[\sum_{i=1}^{\ell}\mathcal B_{\lambda_k,\mathbf x_{k,i}}\right](x)\right|
 \leq Cb\sum_{k>K}r_k^\alpha
 \leq Cb\sum_{k>K}\gamma^{k\alpha}\to0
 \quad\text{as }K\to\infty.
\]
Thus \(U_{\lambda,r}\) converges locally uniformly.
By interior elliptic regularity, we may extract a subsequence for which
the corresponding functions \(\phi\) converge uniformly on compact
subsets of \(\Omega\setminus\{0\}\), and each fixed \(c_{k,\nu}\) converges.
The bounds for \(f_0\) and \(f_{k,i}\) in
Lemma~\ref{lem:ce-regular-source} give uniformly,
for small \(s>0\),
\begin{equation}\label{eq:ce-source-uniform-integrability}
\int_{|x|<s}|x|^{2-n}|\Delta^2\phi|
 \leq C\bigl(b^{\frac{n+4}{n-3}}+b^{\frac{n+4}{n-4}}\bigr)s^2 +Cb^{\frac{n+2}{n-2}}
   \sum_{r_j<Cs}r_j^{2+\alpha+\frac{4(1+\alpha)}{n-2}}
 \to0.
\end{equation}
The limit in \eqref{eq:ce-source-uniform-integrability}
allows us to pass to the limit in
\(\phi=\mathcal T(\Delta^2\phi)\) and in each orthogonality
condition in \eqref{eq:ce-constrained-correction}. Since
\(\|\Delta^2\phi\|_{**}<\infty\), \eqref{eq:ce-potential-energy-bound}
gives \(\phi\in\mathcal H\). The sums in the definition of \(\|\cdot\|_*\) in
\eqref{eq:ce-linear-norms} converge uniformly in \(x\) as
\((\mathbf p,\mathbf r)\) converges in \(Y\), since \(\sum_{k\geq0}\gamma^{k(2+\alpha)}<\infty\).
Local uniform convergence of \(\phi\) and convergence of each
\(c_{k,\nu}\) imply that the limit remains in the small ball
\(\|\phi\|_*+\|\mathbf c\|_\infty\leq4C\varepsilon_0\).
By \eqref{eq:ce-nonlinear-closed-set} and the local uniform
convergence of \(U_{\lambda,r}\) and \(\phi\), we obtain
\(|\phi|\leq \frac{U_{\lambda,r}}{2}\). In particular,
\(U_{\lambda,r}+\phi\geq \frac{U_{\lambda,r}}{2}>0\) by
\eqref{eq:ce-profile-positive}.
Uniqueness in Proposition~\ref{lem:ce-projected} identifies \(\phi\),
so each fixed component of \(E_b\) is continuous.
Since \(\rho<1\) and \(\frac{\gamma^\beta}{\rho}<1\), the right-hand side of
\eqref{eq:ce-reduced-remainder-bound} tends
to zero uniformly as \(k\to\infty\). Thus \(E_b(\lambda,r)\in\mathcal X_\rho\).
The truncation of \(E_b(\lambda,r)\) to \(0\leq k\leq K\), followed by
\(K\to\infty\), gives continuity and relative compactness in
\(\|\cdot\|_{\mathcal X_\rho}\).
\end{proof}

\subsection{Choice of \texorpdfstring{\(\lambda,r\)}{lambda, r}}\label{subsec:ce-parameter-choice}
Fix the right inverse \(\mathcal S\) from Lemma~\ref{lem:ce-sequence}.
To solve \(\mathcal R(\lambda,r)=0\), we seek \(\lambda,r\)
satisfying
\begin{equation}\label{eq:ce-parameter-fixed-point}
 \begin{pmatrix}\lambda_k\\r_k\end{pmatrix}
 =
 \begin{pmatrix}\lambda_k^0\\\gamma^k\end{pmatrix}
 -\bigl(\mathcal S E_b(\lambda,r)\bigr)_k,
 \,\  k\geq-3.
\end{equation}
\begin{lemma}\label{lem:ce-recovery}
For every fixed integer \(n\geq5\) and all small \(b>0\),
there exist \(\lambda_k,r_k>0\), \(k\geq-3\), satisfying
\eqref{eq:ce-parameter-space} such that \(c_{k,1}=c_{k,2}=0\)
for every \(k\geq0\). For this choice of \(\lambda,r\),
\[
 u=U_{\lambda,r}+\phi>0,
 \,\ \Delta^2u=u^{\frac{n+4}{n-4}}\quad\text{in }\Omega\setminus\{0\}.
\]
In \(h_{r,\lambda}\), only \(\lambda_{-j},r_{-j}\), \(j=1,2,3\),
vary; the terms with \(j\geq4\) are fixed as in
Definition~\ref{def:ce-exterior-sum}.
\end{lemma}
\begin{proof}
Fix small \(0<\varepsilon\leq\varepsilon_0\) and consider the closed ball
\begin{equation}\label{eq:ce-parameter-ball}
 \left\{(\mathbf p,\mathbf r)\in Y:
       \|(\mathbf p,\mathbf r)\|_Y\leq\frac{\varepsilon}{2}\right\}.
\end{equation}
For \((\mathbf p,\mathbf r)\) in \eqref{eq:ce-parameter-ball}, the parameters
\(\lambda_k=\lambda_k^0+\mathbf p_k\) and \(r_k=\gamma^k+\mathbf r_k\)
are positive. By the mean value theorem,
\[
 \left|\left(\frac{\lambda_k}{\lambda_k^0}\right)^{\frac{n-4}{2}}-1\right|
  +\left|\frac{r_k}{\gamma^k}-1\right|
 \leq2\left(\frac{n-4}{2}\frac{|\mathbf p_k|}{\lambda_k^0}
                  +\frac{|\mathbf r_k|}{\gamma^k}\right)
 \leq\varepsilon\rho^k,\,\  k\geq-3.
\]
For fixed small \(\varepsilon>0\) and \(b>0\) small enough,
Lemmas~\ref{lem:ce-sequence} and~\ref{lem:ce-reduced} give
\[
 \|\mathcal S E_b(\lambda,r)\|_Y
 \leq C\bigl(\|(\mathbf p,\mathbf r)\|_Y^2+o_b(1)\bigr)
 \leq C\left(\frac{\varepsilon^2}{4}+o_b(1)\right)
 \leq\frac{\varepsilon}{2}.
\]
The map \((\mathbf p,\mathbf r)\mapsto-\mathcal S E_b(\lambda,r)\)
is continuous, has relatively compact image in \(Y\), and maps the
closed ball \eqref{eq:ce-parameter-ball} to itself.
The Schauder fixed point theorem gives \(\lambda,r\)
satisfying \eqref{eq:ce-parameter-fixed-point}.
At the fixed point,
\[
 (\mathbf p,\mathbf r)=-\mathcal S E_b(\lambda,r),
 \,\ \mathcal R(\lambda,r)
 =-L\mathcal S E_b(\lambda,r)+E_b(\lambda,r)=0.
\]
Lemma~\ref{lem:ce-local-identities} gives
\(c_{k,1}=c_{k,2}=0\) for every \(k\geq0\), and
\eqref{eq:ce-constrained-correction} becomes the equation
in \eqref{eq:biharmonic-problem}.

The fixed point estimate also gives
\[
 \|(\mathbf p,\mathbf r)\|_Y
 \leq\frac{C\varepsilon}{2}\|(\mathbf p,\mathbf r)\|_Y+o_b(1).
\]
For the fixed small \(\varepsilon\), absorbing
\(\frac{C\varepsilon}{2}\|(\mathbf p,\mathbf r)\|_Y\) gives
\(\|(\mathbf p,\mathbf r)\|_Y\to0\).
Condition~\eqref{eq:ce-parameter-space} therefore holds with
\(\varepsilon=C\|(\mathbf p,\mathbf r)\|_Y\to0\).
Proposition~\ref{lem:ce-projected} yields
\(\|\phi\|_*+\|\mathbf c\|_\infty\to0\),
\(u=U_{\lambda,r}+\phi\geq \frac{U_{\lambda,r}}{2}>0\), and
\(u\in C^\infty(\Omega\setminus\{0\})\).
Finally, restrict \(u\) to \(\mathbf B_1\setminus\{0\}\).
\end{proof}

\section{Proof of Theorem~\ref{thm:counterexample}}\label{sec:ce-final}

It remains to establish the behavior near zero of the solution \(u\) obtained in
Lemma~\ref{lem:ce-recovery}.
\begin{proposition}\label{prop:ce-existence}
For each fixed integer \(n\geq5\) and every even \(\ell\) large enough
for Lemmas~\ref{lem:ce-certificates} and \ref{lem:ce-midpoint-positive},
fix \(\gamma,\alpha\) as in Lemma~\ref{lem:ce-certificates},
\(1<R<\gamma^{-1}\), and \(\rho\) as in Section~\ref{app:parameter-surjectivity}.
There exists \(b_0>0\) such that, for every \(0<b<b_0\), the construction gives a positive function
\(u\in C^\infty(\mathbf B_1\setminus\{0\})\) satisfying \eqref{eq:biharmonic-problem}.
The singularity at zero is nonremovable, and
\[
 \liminf_{x\to0}u(x)=0,\,\  \limsup_{x\to0}u(x)=\infty.
\]
Moreover,
\[
 \int_{\mathbf B_1}|x|^{2-n}u(x)^{\frac{n+4}{n-4}}\,dx<\infty,
 \,\ \lim_{s\to0^+}(-\Delta u)(se_n)
   =-2(n-2)b^{\frac{n-4}{n-3}}<0,
\]
where \(e_n=(0,\ldots,0,1)\).
\end{proposition}

\begin{proof}
Fix \(n,\ell,\gamma,\alpha,R,\rho\) as in the statement, and let
\(u=U_{\lambda,r}+\phi\) be the solution obtained in
Lemma~\ref{lem:ce-recovery}. The constants in this proof are independent of \(k\).
We work in the fixed domain
\(\Omega\supset\overline{\mathbf B_1}\).

Positivity and smoothness follow from Lemma~\ref{lem:ce-recovery}.

Evaluating at \(\mathbf x_{k,1}\), \eqref{eq:ce-background-taylor} and
Proposition~\ref{lem:ce-projected} give, as \(k\to\infty\),
\[
 u(\mathbf x_{k,1})
 =A_n\lambda_k^{-\frac{n-4}{2}}(1+o(1))\to\infty.
\]
The bounds \eqref{eq:ce-parameter-space} imply
\(r_k\to0\) and \(\lambda_k\to0\), so
\(\mathbf x_{k,1}\in\mathbf B_1\) for all large \(k\).
Thus \(\limsup_{x\to0}u(x)=\infty\), and the singularity at zero is nonremovable.

On \(\mathbf B_{\eta r_k}(\mathbf x_{k,1})\), the weight \(|x|^{2-n}\)
is comparable to \(r_k^{2-n}\). The estimate
\eqref{eq:ce-local-source-integrals} with \(j=0\) and
\(\lambda_k^{\frac{n-4}{2}}\leq Cbr_k^{n-2+\alpha}\), followed by symmetry
and summation over the circles, give
\[
 \int_{\mathbf B_1}|x|^{2-n}u^{\frac{n+4}{n-4}}\,dx
 \leq Cb\left(1+\sum_{k\geq0}r_k^{\alpha}\right)<\infty.
\]
The sum \(\sum_{k\geq0}r_k^\alpha\) converges because \(\alpha>0\).
Outside \(\bigcup_{k,i}\mathbf B_{\eta r_k}(\mathbf x_{k,i})\),
\(|U_{\lambda,r}(x)|\leq Cb^{\frac{n-4}{n-3}}|x|^2\) and
\(|\phi(x)|\leq Cb|x|^{2+\alpha}\), by \eqref{eq:ce-profile-upper}
and \eqref{eq:ce-nonlinear-bound}, give a finite integral of
\(|x|^{2-n}u^{\frac{n+4}{n-4}}\).

On the ray \(x=se_n\), \(s>0\), it suffices by symmetry
to take \(i=1\).
Direct computation gives
\begin{equation}\label{eq:ce-axis-bubble-hessian}
 \sup_{0\leq s\leq s_0}
 \left|D^2\bigl[(I-J_2)\mathcal B_{\lambda_k,\mathbf x_{k,1}}\bigr](se_n)\right|
 \leq C\lambda_k^{\frac{n-4}{2}}r_k^{2-n}\leq Cbr_k^\alpha.
\end{equation}
Estimate~\eqref{eq:ce-axis-bubble-hessian} holds with
\(\mathbf x_{k,i}\) in place of \(\mathbf x_{k,1}\) for every \(i\).
Since \(\ell\) is fixed, \(\sum_k r_k^\alpha<\infty\), and
\(D^2[(I-J_2)\mathcal B_{\lambda_k,\mathbf x_{k,i}}](se_n)\to0\)
for each fixed \(k,i\), dominated convergence applies to the
series defining \(\mathbf h_{r,\lambda}\) and gives
\[
 D^2\bigl(U_{\lambda,r}-b^{\frac{n-4}{n-3}}|x''|^2\bigr)(se_n)
 \to0\,\ (s\to0^+).
\]
The function \(h_{r,\lambda}\) is smooth near zero and
\(J_2h_{r,\lambda}=0\), so \(D^2h_{r,\lambda}(se_n)\to0\).

For \(\phi=\mathcal T(f_0+\sum f_{k,i})\), with \(f_0,f_{k,i}\)
from \eqref{eq:ce-correction-source-decomposition},
\eqref{eq:ce-local-source-mass} gives
\(\sum_{k,i}r_k^{2-n}\int|f_{k,i}|<\infty\). On the support of each \(f_{k,i}\),
\[
 |D^2G_n(se_n-y)-D^2G_n(-y)|\leq Cr_k^{2-n}.
\]
The ray \(\{se_n:s>0\}\) is separated from
\(\operatorname{supp}f_{k,i}\), and
\(D^2G_n(se_n-y)-D^2G_n(-y)\to0\) for each fixed \(y\ne0\). Dominated convergence therefore
gives \(D^2\sum\mathcal Tf_{k,i}(se_n)\to0\).
The function \(f_0\), extended by zero outside \(\Omega\),
is bounded for fixed \(b\). Since \(D^2G_n\) is locally integrable,
\(D^2(G_n*f_0)\) is continuous at zero. Hence
\[
 D^2\mathcal Tf_0(se_n)
 =D^2(G_n*f_0)(se_n)-D^2(G_n*f_0)(0)\to0.
\]
Thus \(D^2\phi(se_n)\to0\).

To estimate \((I-J_2)\mathcal B_{\lambda_k,\mathbf x_{k,1}}(se_n)\),
we differentiate up to order two and obtain
\begin{equation}\label{eq:ce-axis-bubble-value}
 \sup_{0\leq s\leq s_0}
 |(I-J_2)\mathcal B_{\lambda_k,\mathbf x_{k,1}}(se_n)|
 \leq Cb(1+s_0+s_0^2)r_k^\alpha.
\end{equation}
Estimate~\eqref{eq:ce-axis-bubble-value} holds with
\(\mathbf x_{k,i}\) in place of \(\mathbf x_{k,1}\) for every \(i\),
and \((I-J_2)\mathcal B_{\lambda_k,\mathbf x_{k,i}}(se_n)\to0\)
for fixed \(k,i\). Using \eqref{eq:ce-axis-bubble-value},
the smooth function \(h_{r,\lambda}\) and
\(b^{\frac{n-4}{n-3}}s^2\to0\), we obtain \(U_{\lambda,r}(se_n)\to0\). Combining this with
\(|\phi(se_n)|\leq bs^{2+\alpha}\), we obtain \(u(se_n)\to0\).
Since \(u>0\), its lower limit at zero is zero. Finally,
\[
 \Delta\bigl(b^{\frac{n-4}{n-3}}|x''|^2\bigr)=2(n-2)b^{\frac{n-4}{n-3}},
 \,\ \lim_{s\to0^+}(-\Delta u)(se_n)
 =-2(n-2)b^{\frac{n-4}{n-3}}.
\]
All the assertions remain valid after restriction to \(\mathbf B_1\).
\end{proof}

Fixing any admissible parameters in Proposition~\ref{prop:ce-existence}
gives the solution asserted in Theorem~\ref{thm:counterexample}.

\appendix
\section{Geometry satisfying the two balance equations}
\label{app:balanced-geometry}
We prove the balance assertion of Lemma~\ref{lem:ce-certificates}
and the positivity estimates used in Lemmas~\ref{lem:ce-midpoint-positive}
and \ref{lem:ce-trace}.
Fix \(n\geq5\), let \(\ell\) tend to infinity through even integers,
and set
\[
 \kappa:=-\ell\log\gamma,\,\  \kappa\in[2\pi,\kappa_n^*].
\]
The endpoint \(\kappa_n^*>2\pi\) is supplied by Lemma~\ref{lem:large-ell-endpoints}.
All estimates in this appendix are uniform for fixed \(n\) and
\(\kappa\in[2\pi,\kappa_n^*]\). We use \(H\) from
\eqref{eq:ce-normalized-background} and \(\mathcal F\) from
Section~\ref{app:parameter-surjectivity}, with
\[
 \frac1\ell\sum_{i=1}^{\ell}\mathcal K_{\xi_i}(s\xi_1)
 =\mathcal F(s)+1+\frac{(n-4)^2}{4}s^2.
\]
Throughout this appendix, \(H\) is evaluated at \(\gamma=e^{-\frac{\kappa}{\ell}}\);
its dependence on \(\ell,\kappa,\alpha\) is suppressed.

\subsection{An inequality for \(\vartheta_0\) and \(\vartheta_{1/2}\)}
Define
\[
 \vartheta_0(t)=\sum_{j\in\mathbb Z}e^{-tj^2},
 \,\ \vartheta_{1/2}(t)=\sum_{j\in\mathbb Z}e^{-t(j+\frac{1}{2})^2},
 \,\  t>0.
\]
The series defining \(\vartheta_0,\vartheta_{1/2}\) and all
their derivatives converge locally uniformly.

\begin{lemma}\label{lem:large-ell-theta}
For \(A\geq B>0\),
\begin{equation}\label{eq:large-ell-theta-positive}
 \vartheta_{1/2}(A)\vartheta_{1/2}(B)
       +2A\vartheta_0'(A)\vartheta_0(B)>0.
\end{equation}
\end{lemma}
\begin{proof}
First take \(A=B=t\). If \(t\geq4\), the defining series give
\[
 \vartheta_0(t)\leq\frac{1+e^{-t}}{1-e^{-t}},
 \,\ -\vartheta_0'(t)\leq\frac{2e^{-t}(1+e^{-t})}{(1-e^{-t})^3},
 \,\ \vartheta_{1/2}(t)\geq2e^{-\frac{t}{4}}.
\]
Since \(e^{-t}<\frac{1}{16}\),
\(\frac{(1+e^{-t})^2}{(1-e^{-t})^4}<\frac{3}{2}\). Also \(\frac{e^{\frac{t}{2}}}{t}\) increases
for \(t\geq4\), and \(e^2>6\). Consequently
\[
 -2t\vartheta_0'(t)\vartheta_0(t)
      <6te^{-t}<4e^{-\frac{t}{2}}\leq\vartheta_{1/2}(t)^2.
\]
For \(0<t\leq4\), let \(s=\frac{\pi^2}{t}>2\). Poisson summation gives
\begin{equation}\label{eq:large-ell-theta-poisson}
 \vartheta_0(t)=\sqrt{\frac{\pi}{t}}\,\vartheta_0(s),
 \,\ \vartheta_{1/2}(t)=\sqrt{\frac{\pi}{t}}
                         \sum_{j\in\mathbb Z}(-1)^je^{-sj^2}.
\end{equation}
Subtracting \(\sqrt{\frac{t}{\pi}}\,\vartheta_{1/2}(t)\) from
\(\vartheta_0(s)\) gives
\[
0<\vartheta_0(s)-\sqrt{\frac{t}{\pi}}\,\vartheta_{1/2}(t)
   =4\sum_{j\geq1,\ j\ {\rm odd}}e^{-sj^2} \leq-2\vartheta_0'(s).
\]
Also \(\vartheta_0(s)+\sqrt{\frac{t}{\pi}}\,\vartheta_{1/2}(t)
<2\vartheta_0(s)\). Multiplying the sum and difference gives
\[
 \frac t\pi\vartheta_{1/2}(t)^2-\vartheta_0(s)^2
                         >4\vartheta_0(s)\vartheta_0'(s).
\]
Differentiating \(\vartheta_0(t)=\sqrt{\frac{\pi}{t}}\,\vartheta_0(\frac{\pi^2}{t})\)
therefore gives
\[
\begin{aligned}
&\frac t\pi\bigl[\vartheta_{1/2}(t)^2
                        +2t\vartheta_0'(t)\vartheta_0(t)\bigr] ={}\frac t\pi\vartheta_{1/2}(t)^2-\vartheta_0(s)^2
                         -2s\vartheta_0'(s)\vartheta_0(s)\\
  &\quad>(4-2s)\vartheta_0'(s)\vartheta_0(s)>0,
\end{aligned}
\]
because \(s>2\) and \(\vartheta_0'(s)<0\).

Using \eqref{eq:large-ell-theta-poisson} and the Jacobi triple product formula, we obtain
\begin{equation}\label{eq:large-ell-theta-quotient}
 \frac{\vartheta_{1/2}(t)}{\vartheta_0(t)}
  =\prod_{j\geq1}
   \left(\frac{1-e^{-\frac{(2j-1)\pi^2}{t}}}
                   {1+e^{-\frac{(2j-1)\pi^2}{t}}}\right)^2.
\end{equation}
Each factor in \eqref{eq:large-ell-theta-quotient} is positive
and strictly decreasing in \(t\).
For \(0<t_1<t_2\), the ratio of the partial products in
\eqref{eq:large-ell-theta-quotient} at \(t_1\) and \(t_2\) is at least
the ratio of their factors with \(j=1\), which is strictly greater than one.
Taking limits in the products in
\eqref{eq:large-ell-theta-quotient} shows that
\(\frac{\vartheta_{1/2}(t)}{\vartheta_0(t)}\) is strictly decreasing.
Thus for \(A\geq B\),
\[
 \frac{\vartheta_{1/2}(A)}{\vartheta_0(A)}
 \frac{\vartheta_{1/2}(B)}{\vartheta_0(B)}
 \geq\left(\frac{\vartheta_{1/2}(A)}{\vartheta_0(A)}\right)^2
       >-\frac{2A\vartheta_0'(A)}{\vartheta_0(A)}.
\]
Multiplication by \(\vartheta_0(A)\vartheta_0(B)>0\) proves
\eqref{eq:large-ell-theta-positive}, including \(A=B\).
\end{proof}

\subsection{The sums \(E_n,W_n\) and their finite parts}
For \(n>6\), \(\kappa>0\) and \(x\in\mathbb R\), define the
absolutely convergent sums
\begin{equation}\label{eq:large-ell-lattice}
\begin{split}
 E_n(\kappa)&=\sum_{(k,j)\in\mathbb Z^2\setminus\{(0,0)\}}
       (\kappa^2k^2+4\pi^2j^2)^{-\frac{n-4}{2}},\\
 W_n(x;\kappa)&=\sum_{k,j\in\mathbb Z}
       [(x-\kappa k)^2+4\pi^2(j+\frac{1}{2})^2]^{-\frac{n-4}{2}}.
\end{split}
\end{equation}
The function \(W_n(x;\kappa)\) is smooth and periodic in \(x\), with period \(\kappa\).
A prime on \(E_n\) denotes differentiation with respect to \(\kappa\).
The sums in \eqref{eq:large-ell-lattice} diverge for \(n=5,6\);
we define \(E_n,W_n\) in these dimensions by continuation.

For the finite parts at \(n=5,6\), regard \(n>4\) as a real
variable for continuation; the geometric assertions concern integer
\(n\geq5\). Mellin transformation and Poisson summation give
\begin{equation}\label{eq:large-ell-mellin-E}
\begin{split}
 E_n(\kappa)=\frac1{\Gamma(\frac{n-4}{2})}\bigg\{&
 \int_0^1t^{\frac{n-6}{2}}
  \left[\vartheta_0(\kappa^2t)\vartheta_0(4\pi^2t)
                         -1-\frac1{2\kappa t}\right]dt\\
 &+\int_1^\infty t^{\frac{n-6}{2}}
       [\vartheta_0(\kappa^2t)\vartheta_0(4\pi^2t)-1]dt
                         +\frac1{\kappa(n-6)}\bigg\}.
\end{split}
\end{equation}
Here \(\Gamma\) is the Gamma function.
For \(W_n\), replace
\(\vartheta_0(\kappa^2t)\vartheta_0(4\pi^2t)-1\)
in both integrals in \eqref{eq:large-ell-mellin-E} by
\[
 \sum_k e^{-t(x-\kappa k)^2}\vartheta_{1/2}(4\pi^2t),
\]
and use the subtraction \(\frac{1}{2\kappa t}\) and the term
\(\frac{1}{\kappa(n-6)}\). The integrals are locally uniformly convergent for
\(n>4\). As \(t\to0\),
\[
\begin{aligned}
 \vartheta_0(\kappa^2t)\vartheta_0(4\pi^2t)-1-\frac1{2\kappa t}
 &=-1+O(t^{-1}e^{-\frac{c}{t}}),\\
 \sum_k e^{-t(x-\kappa k)^2}\vartheta_{1/2}(4\pi^2t)-\frac1{2\kappa t}
 &=O(t^{-1}e^{-\frac{c}{t}}).
\end{aligned}
\]
For \(n>6\), adding back the subtracted integral recovers
\eqref{eq:large-ell-lattice}.

At \(n=5\) we evaluate \eqref{eq:large-ell-mellin-E} and the
corresponding formula for \(W_n\). At \(n=6\), we denote the coefficients
of \((n-6)^0\) in the expansions below by \(E_6,W_6\). The expansions are
\[
 E_n=\frac1{\kappa(n-6)}+E_6+O(n-6),
 \,\ W_n=\frac1{\kappa(n-6)}+W_6+O(n-6).
\]
In particular \(\partial_\kappa(\kappa E_n)\) has a removable
singularity at dimension six and its value is
\(\partial_\kappa(\kappa E_6)\).

We use Lemma~\ref{lem:large-ell-theta} to obtain a lower bound for
\(W_n(x;\kappa)+\kappa E_n'(\kappa)\).
\begin{lemma}\label{lem:large-ell-lattice-positive}
For each fixed integer \(n\geq5\) and \(\kappa\geq2\pi\),
\[
 W_n(x;\kappa)+\kappa E_n'(\kappa)>0,
 \,\  x\in\mathbb R.
\]
On every compact interval \([2\pi,K]\),
\(W_n(x;\kappa)+\kappa E_n'(\kappa)\) has a positive lower bound
independent of \(x\) and \(\kappa\).
\end{lemma}
\begin{proof}
For \(n>6\), applying Mellin transformation to the series gives
\begin{equation}\label{eq:large-ell-joint-mellin}
\begin{split}
 W_n(x;\kappa)+\kappa E_n'(\kappa)
  =\frac1{\Gamma(\frac{n-4}{2})}\int_0^\infty t^{\frac{n-6}{2}}
  \bigg[&\sum_k e^{-t(x-\kappa k)^2}\vartheta_{1/2}(4\pi^2t)\\
       &+2\kappa^2t\vartheta_0'(\kappa^2t)
                         \vartheta_0(4\pi^2t)\bigg]dt.
\end{split}
\end{equation}
The contributions of \(\frac{1}{2\kappa t}\) to \(W_n\) and
\(\kappa E_n'\) cancel as \(t\to0^+\). In \eqref{eq:large-ell-joint-mellin},
\[
 \sum_k e^{-t(x-\kappa k)^2}\vartheta_{1/2}(4\pi^2t)
 +2\kappa^2t\vartheta_0'(\kappa^2t)\vartheta_0(4\pi^2t)
 =O(t^{-2}e^{-\frac{c}{t}})
\]
as \(t\to0^+\), locally uniformly in \(\kappa\) and uniformly in \(x\). As \(t\to\infty\), the factors
\(\vartheta_{1/2}(4\pi^2t)\) and \(\vartheta_0'(\kappa^2t)\)
decay exponentially.
Consequently the integral converges for every
\(n>4\). Subtracting \(\frac{1}{2\kappa t}\) for both sums in
\eqref{eq:large-ell-mellin-E} shows that
\eqref{eq:large-ell-joint-mellin} equals \(W_n+\kappa E_n'\)
also for the specified values at \(n=5,6\).

For every \(x\in\mathbb R\),
\[
 \sum_k e^{-t(x-\kappa k)^2}\geq
                \vartheta_{1/2}(\kappa^2t).
\]
Poisson summation and the Jacobi triple product formula give this bound: each factor
that depends on \(x\) has the form
\(1+2a\cos(\frac{2\pi x}{\kappa})+a^2\), with \(0<a<1\),
and is minimized when \(\cos(\frac{2\pi x}{\kappa})=-1\).
Taking limits in the convergent positive products preserves the bound.
Substitution in \eqref{eq:large-ell-joint-mellin} bounds
\(\sum_k e^{-t(x-\kappa k)^2}\vartheta_{1/2}(4\pi^2t)
 +2\kappa^2t\vartheta_0'(\kappa^2t)\vartheta_0(4\pi^2t)\) from below by
\begin{equation}\label{eq:large-ell-positive-integrand}
 \vartheta_{1/2}(\kappa^2t)\vartheta_{1/2}(4\pi^2t)
       +2\kappa^2t\vartheta_0'(\kappa^2t)
                                  \vartheta_0(4\pi^2t)>0
\end{equation}
by Lemma~\ref{lem:large-ell-theta}. On a compact interval of \(\kappa\),
the expression in \eqref{eq:large-ell-positive-integrand}
has a positive minimum for \(1\leq t\leq2\).
Integration over \([1,2]\) gives the uniform positive bound in
Lemma~\ref{lem:large-ell-lattice-positive}.
\end{proof}

To solve \(H'(1)=0\) along the branch \(H(1)=0\), we determine the signs of
\(E_n(\kappa)+\kappa E_n'(\kappa)\) at \(\kappa=2\pi\) and at a large value of
\(\kappa\).
\begin{lemma}\label{lem:large-ell-endpoints}
For every fixed integer \(n\geq5\), there is a finite
\(\kappa_n^*>2\pi\) such that
\begin{equation}\label{eq:large-ell-endpoint-signs}
 \left.\partial_\kappa(\kappa E_n)\right|_{\kappa=2\pi}<0,
 \,\ \left.\partial_\kappa(\kappa E_n)\right|_{\kappa=\kappa_n^*}>0.
\end{equation}
\end{lemma}
\begin{proof}
Applying \(1+\kappa\partial_\kappa\) to \eqref{eq:large-ell-mellin-E}
cancels both \(\frac{1}{\kappa(n-6)}\) and the integral of the subtraction \(\frac{1}{2\kappa t}\), giving
\begin{equation}\label{eq:large-ell-endpoint-mellin}
E_n+\kappa E_n'
 =\frac1{\Gamma(\frac{n-4}{2})}\int_0^\infty t^{\frac{n-6}{2}}
 \bigl[\vartheta_0(\kappa^2t)\vartheta_0(4\pi^2t) +2\kappa^2t\vartheta_0'(\kappa^2t)\vartheta_0(4\pi^2t)-1\bigr]dt.
\end{equation}
Formula~\eqref{eq:large-ell-endpoint-mellin} also holds for the finite parts
at \(n=5,6\), since
\[
 \vartheta_0(\kappa^2t)\vartheta_0(4\pi^2t)
 +2\kappa^2t\vartheta_0'(\kappa^2t)\vartheta_0(4\pi^2t)-1
 =-1+O(t^{-2}e^{-\frac{c}{t}})
\]
as \(t\to0\), while
\(\vartheta_0(\kappa^2t)\vartheta_0(4\pi^2t)-1\) and
\(2\kappa^2t\vartheta_0'(\kappa^2t)\vartheta_0(4\pi^2t)\)
decay exponentially as \(t\to\infty\). Poisson summation, with \(s=\frac{\pi^2}{A}\), gives
\[
 \vartheta_0(A)+2A\vartheta_0'(A)
       =-2s\sqrt{\frac{\pi}{A}}\,\vartheta_0'(s)>0.
\]
At \(\kappa=2\pi\), put \(A=4\pi^2t\) in
\eqref{eq:large-ell-endpoint-mellin}. For \(A\geq1\),
\(\vartheta_0(A)^2+A(\vartheta_0(A)^2)'-1\) equals
\[
 \sum_{z\in\mathbb Z^2\setminus\{0\}}
                 (1-A|z|^2)e^{-A|z|^2}<0.
\]
For \(0<A\leq1\), Poisson summation gives
\(-\frac{s^2(\vartheta_0(s)^2)'}{\pi}-1\), where \(s\geq\pi^2>9\).
The elementary series bound
\[
 -\frac{s^2}{\pi}(\vartheta_0(s)^2)'
 \leq\frac{4s^2e^{-s}(1+e^{-s})^2}{\pi(1-e^{-s})^4}<1
\]
proves negativity also in this range. Thus
\eqref{eq:large-ell-endpoint-mellin} is negative at \(2\pi\).
As \(\kappa\to\infty\),
\[
 \vartheta_0(\kappa^2t)\vartheta_0(4\pi^2t)
 +2\kappa^2t\vartheta_0'(\kappa^2t)\vartheta_0(4\pi^2t)-1
 \to\vartheta_0(4\pi^2t)-1>0.
\]
Apply dominated convergence to the negative part of
\(\vartheta_0(\kappa^2t)\vartheta_0(4\pi^2t)
 +2\kappa^2t\vartheta_0'(\kappa^2t)\vartheta_0(4\pi^2t)-1\),
with weight \(t^{\frac{n-6}{2}}\). This negative part is bounded by
\(1\) for \(0<t\leq1\) and by \(Ct e^{-4\pi^2t}\) for \(t\geq1\),
uniformly for \(\kappa\geq2\pi\).
The bounds \(1\) on \((0,1]\) and \(Ct e^{-4\pi^2t}\) on
\([1,\infty)\) are integrable with weight \(t^{\frac{n-6}{2}}\), so the
integral of the negative part tends to zero. The positive contribution on
\(1\leq t\leq2\) tends to a strictly positive integral. Hence
\eqref{eq:large-ell-endpoint-mellin} is positive for all
large \(\kappa\), proving \eqref{eq:large-ell-endpoint-signs}.
\end{proof}

\subsection{Two uniform asymptotic formulas}
Write
\[
 P_\ell(x)=\mathcal P_n(e^{-\frac{x}{\ell}}),\,\  r=e^{-\frac{x}{\ell}}.
\]
We use the range
\begin{equation}\label{eq:large-ell-scales}
 0<\alpha\leq\bar\alpha=1+\frac3n<2,
 \,\ 0\leq x\leq\kappa.
\end{equation}
If the radii have ratio \(s\geq0\) and the angle between the two points is \(\theta\), set
\[
\begin{split}
 k_n(s,\theta)&=(1+s^2-2s\cos\theta)^{-\frac{n-4}{2}},\\
 T_{2,n}(s,\theta)&=1+(n-4)s\cos\theta
   +\left[\frac{(n-4)^2}{4}
       +\frac{(n-4)(n-2)}4\cos(2\theta)\right]s^2.
\end{split}
\]
The cubic Taylor term \(T_{3,n}(s,\theta)\) of
\(k_n(s,\theta)\) at \(s=0\) is a linear combination of \(\cos\theta\)
and \(\cos(3\theta)\). Its averages at the angles
\(\frac{2\pi j}{\ell}\) and \(\frac{(2j+1)\pi}{\ell}\), \(0\leq j<\ell\),
and its continuous angular average vanish for large \(\ell\).
Choose a smooth \(\chi_-\) equal to one for \(t\leq-2\) and zero for
\(t\geq-1\), and put
\[
 f_n(t,\theta)=e^{-(2+\alpha)t}
 [k_n(e^t,\theta)-T_{2,n}(e^t,\theta)
                         -\chi_-(t)T_{3,n}(e^t,\theta)].
\]
All angular sums below are taken over one period, \(-\frac{\ell}{2}\leq j<\frac{\ell}{2}\),
and a prime omits \((k,j)=(0,0)\). Direct substitution gives
\begin{equation}\label{eq:large-ell-low-d-self}
\begin{aligned}
H(1)&=\frac1\ell\sum_{k,j}'
 f_n\left(\frac{\kappa k}{\ell},\frac{2\pi j}{\ell}\right)
                             -\frac{T_{2,n}(1,0)}\ell,\\
 r^{-(2+\alpha)}P_\ell(x)&=\frac1\ell\sum_{k,j}
 f_n\left(\frac{\kappa k-x}{\ell},\frac{2\pi(j+\frac{1}{2})}\ell\right),\\
 H'(1)-(2+\alpha)H(1)&=\frac1\ell\sum_{k,j}'
 \partial_t f_n\left(\frac{\kappa k}{\ell},\frac{2\pi j}{\ell}\right) +\frac{n-2+\alpha T_{2,n}(1,0)}\ell.
\end{aligned}
\end{equation}
Here \(T_{2,n}(1,0)=\frac{(n-3)(n-2)}{2}\). The grid shifted by half an angular spacing is at distance at
least \(\frac{\pi}{\ell}\) from the singular point, including \(x=0,\kappa\).

We use \eqref{eq:large-ell-low-d-self} to express the leading terms of
\(r^{-(2+\alpha)}P_\ell(x)-H(1)\) and \(H'(1)-(2+\alpha)H(1)\) in terms of \(E_n\)
and \(W_n\).
\begin{lemma}\label{lem:large-ell-transfer}
For fixed \(n\geq5\), as \(\ell\to\infty\) through even integers,
\begin{equation}\label{eq:large-ell-transfer}
\begin{split}
 \frac{r^{-(2+\alpha)}P_\ell(x)-H(1)}{\ell^{n-5}}
   &=W_n(x;\kappa)-E_n(\kappa)+o(1),\\
 \frac{H'(1)-(2+\alpha)H(1)}{\ell^{n-5}}
   &=-\left(\alpha+\frac n2\right)
                      [E_n(\kappa)+\kappa E_n'(\kappa)]+o(1).
\end{split}
\end{equation}
The errors are uniform in \eqref{eq:large-ell-scales} and for
\(\kappa\) in each fixed compact subset of \((0,\infty)\).
\end{lemma}
\begin{proof}
Put \(d=n-4\) and choose a smooth cut-off function \(\chi_+\),
zero for \(t\leq1\) and one for \(t\geq2\).
As \(t\to+\infty\), the only term in \(f_n(t,\theta)\) whose
decay is not uniform as \(\alpha\to0\) is the negative of
\[
 \left[\frac{d^2}{4}+\frac{d(d+2)}4\cos(2\theta)\right]e^{-\alpha t}.
\]
Set \(h_\alpha(t)=\chi_+(t)e^{-\alpha t}\) and
\[
 G_\alpha(t,\theta)=f_n(t,\theta)
 +\left[\frac{d^2}{4}+\frac{d(d+2)}4\cos(2\theta)\right]h_\alpha(t).
\]
The function \(G_\alpha\) and its derivatives up to order three
decay exponentially as \(t\to\pm\infty\), uniformly because \(\bar\alpha<2\).
For \(j\geq1\), \(\|h_\alpha^{(j)}\|_{L^1}\leq C_j\), uniformly
for every \(0<\alpha\leq\bar\alpha\). The estimate comparing a sum on a shifted grid with its integral
\begin{equation}\label{eq:large-ell-grid-quadrature}
 \left|\sum_k F(\Delta k+s)-\Delta^{-1}\int_{\mathbb R}F\right|
             \leq C\Delta\|F''\|_{L^1},\,\  \Delta=\frac\kappa\ell,
\end{equation}
is applied to \(h_\alpha\) at
\(t=\frac{\kappa k}{\ell}\) and \(t=\frac{\kappa k-x}{\ell}\), \(k\in\mathbb Z\);
the integral term \(\frac\ell\kappa\int_{\mathbb R}h_\alpha(t)\,dt\)
is the same for both sums and cancels in their difference. The averages of \(\cos(2\theta)\) over
\(\theta=\frac{2\pi j}{\ell}\) and \(\theta=\frac{(2j+1)\pi}{\ell}\), \(0\leq j<\ell\), vanish. The contribution of \(h_\alpha\) to
\(r^{-(2+\alpha)}P_\ell(x)-H(1)\) is therefore \(O(\ell^{-1})\).
For \(H'(1)-(2+\alpha)H(1)\),
\(\int h_\alpha'=0\) and the estimate with
\(\|h_\alpha'''\|_{L^1}\leq C\) gives the same error.
After division by \(\ell^{n-5}\), both errors are \(O(\ell^{-(n-4)})\).

Near \(z=(t,\theta)=0\), write
\[
\begin{aligned}
 K_n(t,\theta)&=[4\sinh^2(\frac{t}{2})+4\sin^2(\frac{\theta}{2})]^{-\frac{d}{2}},\\
 e^{-(2+\alpha)t}k_n(e^t,\theta)&=e^{-(\alpha+\frac{n}{2})t}K_n(t,\theta).
\end{aligned}
\]
All estimates below are uniform for \(\kappa\in I\), with \(I\) a compact subinterval of \((0,\infty)\), \(0<\alpha\leq\bar\alpha\), and \(0\leq x\leq\kappa\).
Put \(h=\ell^{-1}\), \(A=2\pi\kappa\), \(a=\alpha+\frac{n}{2}\), \(\Lambda_\kappa=\{(\kappa k,2\pi j):k,j\in\mathbb Z\}\), and \(v=(-x,\pi)\).
Let \(\Sigma_h^0,\Sigma_h^v\) denote the two grids in \eqref{eq:large-ell-low-d-self}, omitting the zero node from \(\Sigma_h^0\);
the normalized sums are \(h^d\Sigma_h^\sigma\).
On \(\mathcal C=\mathbb R\times(\mathbb R/2\pi\mathbb Z)\), choose a smooth radial cutoff \(\chi\), equal to one near zero and supported in \(|z|<\frac{\pi}{2}\).
The local expansion gives
\[
 G_\alpha=\chi|z|^{-d}+g_\alpha,
 \,\ F_\alpha:=(\partial_tG_\alpha)_{\mathrm{even}}=a\chi L_d+R_\alpha,
\]
where
\[
 aL_d(z)=\left(\alpha+\frac n2\right)\left(-|z|^{-d}+dt^2|z|^{-d-2}\right).
\]
For \(d=1,2\), near zero,
\(|D^jg_\alpha|\leq C|z|^{1-d-j}\) \((0\leq j\leq2)\), \(|R_\alpha|\leq C\), and
\(|D^2R_\alpha|\leq C(1+|z|^{-d})\); away from zero the derivatives have uniform exponential tails.

Set \(Q_h^\sigma(u)=Ah^2\Sigma_h^\sigma u-\int_{\mathcal C}u\).
For centered rectangular cells of area \(Ah^2\), interpreted periodically in \(\theta\), the \(O(1)\) cells meeting \(|z|\leq Ch\),
including the omitted central cell, contribute \(O(h^{3-d})\) for \(g_\alpha\) and \(O(h^2)\) for \(R_\alpha\).
Retained nodes are separated from zero by \(ch\), and shifted nodes by \(\pi h\).
Taylor's theorem on the other cells, whose linear terms integrate to zero, yields
\[
 \left|Q_h^\sigma(g_\alpha)\right|
 \leq C\left(h^{3-d}+h^2\int_h^1r^{-d}\,dr+h^2\right),
 \,\ \left|Q_h^0(R_\alpha)\right|
 \leq Ch^2\left(1+\int_h^1r^{1-d}\,dr\right).
\]
The normalized errors are \(\frac{h^{d-2}Q_h^\sigma}{A}\).

For the first formula, the identity
\[
 W_{d+4}(x;\kappa)-E_{d+4}(\kappa)
 =|v|^{-d}
 +\frac12\sum_{p\in\Lambda_\kappa\setminus\{0\}}
 \left(|p+v|^{-d}+|p-v|^{-d}-2|p|^{-d}\right)
\]
holds for \(d>2\), and hence for \(d=1,2\) by \eqref{eq:large-ell-mellin-E}:
the paired summands are \(O(|p|^{-d-2})\) and define a holomorphic function for \(\operatorname{Re}d>0\).
The cutoff sum differs by \(O(h^d)\), because the Hessian of \((\chi(hz)-1)|z|^{-d}\) vanishes for \(|z|<\frac{c}{h}\)
and is \(O(|z|^{-d-2})\) elsewhere. The two integrals of \(g_\alpha\) cancel, so the cell estimates prove the first formula.

For the second formula, take \(C_p=p+[-\frac{\kappa}{2},\frac{\kappa}{2}]\times[-\pi,\pi]\), let \(\rho(\varphi)\) be the radial boundary of \(C_0\), and set
\[
 J_d=\frac1{2-d}\int_0^{2\pi}(-1+d\cos^2\varphi)\rho(\varphi)^{2-d}\,d\varphi,
 \,\ \mathcal L_d=\sum_{p\ne0}\left(L_d(p)-\frac1A\int_{C_p}L_d\right)-\frac{J_d}{A}.
\]
The cell errors are \(O(|p|^{-d-2})\) and the singularity of \(J_d\) at \(d=2\) is removable.
Thus \(\mathcal L_d\) is holomorphic for \(\operatorname{Re}d>0\).
Radial integration gives \(J_d=\int_{C_0}L_d\) for \(0<d<2\) and \(J_d=-\int_{\mathbb R^2\setminus C_0}L_d\) for \(d>2\);
consequently
\[
 \mathcal L_d=\sum_{p\ne0}L_d(p)=-(1+\kappa\partial_\kappa)E_{d+4}(\kappa)\,\ (d>2).
\]
Absolute convergence permits differentiation here. By \eqref{eq:large-ell-mellin-E}, the identity continues to \(d=1,2\),
since \(1+\kappa\partial_\kappa\) annihilates the pole \(\frac{1}{\kappa(d-2)}\). Moreover,
\[
 J_2=\int_0^{2\pi}\cos(2\varphi)\log\rho(\varphi)\,d\varphi-\pi
 =\operatorname{PV}\int_{C_0}L_2-\pi.
\]
Thus the disk cell constant at \(d=2\) differs from \(\mathcal L_2\) by \(-\frac{\pi}{A}\).
Summing the centered-cell discrepancies of \((\chi(hz)-1)L_d\), whose Hessian is \(O(|z|^{-d-2})\) outside \(\frac{c}{h}\),
gives the finite-part identity
\[
 \sum_{p\ne0}\chi(hp)L_d(p)
 -\frac1A\operatorname{FP}\int_{\mathbb R^2}\chi(hz)L_d(z)\,dz
 =-(E_{d+4}+\kappa E_{d+4}')-\frac{\pi}{A}\mathbf1_{\{d=2\}}+O(h^d).
\]
Here \(\operatorname{FP}\) is the ordinary integral for \(d=1\) and the circular principal value for \(d=2\);
the latter integral is zero by angular cancellation.

On \(z=\delta(\cos\varphi,\sin\varphi)\), \(G_\alpha=\delta^{-d}-a\delta^{1-d}\cos\varphi+O(1)\).
The outward normal to the punctured cylinder has component \(-\cos\varphi\) on its inner boundary, and the end contributions vanish. Hence
\[
 \operatorname{FP}\int_{\mathcal C}F_\alpha
 =\lim_{\delta\to0^+}\left[-\delta\int_0^{2\pi}G_\alpha(\delta\cos\varphi,\delta\sin\varphi)\cos\varphi\,d\varphi\right]
 =a\pi\mathbf1_{\{d=2\}}.
\]
Combining this with the cutoff identity and the remainder estimate,
\[
 \begin{aligned}
 h^d\Sigma_h^0F_\alpha
 &=\frac{h^{d-2}}A\operatorname{FP}\int_{\mathcal C}F_\alpha
 -a(E_{d+4}+\kappa E_{d+4}')-\frac{a\pi}{A}\mathbf1_{\{d=2\}}+O(h)\\
 &=-a(E_{d+4}+\kappa E_{d+4}')+O(h),\,\  d=1,2.
 \end{aligned}
\]
The odd part cancels on the unshifted grid. By \eqref{eq:large-ell-grid-quadrature}, the \(h_\alpha\) terms and the explicit corrections
in \eqref{eq:large-ell-low-d-self} are \(O(h^d)\) after normalization. Thus the first line of \eqref{eq:large-ell-transfer}
has error \(O(\ell^{-1}(1+\log\ell))\) for \(n=5\) and \(O(\ell^{-1})\) for \(n=6\), while the second has error \(O(\ell^{-1})\)
in both cases, with the stated uniformity. For \(d>2\), the local expansions and uniform exponential tails dominate the scaled summands
by \(C(1+k^2+j^2)^{-\frac{d}{2}}\), omitting the zero node on the unshifted grid. Uniform dominated convergence gives both limits directly.

\end{proof}

\subsection{Solving \(H(1)=H'(1)=0\)}
By \eqref{eq:ce-full-circle-coefficients},
\[
\begin{gathered}
 \frac1\ell\sum_{i=1}^{\ell}\mathcal K_{\xi_i}(s\xi_1)
 =\sum_{j\geq0}c_{2j}s^{2j},\\
 c_0=1,\quad c_2=\frac{(n-4)^2}{4},\quad
 c_4=\left[\frac{(n-4)(n-2)}8\right]^2,\quad c_{2j}\geq0.
\end{gathered}
\]
The coefficients \(c_{2j}\) have at most polynomial growth in \(j\). Pairing circles at
reciprocal radii in \eqref{eq:ce-normalized-background}, and writing
\(s_k=e^{-\frac{\kappa k}{\ell}}\), gives
\begin{equation}\label{eq:large-ell-exact-tail}
\begin{aligned}
H(1)={}&\frac1\ell\sum_{j=2}^{\ell}\mathcal K_{\xi_j}(\xi_1)-1-c_2\\
 &+\sum_{k\geq1}\Bigl[
 s_k^{-2-\alpha}\mathcal F(s_k)
 +s_k^{n-2+\alpha}\bigl(\mathcal F(s_k)+1+c_2s_k^2\bigr)
 -s_k^{2+\alpha}\Bigr]\\
 &-\frac{c_2}{e^{\frac{\kappa\alpha}{\ell}}-1}.
\end{aligned}
\end{equation}

To solve \(H(1)=0\) for \(\alpha\), we first use \eqref{eq:large-ell-exact-tail}
to show that \(H(1)\) is strictly increasing in \(\alpha\).
\begin{lemma}\label{lem:large-ell-center-monotone}
For fixed even \(\ell\geq8\), \(0<\gamma<1\), and \(0<\alpha<2\),
one has \(\partial_\alpha H(1)>0\).
\end{lemma}
\begin{proof}
The function \(k_\gamma(z)=\frac{(-\log\gamma)\gamma^z}{(1-\gamma^z)^2}\)
strictly decreases on \(z>0\). Differentiating the power series after evaluating the geometric
sums gives
\begin{equation}\label{eq:large-ell-center-derivative}
\begin{split}
 \partial_\alpha H(1)={}&\sum_{m\geq4}c_m
 [k_\gamma(m-2-\alpha)-k_\gamma(m+n-2+\alpha)]\\
 &+k_\gamma(2+\alpha)-k_\gamma(n-2+\alpha)
 +c_2[k_\gamma(\alpha)-k_\gamma(n+\alpha)]>0.
\end{split}
\end{equation}
The sum runs over even indices and converges locally at an exponential rate.
Each difference of two values of \(k_\gamma\) in
\eqref{eq:large-ell-center-derivative} is positive.
\end{proof}

We first solve \(H(1)=0\) for \(\alpha\) and then use the signs in
Lemma~\ref{lem:large-ell-endpoints} to choose \(\kappa\) so that \(H'(1)=0\).
\begin{lemma}\label{lem:large-ell-root}
For fixed \(n\geq5\), every large even \(\ell\) admits
\[
 \gamma_\ell=e^{-\frac{\kappa_\ell}{\ell}},\quad
 \kappa_\ell\in[2\pi,\kappa_n^*],
 \,\ 0<\alpha_\ell<1+\frac3n,
\]
such that \(H(1)=H'(1)=0\).
\end{lemma}
\begin{proof}
Fix \(\kappa\in[2\pi,\kappa_n^*]\). In
\eqref{eq:large-ell-exact-tail}, the term
\(-\frac{c_2}{e^{\frac{\kappa\alpha}{\ell}}-1}\) gives \(H(1)\to-\infty\)
as \(\alpha\to0^+\), while \(c_4\) gives
\(\frac{c_4}{e^{\frac{\kappa(2-\alpha)}{\ell}}-1}\to+\infty\) as
\(\alpha\to2^-\). The other terms in
\eqref{eq:large-ell-exact-tail} are bounded at the respective endpoints
\(\alpha=0,2\). Lemma~\ref{lem:large-ell-center-monotone}
therefore gives a unique root \(\alpha_\ell(\kappa)\in(0,2)\)
of \(H(1)=0\).
Local smoothness and \(\partial_\alpha H(1)>0\) imply that
\(\kappa\mapsto\alpha_\ell(\kappa)\) is continuous.
Discarding nonnegative terms in \eqref{eq:large-ell-exact-tail} gives
\[
H(1)\geq{}-1-c_2+\frac{c_4}{e^{\frac{\kappa(2-\alpha)}{\ell}}-1}
 +\frac1{e^{\frac{\kappa(n-2+\alpha)}{\ell}}-1} -\frac1{e^{\frac{\kappa(2+\alpha)}{\ell}}-1}
 -\frac{c_2}{e^{\frac{\kappa\alpha}{\ell}}-1}.
\]
At \(\bar\alpha=1+\frac{3}{n}\), the coefficient of \(\frac{\ell}{\kappa}\) obeys
\[
 \frac{c_4}{2-\bar\alpha}-\frac{c_2}{\bar\alpha}
 -\frac{n-4}{(2+\bar\alpha)(n-2+\bar\alpha)}
 \geq\frac{13(n-4)^2}{128}-\frac{n-4}{3(n-1)}
 \geq\frac{7(n-4)}{384}>0.
\]
Thus \(H(1)>0\) at \(\alpha=\bar\alpha\) for large \(\ell\), uniformly
on the compact \(\kappa\) interval, and
\(0<\alpha_\ell(\kappa)<\bar\alpha\).
Along this branch, \(H(1)=0\), so the second line of
\eqref{eq:large-ell-transfer} becomes
\[
 \frac{H'(1)}{\ell^{n-5}}
 =-\left(\alpha_\ell(\kappa)+\frac n2\right)
                  \partial_\kappa(\kappa E_n)+o(1).
\]
Lemma~\ref{lem:large-ell-endpoints} gives opposite signs of
\(H'(1)\) at \(2\pi,\kappa_n^*\). Continuity gives \(\kappa_\ell\)
with \(H'(1)=0\), while \(H(1)=0\) already holds.
The bound \(\alpha_\ell<1+\frac{3}{n}\) also gives
\(2-\alpha_\ell\geq1-\frac{3}{n}>0\) and
\(n-(n-6)\alpha_\ell\geq\min\{n,3+\frac{18}{n}\}>0\).
\end{proof}

\subsection{Positivity when \(H(1)=H'(1)=0\)}
\label{app:geometric-positivity}
We justify Lemmas~\ref{lem:ce-midpoint-positive} and \ref{lem:ce-trace}
for the parameters already chosen in Lemma~\ref{lem:large-ell-root}.
Combining the two lines of \eqref{eq:large-ell-transfer}, with
\(r=e^{-\frac{x}{\ell}}\), gives
\begin{equation}\label{eq:large-ell-joint-cancellation}
\frac{(\alpha+\frac{n}{2})r^{-(2+\alpha)}P_\ell(x)
              -H'(1)-\frac{n-4}{2}H(1)}{\ell^{n-5}} =\left(\alpha+\frac n2\right)
       [W_n(x;\kappa)+\kappa E_n'(\kappa)]+o(1).
\end{equation}
At \(H(1)=H'(1)=0\), Lemma~\ref{lem:large-ell-lattice-positive}
and the uniform \(o(1)\) in
\eqref{eq:large-ell-joint-cancellation} give
\[
 \mathcal P_n(r)\geq c_n \ell^{n-5} r^{2+\alpha_\ell}>0,
 \,\  \gamma_\ell\leq r\leq1.
\]

For a fixed radius ratio \(s>0\), consider
\(f(z)=(1+s^2-2sz)^{-\frac{n-4}{2}}\) on \([-1,1]\), with
\(s=z=1\) excluded. Its Hermite interpolant at
\(\cos(\frac{(2j+1)\pi}{\ell})\), \(0\leq j<\frac{\ell}{2}\),
matches values and first derivatives, has degree
at most \(\ell-1\), and has the nonnegative remainder
\begin{equation}\label{eq:large-ell-hermite-remainder}
 \frac{f^{(\ell)}(\zeta)}{\ell!}
       \prod_{j=0}^{\frac{\ell}{2}-1}
                    [z-\cos(\frac{(2j+1)\pi}{\ell})]^2\geq0,
\end{equation}
where \(\zeta\) is an intermediate point and \(f^{(\ell)}>0\).
The average of the Hermite interpolant over a rotated regular polygon
is independent of the rotation, because every term depending on the angle
is a linear combination of \(\cos(q\theta)\) and \(\sin(q\theta)\), \(1\leq q<\ell\). At the midpoint angles the
polynomial equals \(f\), and \eqref{eq:large-ell-hermite-remainder}
shows that \(f\)
is at least as large as the polynomial elsewhere. Hence the average of \(f\) over every
rotation of the polygon is at least the average at the midpoint angles; taking the
limit as the evaluation point approaches \(\xi_1\), and then using symmetry, preserves the
inequality.
The averaged Taylor polynomial \(1+\frac{(n-4)^2s^2}{4}\) is independent
of the rotation. Subtracting it and summing the absolutely convergent
series for \(\mathcal P_n\) therefore preserves the comparison.

Shifting the index gives
\(\mathcal P_n(\gamma_\ell r)=
\gamma_\ell^{2+\alpha_\ell}\mathcal P_n(r)\).
The minimum of \(\frac{\mathcal P_n(r)}{r^{2+\alpha_\ell}}\) for
\(\gamma_\ell\leq r\leq1\) is positive, so
\begin{equation}\label{eq:large-ell-all-r-positive}
 \mathcal P_n(r)\geq c_{n,\ell}r^{2+\alpha_\ell}>0,
                         \,\  r>0.
\end{equation}
The comparison of the angular sums gives the lower bound
in \eqref{eq:large-ell-all-r-positive} at every nonsingular
planar point. The circles and the weights in
\eqref{eq:ce-midpoint-potential} are fixed before
\(b\) is chosen.

\medskip
\noindent\textbf{Remark.}
Once our work was completed, we learned that Wei, Yang and Yang
\cite{WYY} independently obtained the failure of superharmonicity near
an isolated singularity for the critical biharmonic equation
\eqref{eq:biharmonic-problem}, i.e., the result of Theorem~\ref{thm:counterexample} here.

\section*{Acknowledgements}
X. Huang is supported by the National Natural Science Foundation of China (NSFC, No.~12671132)
and the Science and Technology Commission of Shanghai Municipality (No.~22DZ2229014); Y. Li
is supported by the NSFC (grant No.~12401132).

\section*{Declaration of AI-assisted tools}

This paper was written by the authors and was not generated by artificial
intelligence. During its preparation, the authors used AI tools as auxiliary
aids for exploratory discussions, preliminary proof drafting, and language
editing. The authors take full responsibility for the mathematical validation,
the final proof choices, and the final text.

\end{document}